\documentclass[12pt]{amsart}

\usepackage{amsthm,amsmath,amsfonts,amssymb,amsthm,graphicx,epsf,epstopdf,makecell,multirow,
color,pifont,subfigure,flafter,bm,amscd}
\usepackage{amsmath}
\usepackage{amsthm}
\usepackage{amssymb}
\usepackage{amscd}
\usepackage{amsfonts}
\usepackage{amsbsy}
\usepackage{graphicx,epsf,epstopdf,color,pifont, subfigure,flafter,bm,amscd,caption,booktabs}
\usepackage{enumitem}

\newcommand{\R}{\ensuremath{\mathbb{R}}}

\newcommand{\G}{\Gamma}

\def\p{\partial}
\def\e{\varepsilon}

\def\epsilon{\varepsilon}
\allowdisplaybreaks \allowdisplaybreaks[4]

\def \a {\alpha}

\def \e {\varepsilon}

\def \p {\partial}

\newtheorem{theorem}{Theorem}[section]
\newtheorem{proposition}[theorem]{Proposition}
\newtheorem{corollary} [theorem]{Corollary}
\newtheorem{lemma}[theorem]{Lemma}

\newtheorem{remark}[theorem]{Remark}

\newtheorem{definition}[theorem]{Definition}

\begin{document}

\title[Limit cycles of algebraic foliations]
{Limit cycles of singular algebraic foliations on the real projective plane}

\author[Gan, He and Xiao]
{Tianwei Gan$^1$,   Hongjin He$^2$,  Dongmei Xiao$^{1}$$^{\dag}$}

\address{$^1$    School of Mathematical Sciences, MOE-LSC, Shanghai
	Jiao Tong University, Shanghai, 200240, P.R. China}
\email{gtw1999@sjtu.edu.cn}\email{xiaodm@sjtu.edu.cn}
\address{$^2$ College of Mathematics and Statistics, Key Laboratory of Nonlinear Analysis and its Applications (Ministry of Education),
 Chongqing University, Chongqing, 401331,  China}
\email{hehongjin000@126.com}

\begin{abstract}
   In this paper we consider cycles of singular algebraic foliations on the real projective plane. Our first result shows that there exists a limit cycle of singular algebraic foliations such that the limit cycle cannot be realized in any affine chart, that is, it fails to appear in any affine plane via the usual dehomogenization.
   On the other hand, under an appropriate convexity assumption, we verify that a cycle becomes fully visible, lying entirely in some affine chart. We then analyze the accumulation behavior of limit cycles of singular algebraic foliations on the real projective plane. Last we prove that algebraic foliations of degree at most two on the real projective plane have only finitely many limit cycles, and the maximum number of limit cycle is zero, one and at least six for algebraic foliations of degree zero, one and two on the real projective plane, respectively. This result extends  Bamon's finiteness theorem for quadratic affine vector fields in [Publ. Math. IH\'ES, (64)1986, pp. 111--142].
\end{abstract}

\subjclass[2010]{34C07, 34C08, 37F75}

\keywords{Algebraic foliations, real projective plane,  convexity, the maximum number, limit cycles}
\date{}
\dedicatory{}

\maketitle
\tableofcontents

\section{Introduction}

The subject of planar real polynomial differential equations (or vector fields $\mathcal{X}$)
\begin{equation}\label{DEs}
  \frac{dx}{dt}=X(x,y),\ \  \frac{dy}{dt}=Y(x,y)  \ \ ({\mathcal{X}}={X}(x,y)\frac{\partial}{\partial x}+{Y}(x,y)\frac{\partial}{\partial y}, \ {\rm {resp.}})
\end{equation}
is originated from Poincar\'e's seminal memoirs \cite{Poincare}. His objective was to carry out a global qualitative analysis on the whole plane $\mathbb{R}^2$ of all solutions defined by  \eqref{DEs}. Since this is a broad subject, Poincar\'e initially restricted his attention to the subclass where $X(x,y)$ and  $Y(x,y)$ are polynomials of real variables $x$ and $y$, denoted by $X,Y\in \mathbb{R}[x,y]$, the ring of all polynomials in the variables $x$ and
$y$ and coefficients in the real numbers $\mathbb{R}$.  System \eqref{DEs} (vector field $\mathcal{X}$) is called {\it the real planar polynomial differential equations of degree $n$} ({\it the $n$-degree vector field, resp.}) if the maximum degree of polynomials $X(x,y)$ and  $Y(x,y)$ is $n$, i.e. $n=\max\{\deg X, \deg Y\}$.
A subsequent major impetus came from the second part of {\it Hilbert's 16th problem} \cite{Hilbert1902}.  Hilbert formulated his 16th problem in two parts. The first part is about the relative position of {\it ovals} (nonsingular closed connected components) of  real algebraic curve of order $n$ in the real projective plane (denoted by $\mathbb{RP}^2$) when the number of ovals arrives the maximum $(n-1)(n-2)/2+1$, given by Harnack \cite{harnack}, and the corresponding investigation as to the number, form, and position of the sheets of an algebraic surface in space.  And the second part is {\it the question as to  the maximum number and position of Poincar\'e's boundary cycles (cycles limites) for a differential equation of the first order and degree of the form
$$
\frac{dy}{dx}=\frac{Y}{X},
$$
where $X$ and $Y$ are rational integral functions of the $n$-th degree in $x$ and $y$. Written homogeneously, that is
$$
X(y\frac{dz}{dt}-z\frac{dy}{dt}) + Y(z\frac{dx}{dt}-x\frac{dz}{dt}) + Z(x\frac{dy}{dt}-y\frac{dx}{dt})=0
$$
where $X$, $Y$, and $Z$ are rational integral homogeneous functions of the $n$-th degree in $x$, $y$, $z$, and the latter are to be determined as functions of the parameter $t$.}

In the second part of Hilbert's 16th problem, as originally formulated in \cite{Hilbert1902},  Hilbert proposed considering the homogenization of the vector field $\mathcal{X}$, and polynomials  $X(x,y)$ and $Y(x,y)$ are coprime in $\mathbb{R}[x,y]$, denoted by $\gcd(X,Y)=1$.

Inspired by this homogenization, we aim to study limit cycles of system \eqref{DEs} or planar polynomial vector field $\mathcal{X}$ in $\mathbb{RP}^2$, which arises as the natural compactification of the real plane $\mathbb{R}^2$ via homogenization.
This leads to an equivalent formulation in terms of {\it singular algebraic foliations $\mathcal{F}$ on $\mathbb{RP}^2$}. To avoid confusion with the homogeneous coordinates of the real projective plane $\mathbb{RP}^2$, we replace $X(x,y)$ and $Y(x,y)$ in system \eqref{DEs} or vector fields $\mathcal{X}$ by the polynomials $P(x,y)$ and $Q(x,y)$, respectively. Therefore, in affine coordinates $(x,y)\in \mathbb{R}^2$, system \eqref{DEs} and  vector fields $\mathcal{X}$ can be given by a polynomial $1$-form
\begin{equation}\label{r1form}
\omega :=P(x,y)dy - Q(x,y)dx, \ \gcd(P,Q)=1, \ n=\max\{\deg P, \deg Q\}
\end{equation}
in an equivalent way. This  defines a singular algebraic foliations $\mathcal{F}$ in $\mathbb{R}^2$, denoted by $\mathcal{F}_r$. The singular set of $\mathcal{F}_r$  is
$$
\operatorname{Sing}(\mathcal{F}_r):=\{(x,y)\in \mathbb{R}^2:\ P(x,y)=Q(x,y)=0\},
$$
and the leaves of $\mathcal{F}_r$  outside $\operatorname{Sing}(\mathcal{F}_r)$ is the solutions of the differential equation $\omega=0$.

Let $\sigma: \{(x,y,z)\in\mathbb{R}^3:z\neq 0\}\rightarrow \mathbb{R}^2$ be given by $(x,y,z)\mapsto(\frac{x}{z},\frac{y}{z})$. Then
$$
\sigma^*(\omega)=P\left(\frac{x}{z},\frac{y}{z}\right)d\left(\frac{y}{z}\right)-Q\left(\frac{x}{z},\frac{y}{z}\right)d\left(\frac{x}{z}\right)
$$
has poles at  $\{{z}=0\}$. By multiplying by a suitable factor, we define an analytic homogenized 1-form $${\omega^*}:=z^{n+2}\sigma^*({\omega}) \ (\rm{or}\  {\omega^*}:=z^{n+1}\sigma^*({\omega}))$$ if $xQ_n(x,y)-yP_n(x,y)\not\equiv0$ ($xQ_n(x,y)-yP_n(x,y)\equiv0$, resp.), where $P_n(x,y)$ and $Q_n(x,y)$ are the homogeneous parts of degree $n$ in $P(x,y)$ and $Q(x,y)$, respectively. This process is known as the Poincar\'e compactification (see \cite{DLA2006, Zhang1992}), and  $\omega^*$ defines a singular algebraic foliation $\mathcal{F}$ in $\mathbb{RP}^2$.

Analogous to the definition of singular holomorphic foliations in $\mathbb{CP}^2$ in \cite{Cerveau1991}, we write $\mathcal{F}$
in homogenous coordinates $[X : Y : Z]\in \mathbb{RP}^2$, which  is given by a real polynomial differential $1$-form
$$
\omega^* :=P^*(X,Y,Z)dX + Q^*(X,Y,Z)dY + R^*(X,Y,Z)dZ,
$$
which satisfies {\it the projective condition}
$$
XP^*(X,Y,Z) + YQ^*(X,Y,Z) + ZR^*(X,Y,Z)=0,
$$
where $P^*$, $Q^*$ and $R^*$ are homogeneous polynomials of the same degree and have no common factors, denoted by $\gcd(P^*, Q^*, R^*)=1$.
The singular set of $\mathcal{F}$ is
$$\operatorname{Sing}(\mathcal{F}):=\{[X:Y:Z]\in\mathbb{RP}^2:\ P^*=Q^*=R^*=0 \}.$$
Clearly, $\operatorname{Sing}(\mathcal{F}_r)=\operatorname{Sing}(\mathcal{F})_{Z=1}.$
The leaves of $\mathcal{F}$ are obtained by the intersection of the
leaves of singular holomorphic foliations in $\mathbb{CP}^2$ with $\mathbb{RP}^2$.

Lins Neto \cite[Section 1.4]{LinsNeto1988} characterized  the degree of singular holomorphic foliations in $\mathbb{CP}^2$ by a fundamental geometric property. This idea can be applied to  define  the degree of singular algebraic foliation $\mathcal{F}$ in $\mathbb{RP}^2$.
In algebraic terms, $\mathcal{F}$ is called {\it a degree $d$ foliation  in $\mathbb{RP}^2$} if  in an affine coordinates $(x, y)\in \mathbb{R}^2$, $\mathcal{F}_r$ is given by
$\omega = p(x, y) dy - q(x, y) dx + r(x, y)(x dy - y dx)$, where $ p, q, r\in  R[x, y]$, and either
$r(x,y)\equiv 0$ and $\max\{\deg(p), \deg(q)\} = d$ or $r(x,y)$ is a non-zero homogeneous polynomial of
degree $d$ and $\max\{\deg(p), \deg(q)\} \le d$.

A {\it limit cycle} of a singular algebraic foliation $\mathcal{F}$ in $\mathbb{RP}^2$ is defined as {\it an isolated closed leaf homeomorphic to the circle $\mathbb{S}^1$}. There are two types of limit cycles for $\mathcal{F}$ in $\mathbb{RP}^2$, see \cite{Movasati2008}. One is \textit{affine limit cycle} which is an oval in $\mathbb{RP}^2$, and the other is \textit{projective limit cycle} which is a pseudo-line in $\mathbb{RP}^2$.

The second part of Hilbert's 16th problem has been extensively studied, and yielding numerous significant results, e.g. \cite{ChristopherLi2007, Dulac1923,Ilyashenko1985,Ilyashenko1991,Kaloshin2003,Li1986,Li2003,Roussarie1998,Yeung2025,Ye1986} and references therein. To the best of our knowledge, there were few literature on limit cycles of system \eqref{DEs} in $\mathbb{RP}^2$. The systematic study of projective limit cycles was first undertaken by Movasati and Vieira in \cite{Movasati2008}, they prove that only odd degree foliations may have projective limit cycles and that foliations with exactly one
real simple singularity have a projective limit cycle. Clearly,
the projective limit cycle does not exist for system \eqref{DEs}  in $\R^2$. {\it A natural question} is if all affine limit cycles of $\mathcal{F}$ can become limit cycles of system \eqref{DEs} in $\R^2$, in the other word, whether there exist affine limit cycles of singular algebraic foliations $\mathcal{F}$ in $\mathbb{RP}^2$ that cannot be embedded into any affine chart $\R^2$.  Furthermore, we can ask what kind of constraints can ensure that affine limit cycles of $\mathcal{F}$ in $\mathbb{RP}^2$ can be embedded into  an affine chart $\mathbb{R}^2$, and how many limit cycles singular algebraic foliation $\mathcal{F}$ of degree $d$ have in $\mathbb{RP}^2$?  The objective of this article is to answer these questions and discuss the maximum number of limit cycles for $\mathcal{F}$ of degree $d$  less than two in $\mathbb{RP}^2$.   In particular, we prove that the number of limit cycles for $\mathcal{F}$ with $d\leq2$ is finite in $\mathbb{RP}^2$, and the maximum number of limit cycle is zero, one and at least six for algebraic foliations of degree zero, one and two on the real projective plane, respectively.  These results extend the finiteness theorem and the maximum number of limit cycles for quadratic vector fields in $\mathbb{R}^2$, see  \cite{Bamon1986}, \cite{Chen1979},\cite{CS1983}, \cite{Shi1980}.

The structure of this paper is as follows. In Section 2, we first review some definitions and properties of singular algebraic foliations $\mathcal{F}$ in $\mathbb{RP}^2$ by following  that of singular holomorphic foliations in $\mathbb{CP}^2$  established in \cite{LinsNeto1988}, \cite{Cerveau1991}, and \cite{Movasati2008}. Then we show that there is an affine limit cycle of $\mathcal{F}$ in $\mathbb{RP}^2$ that cannot be embedded into any affine chart $\R^2$. In Section 3, we introduce the notion of strict convexity for cycles of $\mathcal{F}$ in $\mathbb{RP}^2$, and prove that  every strictly convex cycle of an algebraic foliation $\mathcal{F}$ in $\mathbb{RP}^2$ is contained in an affine chart.  Finally, in Section 4, we investigative accumulation,  finiteness and the maximum number of limit cycles for algebraic foliations $\mathcal{F}$ in $\mathbb{RP}^2$, prove that the number of limit cycles for $\mathcal{F}$ with $d\leq2$ is finite in $\mathbb{RP}^2$,  and the maximum number of limit cycle is zero, one and at least six for algebraic foliations of degree zero, one and two on the real projective plane, respectively.

\section{Affine realizability of cycles on the real projective plane}\label{arc}
In this section, we first recall some definitions and properties of algebraic foliations in $\mathbb{RP}^2$ by
following  that of singular holomorphic foliations in $\mathbb{CP}^2$  established  in \cite{Cerveau1991}, \cite{LinsNeto1988} and \cite{Movasati2008}.
We then characterize which cycles can be contained in an affine chart, or equivalently, which cycles admit an affine realization.
Finally, we construct an affine limit cycle of an algebraic foliation that is not contained in any affine chart
$\mathbb{R}^2$.

\subsection{Algebraic foliations in $\mathbb{RP}^2$ and their affine and spherical representations}\label{sec2.1}
Following the definition of holomorphic foliations in $\mathbb{CP}^2$ in \cite{Cerveau1991}, a \emph{singular algebraic foliation $\mathcal F$ of degree $d$} in $\mathbb{RP}^2$ is represented by a homogeneous polynomial $1$-form in $\mathbb R^3\setminus\{0\}$,
\begin{equation}
\label{1form}
    \alpha=A(X,Y,Z)dX+B(X,Y,Z)dY+C(X,Y,Z)dZ,
\end{equation}
where $A,B,C$ are homogeneous coprime polynomials of degree $d+1$ (i.e. $\deg A=\deg B=\deg C=d+1$, $\gcd(A, B, C)=1$), and satisfy the projective condition
\begin{equation}
\label{projective-condition}
XA+YB+ZC=0.
\end{equation}
The singular set of $\mathcal{F}$ is
$$\operatorname{Sing}(\mathcal{F}):=\{[X:Y:Z]\in\mathbb{RP}^2:A=B=C=0\}.$$
It is known that $\mathbb{RP}^2$ is covered by the three standard affine charts
$$\begin{aligned}
    U_X=\{[X:Y:Z]\in\mathbb{RP}^2:X=1\},\\
     U_Y=\{[X:Y:Z]\in\mathbb{RP}^2:Y=1\},\\
     U_Z=\{[X:Y:Z]\in\mathbb{RP}^2:Z=1\},
\end{aligned}$$
each of which is analytically diffeomorphic to $\mathbb{R}^2$.
For example, in the affine coordinates given by
$$\varphi_Z:\mathbb R^2\longrightarrow U_Z,
\qquad (x,y)\longmapsto[x:y:1],$$
the algebraic foliation $\mathcal{F}_r$ in $\mathbb{R}^2$ is represented by
\begin{equation}
\label{affine-restriction}
\alpha|_{Z=1}=A(x,y,1)\,dx+B(x,y,1)\,dy.
\end{equation}
Except $\operatorname{Sing}(\mathcal F)$, the solutions of \eqref{affine-restriction} are the leaves of $\mathcal F$ in $U_Z$.
Due to the homogeneity of \eqref{1form} and the projective condition \eqref{projective-condition}, it is easily checked that this definition is well-defined, that is, in the intersection of two different charts, for instance, $U_{X}$ and $U_{Z}$, the equations $\alpha|_{X=1}=0$ and $\alpha|_{Z=1}=0$ determine the same leaves in $(U_{X}\cap U_{Z})\setminus \operatorname{Sing}(\mathcal{F})$.

Let $\ell:=aX+bY+cZ$, where $a^2+b^2+c^2\neq0$.
A projective line $L=\{[X:Y:Z]\in\mathbb{RP}^2:\ell=0\}$ is called \emph{invariant} under $\mathcal F$ if $(\alpha\wedge d\ell)|_{\ell=0}=0$.
For a projective line $L$, $\mathbb{RP}^2$ can be uniquely decomposed as the disjoint union of $L$ and $U$, $U=\mathbb{RP}^2\setminus L$, and $U$ is homomorphic to $\mathbb{R}^2$. This leads that $U$ can be regarded as an affine chart in $\mathbb{RP}^2$.
Then $L$ is called {\it the line at infinity of this affine chart} $U$ and denote it by  $L_\infty$.

Recall that Lins Neto \cite[Section 1.4]{LinsNeto1988} characterized the degree of an algebraic foliation in $\mathbb{CP}^2$ by a fundamental geometric property.
A point $p\in L$ is called a \emph{tangency point} of $\mathcal F$  with $L$ if either $p\in\operatorname{Sing}(\mathcal F)$,  or $p\not\in\operatorname{Sing}(\mathcal F)$,  $(\alpha\wedge d\ell)|_{\ell=0}\not=0$ and $L$ is tangent to the leaf
of $\mathcal F$ through $p$ at $p$.
The number of {tangency point} of $\mathcal F$ in  $\mathbb{CP}^2$, counted with multiplicity, is a constant $d$.
Geometrically, we call  this number {\it the degree $d$ of algebraic foliation $\mathcal F$ in $\mathbb{CP}^2$}.
However, the number of {tangency point} of $\mathcal F$ in $\mathbb{RP}^2$ is not constant usually. The following proposition gives the upper bound on the number of tangency points of $\mathcal F$ in $\mathbb{RP}^2$.
The proof of this proposition is in  Appendix~\ref{proof2}.

\begin{proposition}[\cite{LinsNeto1988}]
\label{contact}
    Let $\mathcal{F}$ be an algebraic foliation of degree $d$ in $\mathbb{RP}^2$.
    If $L$ is a projective line and it is not invariant under $\mathcal{F}$, then there are at most $d$ tangency points of $\mathcal F$ with $L$ in $\mathbb{RP}^2$.
\end{proposition}

We now give a polynomial vector field defined by an algebraic foliation of degree $d$ in the affine chart.
The complex analogue is stated in \cite{LinsNeto1988}, and we include a proof of the real version in Appendix \ref{proof2}.

\begin{proposition}[\cite{LinsNeto1988,Movasati2022, Zhang2017}]
\label{cases}
    Let  $L$ be a projective line in $\mathbb{RP}^2$. Assume that $\mathcal{F}$ is an algebraic foliation of degree $d$ in $\mathbb{RP}^2$ defined by \eqref{1form} and \eqref{projective-condition}, and consider the affine chart $U=\mathbb{RP}^2\setminus L$ and its line at infinity $L_\infty=L$. Then
    in affine coordinates $(x,y)\in U$, $\mathcal{F}$ is given by a polynomial $1$-form
    $$\omega= (p(x,y)+xr(x,y))dy-(q(x,y)+yr(x,y))dx,$$
    or, equivalently, by a polynomial vector field
    \begin{equation}
        \label{vf_cases}
        \mathcal{X}=(p(x,y)+xr(x,y))\frac{\p}{\p x}+(q(x,y)+yr(x,y))\frac{\p}{\p y},
    \end{equation}
    where $r(x,y)=-\frac{A(x,y,0)}{y}=\frac{B(x,y,0)}{x}$, which is a homogeneous  polynomial of degree $d$,  $p(x,y)=B(x,y,1)-B(x,y,0)$,  $q(x,y)=-A(x,y,1)+A(x,y,0)$, and $\max\{\deg p,\deg q\}\leq d$.
   Furthermore,
    \begin{itemize}
        \item[(i)] if $r(x,y)\equiv0$ and $\max\{\deg p,\deg q\}=d$, then $(\alpha\wedge d\ell)|_{\ell=0}=0$, that is
        the line $L_\infty$ is invariant under $\mathcal F$.
        \item[(ii)] if $r(x,y)\not\equiv0$ and $\max\{\deg p,\deg q\}\leq d$, then  $L_\infty$ is not invariant under $\mathcal F$, i.e. $(\alpha\wedge d\ell)|_{\ell=0}\not=0$.
    \end{itemize}
\end{proposition}

\begin{remark}[An equivalent characterization of an invariant line]
    Assume that $L$ is a projective line, $U$ is an affine chart with $L\cap U\neq\varnothing$, and $\mathcal X$ is the polynomial vector field giving $\mathcal F$ in $U$.
    we choose affine coordinates $(x,y)$ on $U$.
    $L\cap U$ is defined by the equation $l(x,y)=ax+by+c=0$, where $(a,b)\neq(0,0)$.
    Then $L$ is invariant under $\mathcal F$ if and only if $\mathcal X(l)|_{l=0}=0,$ or equivalently, $\mathcal X(l)=kl$ for some polynomials $k\in\mathbb R[x,y]$.
\end{remark}

\begin{definition}\label{def_cases}
      We say that algebraic foliation $\mathcal{F}$ is in {\emph{generic case}} in affine chart $\mathbb{RP}^2\setminus L$ if $L$ is not invariant under $\mathcal{F}$.
    Otherwise, that is, $L$ is invariant under $\mathcal{F}$, we say that the $\mathcal{F}$ is in {\emph{non-generic case}} in the affine chart.
\end{definition}

From Definition \ref{def_cases} and Proposition \ref{cases}, we obtain the following corollary.
\begin{corollary}
\label{non-generic}
    Let $\mathcal F$ be an algebraic foliation of degree $d$ in $\mathbb{RP}^2$, and let $L$ be a projective line. In the affine chart $\mathbb{RP}^2\setminus L$,  polynomial vector field \eqref{vf_cases}
    is generic case if  $r(x,y)\not\equiv0$ and $\max\{\deg p,\deg q\}\leq d$, and \eqref{vf_cases} is non-generic case if $r(x,y)\equiv0$ and $\max\{\deg p,\deg q\}=d$.
\end{corollary}

\begin{remark} [On the terminology generic and non-generic]
The terminologies "generic" and "non-generic" in Definition~\ref{def_cases} is relative to the
choice of the line at infinity. For each fixed line at infinity,
the foliations in the generic case form an open dense subset of
the parameter space of foliations of degree $d$.

Note that the triples $(A, B, C)$ of homogeneous polynomials of degree $d+1$ form a vector space of dimension $3\binom{d+3}{2}$,and the projective condition \eqref{projective-condition} provides $\binom{d+4}{2}$ independent linear conditions. Hence, without imposing the coprimality condition $\gcd(A,B,C)=1$ temporarily, the homogeneous polynomial $1$-forms \eqref{1form} satisfying \eqref{projective-condition}
form a vector space which is isomorphic to $\mathbb{R}^{3\binom{d+3}{2}-\binom{d+4}{2}}=\mathbb{R}^{(d+1)(d+3)}$. Since nonzero real scalar multiples of a defining $1$-form \eqref{1form}
determine the same foliation, the parameter space of algebraic foliations of degree $d$ is the open subset of $\mathbb{RP}^{(d+1)(d+3)-1}=\mathbb{RP}^{d^2+4d+2}$ under coprimality condition on $(A,B,C)$, which is denoted by $\operatorname{Fol}_d$.

Fix a line at infinity  $L_\infty=\{Z=0\}$. The invariance of $L_\infty$ is equivalent to
$Z\mid A$ and $Z\mid B$. Together with
\eqref{projective-condition}, this gives
\[
    A=ZA_1,\qquad B=ZB_1,\qquad C=-XA_1-YB_1,
\]
where $A_1$ and $B_1$ are arbitrary homogeneous polynomials of
degree $d$,  with $2\binom{d+2}{2}=(d+1)(d+2)$ independent coefficients in total. Thus, before imposing coprimality $\gcd(A,B,C)=1$, the 1-forms \eqref{1form}
satisfying this invariance condition constitute a linear
subspace of dimension $(d+1)(d+2)=d^2+3d+2$, which implies the invariance of $L_\infty$ condition defines a
projective linear subspace of
$\mathbb{RP}^{d^2+4d+2}$ of dimension $d^2+3d+1$,
and hence of codimension $d+1$.

Hence, the non-generic case is the intersection of this projective
linear subspace of codimension $d+1\geq1$ with $\operatorname{Fol}_d$. Consequently, its complement, the generic case for the fixed
line at infinity $L_\infty$, is open and dense
in the usual topology.
\end{remark}

Let $\mathbb S^2=\{(X,Y,Z)\in\mathbb R^3:X^2+Y^2+Z^2=1\}$
denote the unit sphere.
Let $\pi:\mathbb{S}^2\rightarrow \mathbb{RP}^2$ be the canonical covering map defined by $\pi(X,Y,Z)=[X:Y:Z]$, and let $i: \mathbb{S}^2 \hookrightarrow \mathbb{R}^3$ be the inclusion map.
We next describe the pull-back of an algebraic foliation in $\mathbb{RP}^2$ to $\mathbb S^2$ via $\pi$ and give an explicit vector field tangent to the pull-back foliation.

\begin{proposition}[Spherical representation]
\label{spherical-representation}
Let $\mathcal{F}$ be an algebraic foliation in $\mathbb{RP}^2$ represented by the $1$-form $\alpha$ in \eqref{1form} satisfying the projective condition \eqref{projective-condition}.
Then the pull-back foliation $\pi^*\mathcal F$ is given by $i^*\alpha$.
Except $\pi^{-1}(\operatorname{Sing}(\mathcal F))$, $\pi^*\mathcal F$ is determined by the vector field
\begin{equation}
\label{vfs2}
\begin{aligned}
V&=(A,B,C)\times(X,Y,Z)\\
 &= (ZB-YC)\frac{\partial}{\partial X}
 +(XC-ZA)\frac{\partial}{\partial Y}
 +(YA-XB)\frac{\partial}{\partial Z},
\end{aligned}
\end{equation}
where $\times$ denotes the cross product in $\mathbb R^3$ and the polynomials $A,B,C$ are evaluated at $(X,Y,Z)\in\mathbb S^2$.
Moreover, the singular set of $V$ is $\pi^{-1}(\operatorname{Sing}(\mathcal F))$.
\end{proposition}

\begin{proof}
Let $\rho:\mathbb R^3\setminus\{0\}\rightarrow\mathbb{RP}^2$ be the canonical projection map. Then $\pi=\rho\circ i$.
Since $\alpha$ satisfies the projective condition \eqref{projective-condition} and $A,B,C$ are homogeneous, it satisfies the integrability condition $\alpha\wedge d\alpha=0$ \cite{Cerveau1991}.
Consequently, except $\rho^{-1}(\operatorname{Sing}(\mathcal F))$, the distribution $\ker(\alpha)$ is a completely integrable $2$-dimensional distribution on $\mathbb R^3\setminus\{0\}$, which integrates into a 2-dimensional foliation denoted by $\widetilde{\mathcal{F}}$ (see Appendix 2 in \cite{CamachoLinsNeto1985}).

For $s=(X,Y,Z)\in\mathbb R^3\setminus\{0\}$ such that $\rho(s)\notin\operatorname{Sing}(\mathcal F)$, we consider
$$d\rho_s:T_s(\mathbb R^3\setminus\{0\})\rightarrow T_{\rho(s)}\mathbb{RP}^2,$$
 which is  the differential map of $\rho$ at $s$.
We define the tangent space to $\mathcal F$ at $\rho(s)$ by $T_{\rho(s)}\mathcal F=d\rho_s\bigl(\ker(\alpha_s)\bigr)$, where $\alpha_s: T_s(\mathbb R^3\setminus\{0\})\to\mathbb R$  denotes the covector obtained by evaluating $\alpha$ at $s$.
We claim that this definition is independent of the choice of the representative $s$.

Indeed, let
$E=X\frac{\partial}{\partial X}
+Y\frac{\partial}{\partial Y}
+Z\frac{\partial}{\partial Z}$
be the radial vector field, and let $E_s$ denote its value at $s$.
Since $\rho$ is constant along each ray through the origin, $\ker(d\rho_s)=\operatorname{span}\{E_s\}.$
Moreover, the projective condition \eqref{projective-condition} gives $\alpha_s(E_s)=XA(s)+YB(s)+ZC(s)=0,$
and hence $\ker(d\rho_s)\subset\ker(\alpha_s).$
Since $\ker(\alpha_s)$ is two-dimensional, its image under $d\rho_s$ is one-dimensional.
Finally, if $\bar s=\lambda s$ with $\lambda\neq0$ is another representative of $\rho(s)$, then the homogeneity of $A,B,C$, together with the invariance of $\rho$ under scalar multiplication, gives
$d\rho_{\bar s}\bigl(\ker(\alpha_{\bar s})\bigr)
=d\rho_s\bigl(\ker(\alpha_s)\bigr).$
Thus the  tangent line is independent on the chosen representative $s$.

We now verify that $\widetilde{\mathcal F}$ is the pull-back foliation $\rho^*\mathcal F$.
By the definition of a pull-back distribution, a vector $v\in T_s(\mathbb R^3\setminus\{0\})$ belongs to $T_s(\rho^*\mathcal F)$ if and only if $d\rho_s(v)\in T_{\rho(s)}\mathcal F=d\rho_s\bigl(\ker(\alpha_s)\bigr).$
It follows that
\begin{align*}
T_s(\rho^*\mathcal F)
&=(d\rho_s)^{-1}\!\left(d\rho_x\bigl(\ker(\alpha_s)\bigr)\right)\\
&=\ker(\alpha_s)+\ker(d\rho_s)\\
&=\ker(\alpha_s),
\end{align*}
where the last equality follows from $\ker(d\rho_s)\subset\ker(\alpha_s)$.
Hence $\widetilde{\mathcal F}=\rho^*\mathcal F$ except  $\operatorname{Sing}(\mathcal F)$.

We then prove that $\pi^*\mathcal F$ is given by $i^*\alpha$.
Let $q\in\mathbb S^2$ satisfy $\pi(q)\notin\operatorname{Sing}(\mathcal F)$, and let $v\in T_q\mathbb S^2$.
Differentiating $\pi=\rho\circ i$ gives $d\pi_q=d\rho_q\circ di_q$, where $d\pi_q$ and $di_q$ denote the differentials of $\pi$ and $i$ at $q$, respectively.
Using $\widetilde{\mathcal F}=\rho^*\mathcal F$, we obtain
\begin{align*}
v\in T_q(\pi^*\mathcal F)
&\iff d\pi_q(v)\in T_{\pi(q)}\mathcal F\\
&\iff di_q(v)\in (d\rho_q)^{-1}(T_{\pi(q)}\mathcal{F})=\ker(\alpha_q)\\
&\iff(i^*\alpha)_q(v)=0.
\end{align*}
Therefore,
$$T_q(\pi^*\mathcal F)=\ker((i^*\alpha)_q),$$
which implies that the foliation $\pi^*\mathcal F$ is given by $i^*\alpha$.

Finally, we prove that $\pi^*\mathcal F$ is determined by the vector field $V(q)$ in \eqref{vfs2} except for $\operatorname{Sing}(\mathcal F)$.
Note that $\ker((i^*\alpha)_q)=\ker(\alpha_q)\cap T_q\mathbb S^2$ for every $q=(X,Y,Z)\in\mathbb S^2$.
We divide the discussion into two cases.
For $\pi(q)\notin\operatorname{Sing}(\mathcal F)$, $(A(q),B(q),C(q))$ is nonzero.
Moreover, by the projective condition \eqref{projective-condition}, $(A(q),B(q),C(q))$ is orthogonal to $q$.
Thus $\ker(\alpha_q)$ and $T_q\mathbb S^2$ are distinct $2$-dimensional planes, and their intersection is $1$-dimensional.
By computation, the vector $V(q)$ orthogonal to both $(A(q),B(q),C(q))$ and the normal vector $q=(X,Y,Z)$ to $\mathbb S^2$.
Hence $V(q)\in\ker(\alpha_q)\cap T_q\mathbb S^2=T_q(\pi^*\mathcal F).$
At every regular point, these two vectors are nonzero and orthogonal, so $V(q)\neq0$ and therefore determines $T_q(\pi^*\mathcal F)$.

Now consider the case that $\pi(q)\in\operatorname{Sing}(\mathcal F)$.
By \eqref{vfs2}, $V(q)=0$ if and only if $(A(q),B(q),C(q))$ and $q$ are linearly dependent.
The projective condition makes $(A(q),B(q),C(q))$ and $q$ orthogonal.
Since $q\neq0$, they can be both linearly dependent and orthogonal only when $A(q)=B(q)=C(q)=0$.
By the definition of $\operatorname{Sing}(\mathcal F)$, this is equivalent to $\pi(q)\in\operatorname{Sing}(\mathcal F)$.
Thus the singular set of $V$ is precisely $\pi^{-1}(\operatorname{Sing}(\mathcal F))$.
\end{proof}

\subsection{Two types of cycles and their spherical lifts}
In this subsection, following \cite{Movasati2008}, we recall the notion of cycles for algebraic foliations in $\mathbb{RP}^2$ and distinguish between two types, that is, affine cycles and projective cycles.
We then describe their lifts to $\mathbb S^2$.

Note that in the terminology of real algebraic geometry, there are two possible types for a simple smooth closed curve $\gamma$ in $\mathbb{RP}^2$, which is homeomorphic to the circle $\mathbb{S}^1$ \cite{bochnak2013real}:
\begin{itemize}
    \item $\gamma$ is called an {\it oval}, if its complement in $\mathbb{RP}^2$ has two connected components, one of which is homeomorphic to an open disk (called the interior of $\gamma$), while the other is homeomorphic to a M\"obius band (called the exterior of $\gamma$).
    \item $\gamma$ is called a {\it pseudo-line}, if its complement in $\mathbb{RP}^2$ has one connected component, which is homeomorphic to an open disk.
\end{itemize}
Similarly, for an algebraic foliation in $\mathbb{RP}^2$, the closed leaves homeomorphic to $\mathbb{S}^1$ can be classified into two types as follows \cite{Movasati2008}, and the same classification applies to limit cycles.

\begin{definition}[cycle \cite{Movasati2008} and limit cycle \cite{Ilyashenko2008}]\label{cycle}
    Let $\gamma$ be a closed leaf, which is homeomorphic to $\mathbb{S}^1$, of an algebraic foliation $\mathcal{F}$ in $\mathbb{RP}^2$. Then
    $\gamma$ is called {an \emph{affine cycle}} (resp. {a \emph{projective cycle}}) if it is an oval (resp. a pseudo-line) in $\mathbb{RP}^2$. Moreover, if $\gamma$ is isolated, i.e. it has a neighborhood containing no closed leaf of $\mathcal F$ in $\mathbb{RP}^2$ except for $\gamma$ itself, then $\gamma$ is called
    {an \emph{affine limit cycle}} (resp. {a \emph{projective limit cycle}}) if it is an oval (resp. a pseudo-line) in $\mathbb{RP}^2$.
\end{definition}

It is worth noting that the neighborhoods of an affine cycle and a projective cycle are different in $\mathbb{RP}^2$. More precisely, the neighborhood of an affine cycle is homeomorphic to an annulus $\mathbb{S}^1\times (0,1)$, and the neighborhood of a projective cycle is homeomorphic to a M\"obius band.

We first discuss projective cycles following \cite{Movasati2008}, considering the generic and non-generic cases introduced in the preceding subsection.
Recall that any algebraic foliation $\mathcal{F}$ in $\mathbb{RP}^2$ has at most one projective cycle \cite{Movasati2008}.
Indeed, if $\gamma_0$ is a projective cycle of $\mathcal{F}$, then any other cycle $\gamma_1\neq \gamma_0$ cannot intersect $\gamma_0$.
Hence, $\gamma_1\subset (\mathbb{RP}^2\setminus \gamma_0)\cong \mathbb{R}^2$.
By the Jordan Curve Theorem again, $\gamma_1$ must separate $\mathbb{RP}^2\setminus \gamma_0$ (and thus $\mathbb{RP}^2$) into two connected components, which implies $\gamma_1$ is an affine cycle.

For a non-generic algebraic foliation $\mathcal F$ of degree $d$ in $\mathbb{RP}^2$, the line at infinity $L_\infty$ is invariant under $\mathcal F$.
Since a projective cycle meets every projective line and distinct leaves cannot intersect, the only possible projective cycle is $L_\infty$ itself. Therefore,
$L_\infty$ is a projective cycle if and only if it contains no singularities of $\mathcal F$.
A classic example is the Hamiltonian foliation given by the 1-form $d(x^{2k}+y^{2k})$, where $k$ is a positive integer.
However, a generic algebraic foliation $\mathcal F$ may have some non-trivial projective cycles that are not projective lines.
For instance, for every odd integer $k\ge3$, the Jouanolou foliation given by the vector field $\mathcal{X}_J=(y^k-x^{k+1})\frac{\partial}{\partial x}+(1-yx^k)\frac{\partial}{\partial y}$ has a projective limit cycle that is not a projective line (see \cite{Movasati2008} for details).
We also recall the following result from \cite{Movasati2008}.
\begin{theorem}[\cite{Movasati2008}]
\label{even}
    An algebraic foliation of even degree in $\mathbb{RP}^2$ cannot possess a projective cycle.
\end{theorem}

To facilitate the geometric description of these two types of cycles, we recall the following elementary result in \cite{bochnak2013real}.
For clarity, we define the antipodal map on $\mathbb S^2$ by $\mathcal A:\mathbb{S}^2\to\mathbb{S}^2,\ \mathcal A(q)=-q$.

\begin{proposition}[Spherical lifts of cycles \cite{bochnak2013real}]
\label{cycle-lifts}
Let $\gamma$ be a cycle of an algebraic foliation in $\mathbb{RP}^2$, and $\pi:\mathbb S^2\to\mathbb{RP}^2$ be the canonical covering map.
\begin{itemize}
    \item If $\gamma$ is an affine cycle, then $\pi^{-1}(\gamma)$ consists of two disjoint simple closed curves that are exchanged by the antipodal map.
    Denote one of them by $\widetilde\gamma$, then $\pi^{-1}(\gamma)=\widetilde\gamma\cup\mathcal A(\widetilde\gamma)$, and $\widetilde\gamma\cap\mathcal A(\widetilde\gamma)=\varnothing$.
    \item If $\gamma$ is a projective cycle, then $\pi^{-1}(\gamma)$ is a connected simple closed curve invariant under the antipodal map.
\end{itemize}
\end{proposition}

\subsection{A criterion for affine realizability and a counterexample}
Note that the term ``affine'' in Definition \ref{cycle} refers to the topological type of the cycle, namely an oval, and does not mean that the cycle is necessarily contained in an affine chart.
The following elementary criterion distinguishes these two notions.

\begin{proposition}[Affine realizability criterion]
\label{affine-realizability-criterion}
An affine cycle $\gamma$ is contained in an affine chart if and only if it is disjoint from some projective line $L\subset\mathbb{RP}^2$.
\end{proposition}

\begin{proof}
If $\gamma$ is contained in an affine chart, then it is disjoint from the line at infinity of that chart.
On the other hand, suppose that $\gamma$ is disjoint from a projective line $L$.
Then $\gamma$ is contained in the affine chart $\mathbb{RP}^2\setminus L$.
In this chart, $\gamma$ is a closed curve and is realized as a closed orbit of the polynomial vector field given by \eqref{vf_cases} according to Proposition \ref{cases}.
This is precisely the situation usually considered for planar polynomial vector fields in $\mathbb{R}^2$.
\end{proof}

In the study of affine limit cycles, the non-generic case is the one most directly related to the classical theory of planar polynomial vector fields in $\mathbb R^2$.
Indeed, consider a polynomial vector field ${\mathcal X}$ associated with the polynomial $1$-form \eqref{r1form}.
If it satisfies the non-degenerate condition
\begin{equation}
    \label{non-degenerate condition}
    y  P_n(x,y)-x  Q_n(x,y)\not\equiv0,
\end{equation}
where $ P_n$ and $ Q_n$ are the homogeneous parts of degree $n$ in $P(x,y)$ and $Q(x,y)$,  then its projective compactification does not make the entire line at infinity $L_\infty$ singular, and $L_\infty$ remains an invariant line of the compactified algebraic foliation in $\mathbb{RP}^2$.
Hence, by Corollary \ref{non-generic}, this foliation is given in the affine chart $\mathbb{RP}^2\setminus L_\infty$ by a vector field of the non-generic case in Definition \ref{def_cases}.
If, on the other hand, $y  P_n-x  Q_n\equiv0$, then, after removing the common factor appearing in the projective compactification, the resulting algebraic foliation in $\mathbb{RP}^2$ is given by the generic case in Definition \ref{def_cases}.

Thus, the planar polynomial vector fields satisfying the non-degenerate condition \eqref{non-degenerate condition} at infinity correspond, after projective compactification, to the non-generic case.
Moreover, no affine limit cycle can intersect the invariant line $L_\infty$.
Therefore, all affine limit cycles of such a non-generic foliation lie in the same affine chart $\mathbb{RP}^2\setminus L_\infty$, and their study reduces to that of limit cycles of a polynomial vector field in $\mathbb{R}^2$.

However, not every affine limit cycle of the generic case is contained in an affine chart.
Equivalently, not every oval in $\mathbb{RP}^2$ is disjoint from a projective line.
The following proposition provides an example.

\begin{proposition}
\label{non-affinely-realizable-cycle}
There exists an algebraic foliation in $\mathbb{RP}^2$ having an affine limit cycle that is not contained in any affine chart.
\end{proposition}

\begin{proof}
Consider an algebraic foliation $\mathcal{F}_0$ which is determined by 1-form
\begin{equation}\label{counterexample}
    \omega_0=A_0(X,Y,Z)dX+B_0(X,Y,Z)dY+C_0(X,Y,Z)dZ,
\end{equation}
where
$$
\begin{aligned}
A_0(X,Y,Z)&=(X^2+Y^2+Z^2)\frac{\partial H_0}{\partial X}-(6X-Y)H_0(X,Y,Z),\\
B_0(X,Y,Z)&=(X^2+Y^2+Z^2)\frac{\partial H_0}{\partial Y}-(6Y+X)H_0(X,Y,Z),\\
C_0(X,Y,Z)&=(X^2+Y^2+Z^2)\frac{\partial H_0}{\partial Z}-6ZH_0(X,Y,Z),
\end{aligned}
$$
here
\begin{equation}
    \label{barH_0} H_0(X,Y,Z)=\frac14(X^2+Y^2+Z^2)^3-\bigl(X^3-3XY^2-Z(X^2+Y^2+Z^2)\bigr)^2.
\end{equation}
Clearly, $\mathcal{F}_0$ satisfies the projective condition \eqref{projective-condition} since
$$
\begin{aligned}
XA_0+YB_0+ZC_0&=(X^2+Y^2+Z^2)(X\frac{\partial H_0}{\partial X}+Y\frac{\partial H_0}{\partial Y}+Z\frac{\partial H_0}{\partial Z})-6(X^2+Y^2+Z^2)H_0\\
&=-6(X^2+Y^2+Z^2)H_0+6(X^2+Y^2+Z^2)H_0=0
\end{aligned}
$$
by Euler's identity.

According to Proposition \ref{cases}, in the affine chart $U_Z$ with coordinates $(x,y)=(\frac{X}{Z},\frac{Y}{Z})$, the foliation $\mathcal F_0$ is given by a vector field
$$(P_0(x,y)+xR_0(x,y))\frac{\partial}{\partial x}+(Q_0(x,y)+yR_0(x,y))\frac{\partial}{\partial y},$$
where
$$\begin{aligned}
    R_0(x,y) =&\frac34x^6+18x^5y-\frac{27}{4}x^4y^2-60x^3y^3
+\frac{33}{4}x^2y^4+18xy^5-\frac14y^6,\\
P_0(x,y)=&\frac34x+2y-12xy+\frac54x^3+4x^2y+\frac54xy^2+4y^3-2x^4-32x^3y+6x^2y^2\\
&+\frac14x^5+14x^4y+\frac12x^3y^2-32x^2y^3
+\frac14xy^4+2y^5\\
&-2x^6-20x^5y+4x^4y^2-8x^3y^3+6x^2y^4+12xy^5,\\
    Q_0(x,y)=&-2x+\frac34y-6x^2+6y^2-4x^3+\frac54x^2y-4xy^2+\frac54y^3\\
&-4x^4-2x^3y-24x^2y^2+6xy^3+12y^4\\
&+4x^5+\frac14x^4y-28x^3y^2+\frac12x^2y^3
+16xy^4+\frac14y^5\\
&+2x^6-2x^5y-22x^4y^2+4x^3y^3-18x^2y^4+6xy^5+6y^6.
\end{aligned}$$
Note that $R_0$ is a homogeneous polynomial of degree six, and $P_0$ and $Q_0$ are polynomials of degree six.
Using the formula \eqref{vfs2}, it is calculated directly that the corresponding vector field $\mathcal{X}_0$ on $\mathbb{S}^2$ determined by $\omega_0$ is given by
$$
\begin{aligned}
\dot X&=(X^2+Y^2+Z^2)(Z\frac{\partial H_0}{\partial Y}-Y\frac{\partial H_0}{\partial Z})-XZH_0,\\
\dot Y&=(X^2+Y^2+Z^2)(X\frac{\partial H_0}{\partial Z}-Z\frac{\partial H_0}{\partial X})-YZH_0,\\
\dot Z&=(X^2+Y^2+Z^2)(Y\frac{\partial H_0}{\partial X}-X\frac{\partial H_0}{\partial Y})+(X^2+Y^2)H_0.
\end{aligned}
$$
It is easily checked that
$$
\mathcal{X}_0 H_0=H_0((X^2+Y^2)\frac{\partial H_0}{\partial Z}-Z(X\frac{\partial H_0}{\partial X}+Y\frac{\partial H_0}{\partial Y}))=H_0((X^2+Y^2+Z^2)\frac{\partial H_0}{\partial Z}-6ZH_0),
$$
i.e. $\mathcal{X}_0 H_0=K_0 H_0$, and
$$
K_0=(X^2+Y^2+Z^2)\frac{\partial H_0}{\partial Z}-6ZH_0,
$$
which means that $H_0=0$ is an invariant algebraic curve of vector field $\mathcal{X}_0$.

We claim that $ H_0=0$ is an affine cycle of $\mathcal F_0$, that is, $H_0=0$ is an oval in $\mathbb{RP}^2$ and the foliation $\mathcal{F}_0$ (or vector field $\mathcal{X}_0$) has no singularities on $H_0=0$.

We first prove that the vector field $\mathcal{X}_0$ has no singularities on $H_0=0$.

On the sphere $\mathbb{S}^2:X^2+Y^2+Z^2=1$ and the algebraic curve $H_0=0$, we have
$$
\begin{aligned}
X\dot Y-Y \dot X&=(X^2+Y^2)\frac{\partial H_0}{\partial Z}-Z(X\frac{\partial H_0}{\partial X}+Y\frac{\partial H_0}{\partial Y})\\
&=(X^2+Y^2+Z^2)\frac{\partial H_0}{\partial Z}-Z(X\frac{\partial H_0}{\partial X}+Y\frac{\partial H_0}{\partial Y}+Z\frac{\partial H_0}{\partial Z})\\
&=(X^2+Y^2+Z^2)\frac{\partial H_0}{\partial Z}-6ZH_0=K_0.
\end{aligned}
$$
Let $g_0(X,Y,Z)=X^3-3XY^2-Z$. Since $ H_0|_{\mathbb S^2}=\frac14-g_0^2$ and $g_0(-q)=-g_0(q)$ for any $q\in\mathbb{S}^2$,
 $H_0=0$ has two components $ \pm\tilde\gamma: g_0=X^3-3XY^2-Z=\pm\frac{1}{2}$ on $\mathbb{S}^2$. On curve $\tilde\gamma$, we have
\begin{equation}\label{K0}
K_0|_{X^2+Y^2+Z^2=1,g_0=\frac{1}{2}}=\frac{\partial H_0}{\partial Z}|_{X^2+Y^2+Z^2=1,g_0=\frac{1}{2}}=2(Z+\frac{3}{8})^2+\frac{23}{32}>0.
\end{equation}
Similarly, $K_0<0$ on curve $-\tilde\gamma$. Hence, $X\dot Y-Y \dot X\neq 0$ on $H_0=0$, which means that $\dot X$ and $\dot Y$ cannot vanish simultaneously on $H_0=0$. Thus the vector field $\mathcal{X}_0$ has no singularities on $H_0=0$.

Next we prove $H_0=0$ is an affine cycle in $\mathbb{RP}^2$. Let $\pi:\mathbb S^2\to\mathbb{RP}^2$ be the canonical covering map. The inverse image of the projective curve \eqref{barH_0} in $\mathbb S^2$ is $\widetilde\gamma\cup(-\widetilde\gamma)$, where $\widetilde\gamma=\{(X,Y,Z)\in\mathbb S^2:g_0(X,Y,Z)=\frac12\}$.
To describe the level set $\widetilde\gamma$, we write $$X=\sqrt{1-Z^2}\cos\theta, \  Y=\sqrt{1-Z^2}\sin\theta$$ due to $(X,Y,Z)\in\mathbb S^2$.
Then the equation $g_0(X,Y,Z)=\frac12$ becomes
\begin{equation}
    \label{phiZ}
    \cos(3\theta)=\phi(Z):=\frac{Z+\frac12}{(1-Z^2)^{\frac32}},
\end{equation}
and $\phi'(Z)=\frac{2Z^2+\frac32Z+1}{(1-Z^2)^{5/2}}>0.$
Thus, the admissible values of $Z$ form an interval $[Z_-,Z_+]$, where $Z_-<0<Z_+$ and $\phi(Z_-)=-1$, $\phi(Z_+)=1$.
For each $Z\in(Z_-,Z_+)$, the equation \eqref{phiZ} has the six solutions $\theta=\frac{2k\pi}{3}\pm\frac13\arccos(\phi(Z)),\ k=0,1,2.$
And so, $\widetilde\gamma$ is composed of six arcs as $Z$ varies.
When $Z\to Z_+$, we have $\phi(Z)\to1$ and hence $\arccos(\phi(Z))\to0$.
Thus the two arcs corresponding to the same value of $k$ meet at $\theta=\frac{2k\pi}{3}.$
When $Z\to Z_-$, we have $\phi(Z)\to-1$ and hence $\arccos(\phi(Z))\to\pi$.
In this case, $\frac{2k\pi}{3}+\frac{\pi}{3}=\frac{2(k+1)\pi}{3}-\frac{\pi}{3},$ so the one arc corresponding to $k$ meets the other arc corresponding to $k+1$, where $k+1$ is taken modulo $3$.
These identifications join all six arcs into a single closed curve, that is, $\widetilde\gamma$.
Since the projective algebraic curve is invariant under $\mathcal F_0$ and contains no singularity of $\mathcal F_0$, Proposition \ref{spherical-representation} implies that $\widetilde\gamma$ contains no singularity of the pull-back foliation and is everywhere tangent to it.
As $\widetilde\gamma$ is compact and connected, it is a closed regular leaf and hence a smooth simple closed curve.
Since $\widetilde\gamma$ and $-\widetilde\gamma$ are disjoint,  $\pi|_{\widetilde\gamma}$ is injective.
Therefore, $\gamma_0:=\pi(\widetilde\gamma)=\pi(-\widetilde\gamma)$ is the only connected component of $H_0=0$ and is a cycle of $\mathcal F_0$.
Since $\pi^{-1}(\gamma_0)=\widetilde\gamma\cup(-\widetilde\gamma)$ has two connected components, Proposition \ref{cycle-lifts} implies that $\gamma_0$ is an affine cycle in $\mathbb{RP}^2$.

We claim that $H_0=0$ is an affine limit cycle of the foliation $\mathcal{F}_0$. In fact, we only need to prove that $\widetilde\gamma$ is a limit cycle of the vector field $\mathcal{X}_0$ on $\mathbb{S}^2$. Fix a point
$$
q_0=(\sqrt{1-Z_+^2},0,Z_+)\in\widetilde\gamma
$$
and a set
$$
\Sigma=\{(\sqrt{1-Z^2},0,Z):|Z-Z_+|<\delta\},
$$
where $\delta>0$ is sufficiently small. Euler's identity therefore gives
$$
\sqrt{1-Z_+^2}\frac{\partial H_0}{\partial X}(q_0)
+Z_+\frac{\partial H_0}{\partial Z}(q_0)
=6H_0(q_0)=0.
$$
Hence
$$
\frac{\partial H_0}{\partial X}(q_0)
=-\frac{Z_+}{\sqrt{1-Z_+^2}}
  \frac{\partial H_0}{\partial Z}(q_0).
$$
It follows that
$$
\begin{aligned}
\left.\frac{d}{dZ}H_0(\sqrt{1-Z^2},0,Z)\right|_{Z=Z_+}
&=-\frac{Z_+}{\sqrt{1-Z_+^2}}
  \frac{\partial H_0}{\partial X}(q_0)
  +\frac{\partial H_0}{\partial Z}(q_0)\\
&=\left(\frac{Z_+^2}{1-Z_+^2}+1\right)
  \frac{\partial H_0}{\partial Z}(q_0)\\
&=\frac{1}{1-Z_+^2}
  \frac{\partial H_0}{\partial Z}(q_0)>0,
\end{aligned}
$$
where the positivity follows from \eqref{K0}. Moreover, we can calculate similarly that
$$
\begin{aligned}
\dot Y(q_0)
&=\sqrt{1-Z_+^2}\frac{\partial H_0}{\partial Z}(q_0)
  -Z_+\frac{\partial H_0}{\partial X}(q_0)\\
&=\left(\sqrt{1-Z_+^2}
  +\frac{Z_+^2}{\sqrt{1-Z_+^2}}\right)
  \frac{\partial H_0}{\partial Z}(q_0)\\
&=\frac{1}{\sqrt{1-Z_+^2}}
  \frac{\partial H_0}{\partial Z}(q_0)>0.
\end{aligned}
$$
Thus $\Sigma$ is a transverse section of vector field $\mathcal{X}_0$, and $c=H_0|_\Sigma$ is a local coordinate on it.
Let $q(c)\in\Sigma$ be the unique point satisfying $H_0(q(c))=c$, and denote by $q(t;c)$ the solution of $\mathcal{X}_0$ starting at $q(c)$.
Let $T>0$ be the period of $\widetilde\gamma$ and $\tau(c)$ the local return time of solution $q(t;c)$ to $\Sigma$. It is clear that $\tau(0)=T$.
Transversality and smooth dependence on initial conditions imply that $\tau(c)$ is smooth near $c=0$.
Since $\mathcal{X}_0H_0=K_0H_0$, we have
$$
\frac{d }{d t}H_0(q(t;c))=(\mathcal{X}_0H_0)(q(t;c))=K_0(q(t;c))H_0(q(t;c)).
$$
Then the integration along the solution gives the Poincar\'e return map in the coordinate $c$ as
$$
P(c)=H_0(q(\tau(c);c))
=c\exp\left(\int_0^{\tau(c)}K_0(q(t;c))\,dt\right).
$$
Consequently, using $P(0)=0$, $\tau(0)=T$ and \eqref{K0}, we obtain
$$
P'(0)=\lim_{c\to0}\frac{P(c)}{c}
=\exp\left(\int_0^T K_0(q(t;0))\,dt\right)>1.
$$
For this analytic vector field, a nonisolated cycle would have the identity as its local Poincar\'e return map, which would imply $P'(0)=1$.
Therefore $\widetilde\gamma$ is isolated, and $P'(0)>1$ also proves that it is hyperbolic. Hence $H_0=0$ is a limit cycle of foliation $\mathcal F_0$.

Finally, we prove that $\gamma_0$ can not be contained in any affine chart.
Let
$$
L=\{[X:Y:Z]\in\mathbb{RP}^2:aX+bY+cZ=0,(a,b,c)\neq(0,0,0)\}
$$
be an arbitrary projective line, and write $\ell(X,Y,Z)=aX+bY+cZ$.
By the analysis in the preceding paragraph, the curve $\widetilde\gamma$ on $\mathbb S^2$ contains six points
$$
q_{+,k}=(\sqrt{1-Z_+^2}\cos\frac{2k\pi}{3},\,\sqrt{1-Z_+^2}\sin\frac{2k\pi}{3},\,Z_+)
$$
and
$$
q_{-,k}=(\sqrt{1-Z_-^2}\cos\frac{(2k+1)\pi}{3},\,\sqrt{1-Z_-^2}\sin\frac{(2k+1)\pi}{3},\,Z _-)
$$
for $k=0,1,2.$

\begin{figure}[h]
    \centering
    \includegraphics[width=0.3\linewidth]{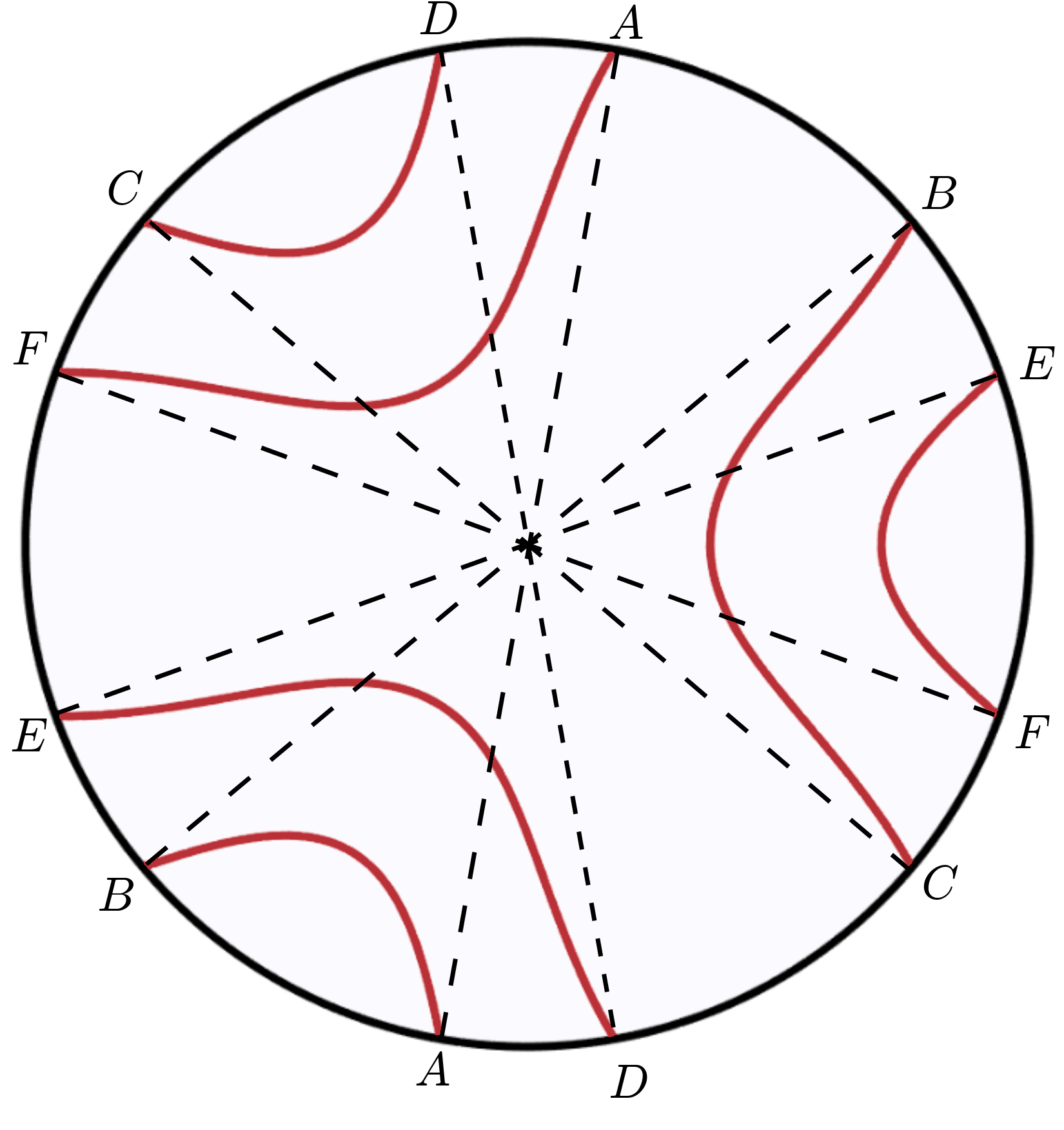}
    \caption{An affine cycle defined in \eqref{barH_0} that is not contained in any affine chart. The black circle represents the Poincar\'e disk, whose antipodal boundary points are identified to obtain $\mathbb{RP}^2$. Under this identification, the red curve represents the affine limit cycle, oriented in the cyclic order $A\to B\to C\to D\to E\to F\to A$.}
    \label{fig_example}
\end{figure}

Using the identities $\sum_{k=0}^2\cos\frac{2k\pi}{3}=\sum_{k=0}^2\sin\frac{2k\pi}{3}=0$ and the analogous identities for the angles $\frac{(2k+1)\pi}3$, we obtain
$\sum_{k=0}^2\ell(q_{+,k})=3cZ_+$, and $\sum_{k=0}^2\ell(q_{-,k})=3cZ_-.$
If $c\neq0$, these two sums have opposite signs because $Z_-<0<Z_+$.
Hence $\ell$ takes values of opposite signs on $\widetilde\gamma$.
If $c=0$, each sum is zero.
Therefore, either $\ell$ vanishes at one of these six points, or it again takes values of opposite signs.
Since $\widetilde\gamma$ is connected, the intermediate value theorem gives a point $q\in\widetilde\gamma$ such that $\ell(q)=0$.
Thus $\pi(q)\in\gamma_0\cap L$.
Therefore, the projective curve $ H_0=0$ intersects every projective line, and Proposition \ref{affine-realizability-criterion} implies that $\gamma_0$ is not contained in any affine chart (see Figure \ref{fig_example}).
\end{proof}

\section{Strict convexity and affine visibility}
\label{strict-convexity}

In this section, we introduce the notion of strict convexity for cycles in $\mathbb{RP}^2$ and investigate its geometric properties.
We first relate strict convexity in $\mathbb{RP}^2$ to convexity on $\mathbb{S}^2$.
We then show that every strictly convex cycle of an algebraic foliation in $\mathbb{RP}^2$ is contained in an affine chart, and hence is fully visible in some affine chart.
Finally, we prove that every affine cycle of an algebraic foliation of degree at most two in $\mathbb{RP}^2$ is strictly convex.
These results will play a crucial role in the proof of the finiteness of limit cycles in Section \ref{proof_main}.
The auxiliary results concerning convex subsets of $\mathbb{S}^2$ and Hausdorff convergence are collected in Appendices \ref{convexity-s2} and \ref{hausdorff-convergence}, respectively.
Here, the concept of convex subsets of $\mathbb{S}^2$ follows the definition of weak convexity in \cite{Danzer1963}.

\subsection{Projective strict convexity}

We begin with a notion of strict convexity that is relative to the projective plane $\mathbb{RP}^2$.

\begin{definition}
\label{defconvex}
A cycle in $\mathbb{RP}^2$ is called \emph{strictly convex} if every projective line intersects it in at most two distinct points.
\end{definition}

The following proposition shows that strict convexity can occur only for affine cycles.

\begin{proposition}
\label{non_projective}
A projective cycle in $\mathbb{RP}^2$ cannot be strictly convex.
\end{proposition}

\begin{proof}
Since a projective cycle is homotopic to a projective line in $\mathbb{RP}^2$, its intersection number modulo two with any projective line is one (see \cite[Section 2.4]{GP1974}).
This implies that the number of transversal intersections of a projective cycle with a projective line is odd.
If a projective cycle is strictly convex, this number can only be one.
However, a projective cycle that is not itself a projective line intersects some projective line transversally in at least three points.
If it is itself a projective line, it intersects that line in infinitely many points, which also
violates Definition \ref{defconvex}.
Thus, a projective cycle cannot be strictly convex.
\end{proof}

We next relate Definition \ref{defconvex} to the notion of convexity in $\mathbb S^2$.

Let $\mathbb{S}^2=\{x\in\mathbb{R}^3:\|x\|=1\}$ be the unit sphere, endowed with its spherical distance $d_{\mathbb{S}^2}$. We use the following notion of convexity, following the definition of weak convexity in \cite{Danzer1963}.

\begin{definition}[{\cite[Section 9.1]{Danzer1963}}]
\label{defconvexset}
A subset $S\subset\mathbb{S}^2$ is called \emph{convex} if for any two non-antipodal points $p,q\in S$, the minor arc connecting $p$ and $q$ is contained in $S$, and if $p$ and $q$ are antipodal, it is required that one of the semicircles connecting them be contained in $S$.
\end{definition}

In particular, the intersection of two convex subsets of $\mathbb{S}^2$ that
contain no pair of antipodal points is convex.

\begin{proposition}[Relationship with convexity on $\mathbb{S}^2$]
\label{relationship}
Let $\gamma$ be a strictly convex affine cycle, and let $\tilde{\gamma}\subset\mathbb{S}^2$ be a connected component of $\pi^{-1}(\gamma)$, where $\pi:\mathbb{S}^2\to\mathbb{RP}^2$ is the canonical covering map.
Then the closed region $K\subset\mathbb{S}^2$ bounded by $\tilde{\gamma}$ that contains no pair of antipodal points is convex on $\mathbb{S}^2$.
\end{proposition}

\begin{proof}
Let $\gamma$ be a strictly convex affine cycle, and let $\tilde{\gamma}\subset\mathbb{S}^2$ be a connected component of $\pi^{-1}(\gamma)$.
By the Jordan Curve Theorem on $\mathbb{S}^2$, $\tilde{\gamma}$ bounds two non-empty regions.
Let $K$ be the one that contains no pair of antipodal points.
Such a region exists because $\gamma$ is an affine cycle and its lift consists of two disjoint antipodal connected components.

Suppose, for the sake of contradiction, that $K$ is not convex.
By Definition \ref{defconvexset}, there exist points $p, q \in \operatorname{int}K$ such that the minor arc $[p,q]$ of a great circle $\tilde{L}$ passing through them is not entirely contained in $K$.
Thus, the minor arc must exit and re-enter $K$, intersecting the boundary $\tilde{\gamma}$ at least twice.
Furthermore, $\tilde{L}\setminus(p,q)$ is a major arc of length greater than $\pi$, and passes through the antipode $-p$.
Since $K\cap(-K)=\varnothing$, we have $-p\notin K$.
Because this major arc connects $p,q\in K$ and passes through a point outside $K$, it must also intersect the boundary $\tilde{\gamma}$ at least twice.
Consequently, $\tilde{L}$ intersects $\tilde{\gamma}$ in at least four distinct points.
Since $\tilde{\gamma}$ contains no pair of antipodal points, the restriction $\pi|_{\tilde{\gamma}}$ is injective.
Hence the projective line $L=\pi(\tilde{L})$ intersects $\gamma$ in at least four distinct points,
contradicting the strict convexity of $\gamma$.
\end{proof}

\subsection{Affine realizability under strict convexity}

Recall that by Proposition \ref{affine-realizability-criterion}, a subset of $\mathbb{RP}^2$ is contained in an affine chart if and only if it is disjoint from some projective line.
Although an affine cycle need not be contained in a single affine chart, the following theorem shows that strict convexity guarantees affine realizability.

\begin{theorem}
\label{embed_cycle}
Any strictly convex cycle of an algebraic foliation in $\mathbb{RP}^2$ is contained in an affine chart.
\end{theorem}

\begin{proof}
We prove this by contradiction.
Suppose that $\gamma$ is a strictly convex cycle of an algebraic foliation in $\mathbb{RP}^2$.
By Proposition \ref{non_projective}, $\gamma$ is an affine cycle.
Consider its inverse image under the canonical covering map $\pi:\mathbb{S}^2\to\mathbb{RP}^2$, and let $\tilde{\gamma}$ be a connected component of $\pi^{-1}(\gamma)$.
Then $\tilde{\gamma}$ contains no pair of antipodal points by Proposition \ref{cycle-lifts}.

Assume that $\gamma$ cannot be fully contained in any affine chart.
By Proposition \ref{affine-realizability-criterion}, a subset of $\mathbb{RP}^2$ is not contained in any affine chart if and only if it intersects every projective line.
It follows that $\tilde{\gamma}$ intersects every great circle on $\mathbb{S}^2$.
Therefore, by the Strong Separation Theorem \cite[Section 3.6]{LL2016}, the origin $O$ of $\mathbb{R}^3$ belongs to the convex hull of $\tilde{\gamma}$, that is,
$$O\in\operatorname{conv}(\tilde{\gamma}) =\left\{\sum_{i=1}^{n}k_iq_i:
k_i\geq0,\ \sum_{i=1}^{n}k_i=1,\ q_i\in\tilde{\gamma}\right\}.$$

By Carath\'eodory's Theorem in $\mathbb{R}^3$ \cite[Theorem 3.3.10]{LL2016}, the origin can be written as a convex combination of at most four points of $\tilde{\gamma}$.
Then we choose such a representation with the minimum number of points. Clearly,
a representation with one point is impossible since $O\notin\tilde \gamma$.
If the representation  contains two points, they must be antipodal, which is impossible because $\tilde{\gamma}$ contains no pair of antipodal points.
If it contains three points, these points lie in a plane $\Pi$ passing through the origin, and the great circle $\tilde L_0=\Pi\cap\mathbb S^2$ intersects $\tilde{\gamma}$ in at least three points. Thus, $\gamma$ intersects a projective line $L_0=\pi(\tilde L_0)$ in at least three points, contradicting strict convexity.
Therefore, the remaining case is that there exist exactly four points $P_1, P_2, P_3, P_4\in\tilde{\gamma}$ and positive coefficients $k_1,k_2,k_3,k_4$ such that
\begin{equation}
    \label{convex-combination}
    k_1P_1+k_2P_2+k_3P_3+k_4P_4=0,\qquad \sum_{i=1}^{4}k_i=1.
\end{equation}
We may label the points so that $P_1$ and $P_2$ are adjacent among these four
points along $\tilde{\gamma}$.

\begin{figure}[h]
    \centering
    \includegraphics[width=0.4\linewidth]{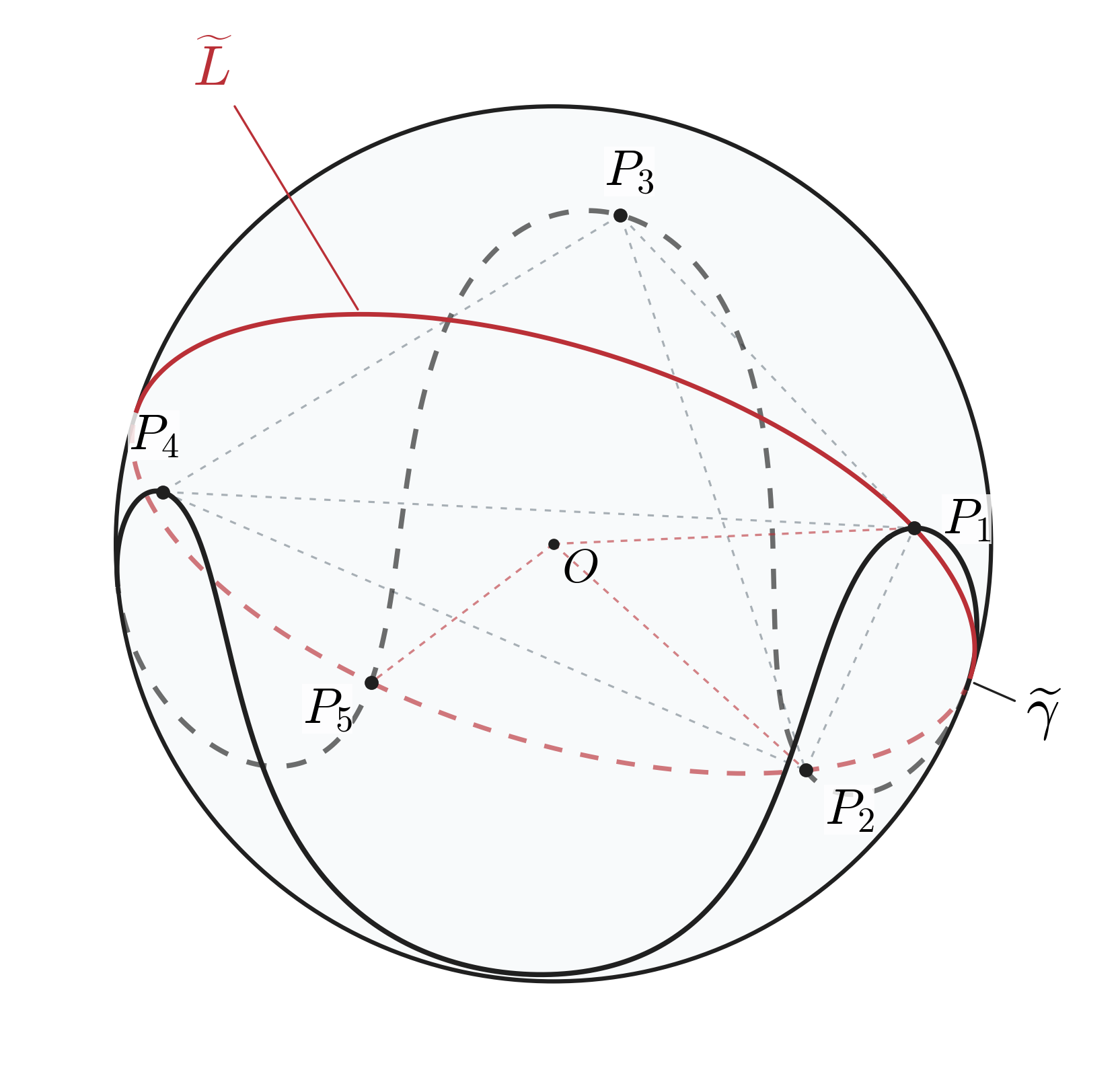}
    \caption{The intersection of the cycle $\tilde{\gamma}$ with the great
    circle $\tilde{L}$.}
    \label{fig_embed}
\end{figure}

Consider the plane in $\mathbb{R}^3$ passing through $O,P_1$, and $P_2$, see
Figure \ref{fig_embed}.
Its intersection with $\mathbb{S}^2$ is a great circle, denoted by  $\tilde{L}$.
Let $f(P)=0$ be a linear equation defining this plane.
By construction, $f(P_1)=f(P_2)=0$.
The convex combination \eqref{convex-combination} gives
\begin{equation}
    \label{p3p4}
    k_3f(P_3)+k_4f(P_4) =-k_1f(P_1)-k_2f(P_2)=0.
\end{equation}
If $f(P_3)=0$ or $f(P_4)=0$, then $\tilde{L}$ already intersects $\tilde{\gamma}$ in at least three distinct points.
Hence we suppose that $f(P_3)\neq0$ and $f(P_4)\neq0$.
Since $k_3,k_4>0$, by \eqref{p3p4}, the values $f(P_3)$ and $f(P_4)$ have opposite signs.
Let $\varphi(s)$, $s\in[s_3,s_4]$, be a parameterization of the arc of
$\tilde{\gamma}$ connecting $P_3$ and $P_4$ that does not contain $P_1$ and
$P_2$.
Then the function $g(s)=f(\varphi(s))$ is continuous.
Moreover, $g(s_3)=f(P_3)$ and $g(s_4)=f(P_4)$ have opposite signs.
By the Intermediate Value Theorem, there exists $s_*\in(s_3,s_4)$ such that $g(s_*)=0$.
Let $P_5=\varphi(s_*)$. Then $P_5$ lies on $\tilde{\gamma}$ and, by construction, also lies on the great circle $\tilde{L}$ since $f(P_5)=0$.
Thus, the great circle $\tilde{L}$ intersects $\tilde{\gamma}$ in at least three distinct points, $P_1$, $P_2$, and $P_5$, no two of which are antipodal.
Projecting this configuration to $\mathbb{RP}^2$, the projective line $L=\pi(\tilde{L})$ intersects $\gamma$ in at least three distinct points, namely, $\pi(P_1)$, $\pi(P_2)$, and $\pi(P_5)$.
This contradicts the strict convexity of $\gamma$.
Therefore, $\gamma$ is contained in an affine chart.
\end{proof}

\subsection{Cycles of foliations of degree at most two}

We now apply the preceding geometric results to algebraic
foliations of low degree in $\mathbb{RP}^2$. The following theorem extends the corresponding
convexity result for  algebraic foliations in $\mathbb{R}^2$ \cite{Ye1986}.

\begin{theorem}
\label{prop_convex}
Any affine cycle of an algebraic foliation of degree at most two in $\mathbb{RP}^2$ is strictly convex.
\end{theorem}

\begin{proof}
We provide the proof for the degree-two case, as the cases of degrees zero and one are analogous.
Let $\mathcal{F}$ be an algebraic foliation of degree two in $\mathbb{RP}^2$.
Consider the pull-back of $\mathcal{F}$ to the sphere $\mathbb{S}^2$ via the canonical covering map $\pi:\mathbb{S}^2\to\mathbb{RP}^2$, where $\pi(x,y,z)=[x:y:z]$.
Let $\gamma$ be an affine cycle of $\mathcal{F}$ in $\mathbb{RP}^2$. By Proposition \ref{cycle-lifts}, we can set $\pi^{-1}(\gamma)=\tilde{\gamma}\cup(-\tilde{\gamma}),$ where $\tilde{\gamma}\cap(-\tilde{\gamma})=\varnothing.$

We proceed by contradiction.
If $\gamma$ is not strictly convex, there exists a projective line intersecting it in more than two points.
By a sufficiently small perturbation, we may choose such a line $L$ so that it intersects $\gamma$ transversally.
Since the complement of an affine cycle has two components, the line must leave the interior of the affine cycle as often as it enters.
Hence the number of transversal intersections between $L$ and $\gamma$ is even, and therefore $L$ intersects $\gamma$ in at least four points.
The inverse image $\tilde{L}=\pi^{-1}(L)$ is a great circle on $\mathbb{S}^2$, and $\tilde{L}$ consequently intersects $\tilde{\gamma}$ in at least four points.
After applying a rotation of $\mathbb S^2$, we may assume that $\widetilde L$ is the equator.
Let $K$ be the closed region bounded by $\widetilde\gamma$ that contains no pair of antipodal points, so that $K\cap(-K)=\varnothing$.
Since $\widetilde L$ intersects $\widetilde\gamma$ transversally, the open arcs between consecutive intersection points along $\widetilde L$ alternate between the interior and the exterior of $K$.
There are at least two exterior arcs, and at most one can have length at least $\pi$.
Choose an exterior arc of length less than $\pi$ and denote its endpoints by $p_2,p_3$ in the chosen direction along $\widetilde L$.
Let $p_1$ and $p_4$ be the immediately preceding and following intersection points, respectively, in the same direction along $\widetilde L$.
Then $(p_1,p_2)\cup(p_3,p_4)\subset\operatorname{int}(K)$
and $(p_2,p_3)\subset\mathbb S^2\setminus K$.
The two arcs contained in $K$ also have length less than $\pi$, since $K$ contains no pair of antipodal points.
Thus these three consecutive arcs are minor arcs, denoted by $[p_1,p_2]$, $[p_2,p_3]$, and $[p_3,p_4]$, with $[p_1,p_2],[p_3,p_4]\subset K$, see Figure \ref{fig_convex}.

\begin{figure}[htbp]
    \centering
    \includegraphics[width=0.95\linewidth]{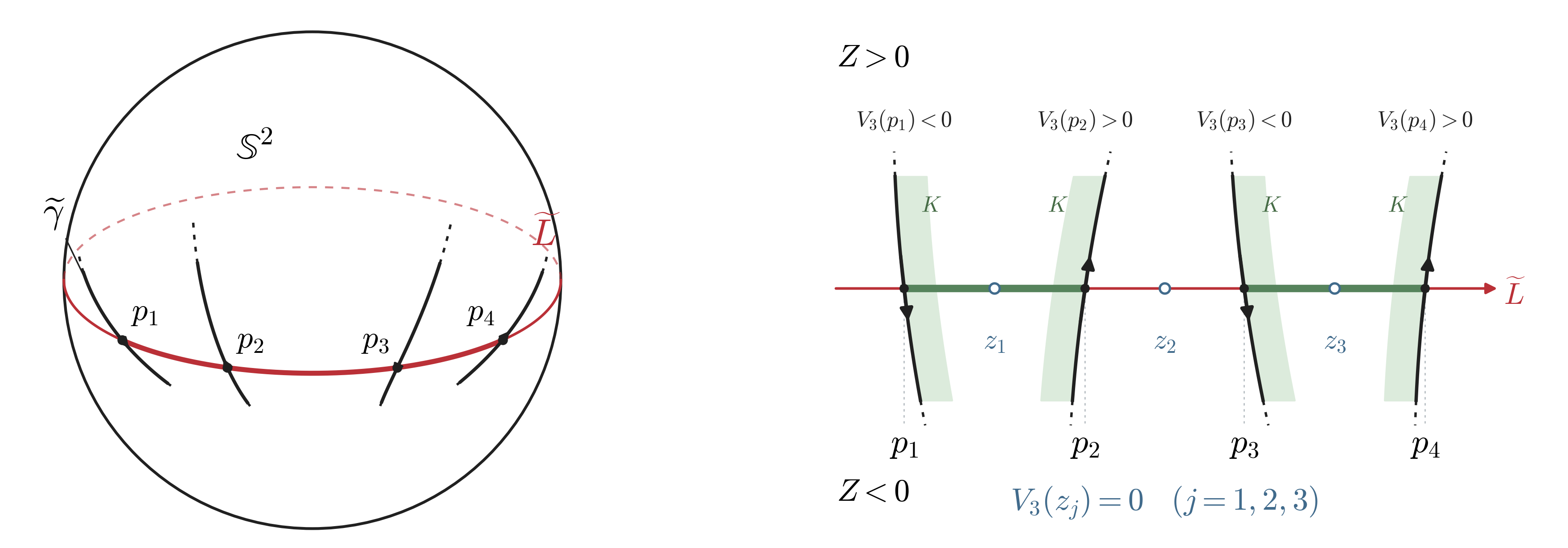}
    \caption{The left panel shows four consecutive intersections along $\widetilde L$, where only local portions of $\widetilde\gamma$ are shown.
    In the right panel, the shaded strips indicate the side belonging to $K$, and the green intervals lie in $\operatorname{int}(K)$.
    The points $z_j$ satisfy $V_3(z_j)=0$, $j=1,2,3$.}
    \label{fig_convex}
\end{figure}

We claim that the arc connecting $p_4$ back to $p_1$ along $\widetilde L$ must be a major arc.
Its length cannot equal $\pi$, since $p_1$ and $p_4$ are not antipodal.
Suppose, for the sake of contradiction, that it is a minor arc, denoted by $[p_4,p_1]$.
Since the lengths of both $[p_2,p_3]$ and $[p_4,p_1]$ are strictly less than $\pi$, we have $p_4\in[-p_1,-p_2]\cup[-p_2,p_1],$ and $p_3\in[p_2,-p_1]\cup[-p_1,-p_2]$.
If either $p_3$ or $p_4$ belongs to $[-p_1,-p_2]$, then $[-p_1,-p_2]\cap[p_3,p_4]\neq\varnothing$.
Otherwise, $p_3\in[p_2,-p_1)$ and $p_4\in(-p_2,p_1]$.
The cyclic order is therefore $p_3,-p_1,-p_2,p_4$, so the arc $[p_3,p_4]$ contains $[-p_1,-p_2]$.
Hence, in either case, $[-p_1,-p_2]\cap[p_3,p_4]\neq\varnothing$.
However, since $[p_1,p_2]\subset K$, its antipodal arc
satisfies $[-p_1,-p_2]\subset-K$.
On the other hand, $[p_3,p_4]\subset K$.
Thus $(-K)\cap K\neq\varnothing$, a contradiction.
Therefore, the remaining arc from $p_4$ to $p_1$ is a major arc.
Consequently, $[p_1,p_2]\cup[p_2,p_3]\cup[p_3,p_4]$ has total length less than $\pi$.

Let $V$ be the vector field on $\mathbb S^2$ given by Proposition \ref{spherical-representation}, and rewrite $V=V_1\frac{\partial}{\partial X} +V_2\frac{\partial}{\partial Y}
 +V_3\frac{\partial}{\partial Z}$ for short.
We now show that $V$ crosses $\widetilde L$ in alternating directions at $p_1,p_2,p_3,p_4$.
By Proposition \ref{spherical-representation}, $V$ is continuous, nonzero, and tangent to the connected closed curve $\widetilde\gamma$.
And so, it defines an orientation of $\widetilde\gamma$.
Therefore, when following $V$, the region $K$ remains on the same side of $\widetilde\gamma$.
Moreover, since $(p_1,p_2),(p_3,p_4)\subset\operatorname{int}(K)$, $(p_2,p_3)\subset\mathbb S^2\setminus K$, and taking into account transversality, it imply that $\widetilde L$ enters $K$ at $p_1,p_3$ and leaves $K$ at $p_2,p_4$.
Thus, the vector field $V$ must cross $\widetilde L$ in one direction at $p_1,p_3$ and in the opposite direction at $p_2,p_4$ to keep $K$ on the same side when following $V$, see Figure \ref{fig_convex} for instance.
Since $\widetilde L=\{Z=0\}\cap\mathbb S^2$ and all four intersections are transverse, this gives
$V_3(p_j)V_3(p_{j+1})<0,\ j=1,2,3$.
By continuity, each open arc $(p_j,p_{j+1})$ contains a point $z_j$ such that $V_3(z_j)=0$.
Since $\widetilde L=\{Z=0\}\cap\mathbb S^2$, $\pi(z_j)$ is either a tangency point of $\mathcal F$ with $L$ or a singularity of $\mathcal F$.
Since the arc from $p_1$ to $p_4$ has length less than $\pi$, it contains no pair of antipodal points.
Thus $\pi$ is injective on this arc, and $\pi(z_1),\pi(z_2),\pi(z_3)$ are three distinct tangent
points of $\mathcal F$ with $L$ by Proposition \ref{contact}.
Moreover, $L$ is not invariant, since it intersects $\gamma$ transversally.
This contradicts Proposition \ref{contact}, which states that a non-invariant projective line has at most two tangent points with a degree-two foliation.
Therefore, $\gamma$ is strictly convex.
\end{proof}

\begin{remark}
    Note that Theorem \ref{prop_convex} can not be extended to algebraic foliations of degree greater than three in $\mathbb{RP}^2$.
    For instance, \cite{DesrochesJeffrey2011} gives an example of a cubic polynomial vector field on the affine plane with a non-convex limit cycle, and so, the projective compactification defines an algebraic foliation in $\mathbb{RP}^2$ with an affine cycle that is not strictly convex.
    Also, in Proposition \ref{non-affinely-realizable-cycle}, there exists an algebraic foliation $\mathcal F_0$ of degree six, given by \eqref{counterexample}, with an affine limit cycle $\{[X:Y:Z]\in\mathbb {RP}^2: H_0(X,Y,Z)=0\}$ defined in \eqref{barH_0} that is not strictly convex.
    Indeed, in the proof of Proposition \ref{non-affinely-realizable-cycle}, we show that the line $\{[X:Y:Z]\in\mathbb{RP}^2:Z=0\}$ intersects the affine cycle $H_0=0$ in six distinct points, so the affine cycle $ H_0=0$ is not strictly convex.
\end{remark}

Combining Theorems \ref{embed_cycle} and \ref{prop_convex}, we immediately
obtain the affine visibility of all affine cycles of algebraic foliations with  degree at most two.

\begin{corollary}
\label{degree-two-affine-chart}
Any affine cycle of an algebraic foliation of degree at most two in $\mathbb{RP}^2$ can be contained in an affine chart.
\end{corollary}

\begin{proof}
By Theorem \ref{prop_convex}, the affine cycle of an algebraic foliation of degree at most two is strictly convex.
The conclusion then follows from Theorem \ref{embed_cycle}.
\end{proof}

The behavior of Hausdorff limits of such cycles will be considered in Section
\ref{proof_main}.

\section{Accumulation, finiteness and the number of limit cycles in $\mathbb{RP}^2$}

In this section, we study the accumulation,  finiteness and number of limit cycles of algebraic foliations in $\mathbb{RP}^2$. We first recall the definition of
Poincar\'e return map and  limit cycles and polycycles. We then describe the Hausdorff limits of strictly convex cycles. Finally, we determine
the maximum number of limit cycles for algebraic foliations of degrees zero and one,  and prove the finiteness of limit cycles
and give the maximum number of limit cycles at least six for algebraic foliations of degree two.

\subsection{Limit cycles and polycycles in $\mathbb{RP}^2$}

\begin{figure}[h]
    \centering
    \includegraphics[width=0.5\linewidth]{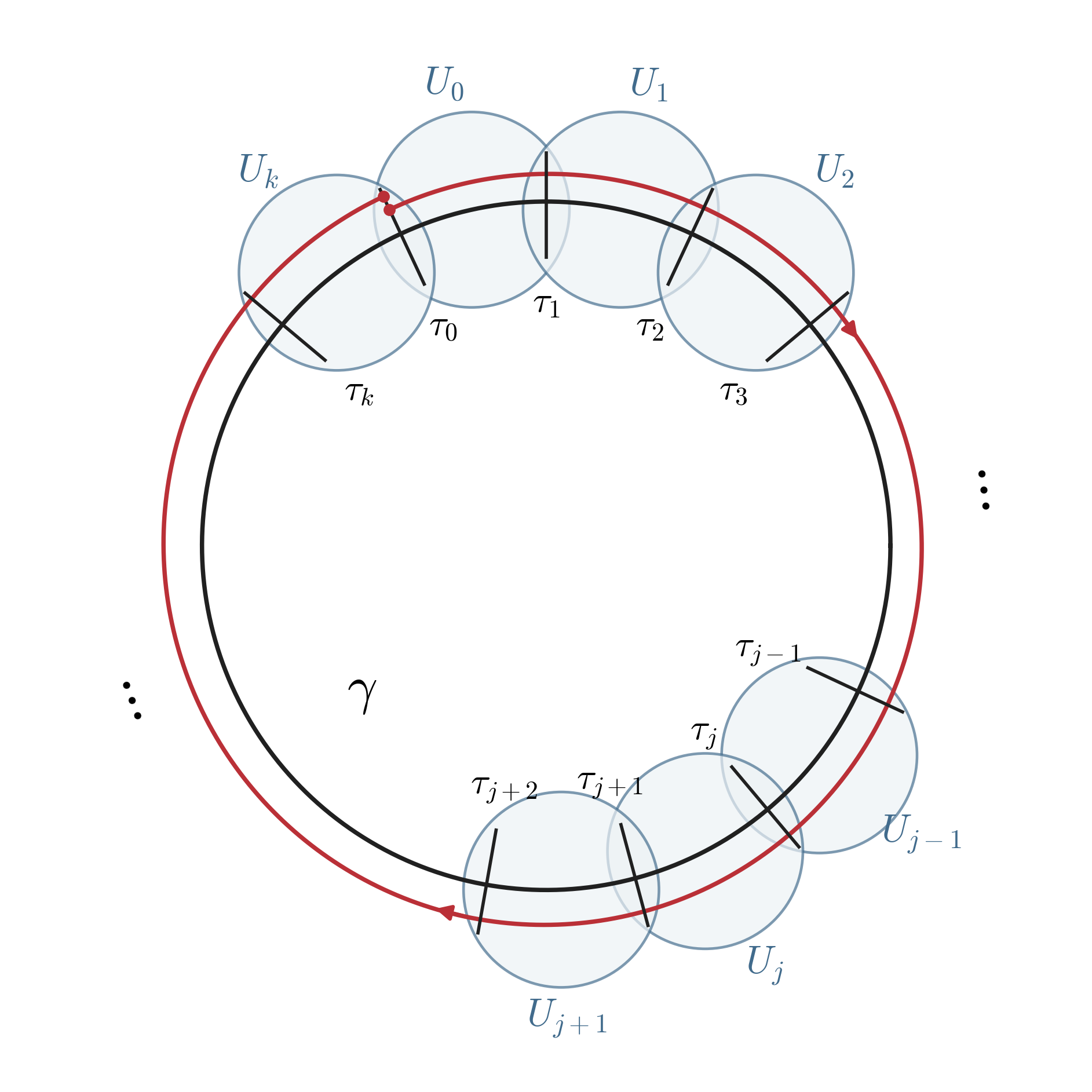}
    \caption{The construction of the holonomy map of the cycle $\gamma$.}
    \label{holonomy}
\end{figure}

To study the isolation of cycles, we recall the \emph{Poincar\'e return map} using the \emph{holonomy map} (see \cite[Chapter 1, Section 2]{Ilyashenko2008}).
Let $\gamma$ be a cycle of an algebraic foliation $\mathcal{F}$ in $\mathbb{RP}^2$, and let $p\in\gamma$.
Let $\tau$ be a transversal section through $p$, parameterized by a local analytic coordinate
$$h:(\tau,p)\longrightarrow(\mathbb{R},0),
\qquad h(p)=0.$$
The holonomy map $\Delta_\gamma:(\tau,p)\rightarrow(\tau,p)$ is defined as follows.
Since $\gamma$ is compact, it can be covered by finitely many open sets $\{U_j\}_{j=0}^{k}$ containing no singularities of $\mathcal{F}$.
By the flow-box theorem \cite[Theorem 1.14]{Ilyashenko2008}, in each $U_j$ the foliation is locally analytically equivalent to the trivial foliation $\mathcal{F}_{\alpha_0}$ defined by $\alpha_0=dx$.
More precisely, there is a local analytic map $H_j=(h_j^1,h_j^2)$ such that $\mathcal{F}|_{U_j}=H_j^*\mathcal{F}_{\alpha_0}.$
Choose intermediate transversal sections $\tau_j$, $j=1,\ldots,k$, with $\tau_j\subset U_{j-1}\cap U_j,\ \tau_0=\tau_{k+1}=\tau\subset U_k\cap U_0$ (see Figure \ref{holonomy}).
Since they are transversal to the foliation, the restrictions $h_j^1|_{\tau_j}$ and $h_j^1|_{\tau_{j+1}}$ are analytic local diffeomorphisms.
Therefore, after restricting the transversals to sufficiently small neighborhoods of their intersection points with $\gamma$, the local transition map along the leaves is given by
$\Delta_{\tau_j,\tau_{j+1}}
=
\left(h_j^1|_{\tau_{j+1}}\right)^{-1}
\circ h_j^1|_{\tau_j}$, which is an analytic local diffeomorphism.

{\it The Poincar\'e return map of $\gamma$ on $\tau$} is the composition
$$
\Delta_\gamma
=\Delta_{\tau_k,\tau_0}
 \circ\Delta_{\tau_{k-1},\tau_k}
 \circ\cdots\circ
 \Delta_{\tau_0,\tau_1}
$$
and is therefore an analytic germ. This composition is independent of the choice of the charts $U_j$ and the intermediate transversal sections $\tau_j$.
Via the coordinate $h$, it is given by an analytic germ on $(\mathbb{R},0)$ of the form
$$
\Delta_\gamma(h)=k_1h+k_2h^2+\cdots,
$$
where $k_1\neq0$ because $\Delta_\gamma$ is a local diffeomorphism.

\begin{figure}[h]
    \centering
    \includegraphics[width=0.5\linewidth]{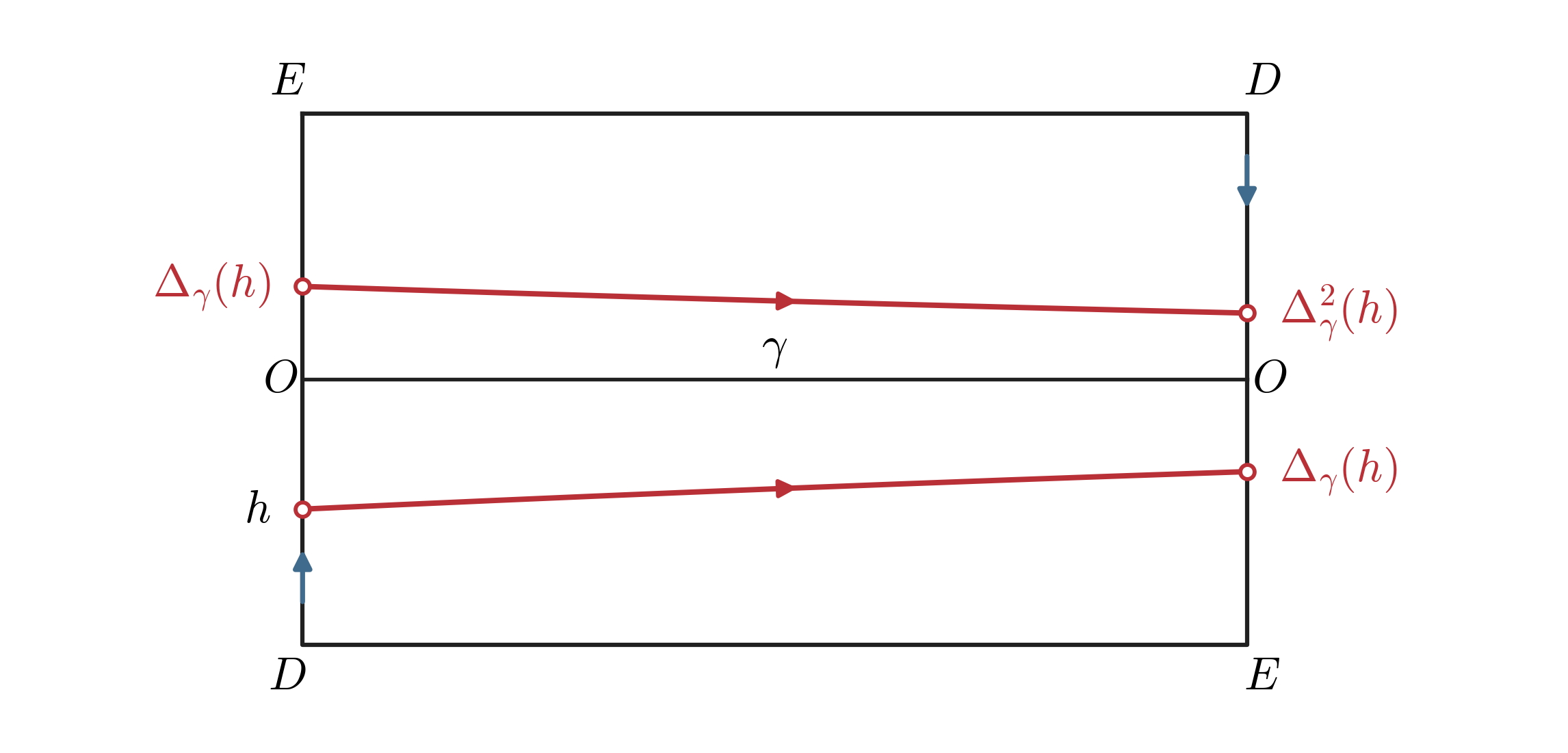}
    \caption{The Poincar\'e return map and its second iteration for a projective cycle.
    The M\"obius band is obtained by identifying the opposite vertical edges of the rectangle with reversed orientations, as indicated by the inverted labels $E$ and $D$.}
    \label{pm}
\end{figure}

Recall from Section \ref{arc} that a neighborhood of an affine cycle is homeomorphic to an annulus, whereas a neighborhood of a projective cycle is homeomorphic to a M\"obius band.
Thus, for an affine cycle, the Poincar\'e return map preserves the orientation of the transversal section.
In particular, if $h>0$ is sufficiently small, then $\Delta_\gamma(h)>0$, and hence $k_1>0$.
For a projective cycle, the return map reverses the orientation of the transversal section.
Consequently, $\Delta_\gamma(h)<0$ for every sufficiently small $h>0$, and hence $k_1<0$.
This distinction leads to the following definition.
Recall Definition \ref{cycle}: a limit cycle is defined as an isolated closed leaf in $\mathbb{RP}^2$.
We now give an equivalent definition in terms of the Poincar\'e return map, following \cite[Definition~2]{Movasati2008}.
\begin{definition}[Limit cycle \cite{Movasati2008}]
\label{definition-limit-cycle}
Let $\gamma$ be a cycle of an algebraic foliation $\mathcal{F}$ in $\mathbb{RP}^2$.
Then $\gamma$ is called an \emph{affine limit cycle} if it is an affine cycle and the germ of its Poincar\'e return map is not the identity.
$\gamma$ is called a \emph{projective limit cycle} if it is a projective cycle and the germ of the second iteration of its Poincar\'e return map is not the identity (see Figure \ref{pm}).
\end{definition}

We next extend the construction of the Poincar\'e return map to a set composed of a finite number of singularities and the leaves connecting them.
Let $\Gamma$ consist of singularities $\{s_1,\ldots,s_n\}$ and connecting leaves $\{\gamma_1,\ldots,\gamma_n\}$, oriented so that $\gamma_i$ runs from $s_i$ to $s_{i+1}$, where $s_{n+1}=s_1$.
Let $V_i$ be a sufficiently small coordinate neighborhood of $s_i$.
In $V_i$, choose a transversal section $\tau_i^+$ through a point $a_i^+$ on the incoming leaf $\gamma_{i-1}$ and a transversal section $\tau_i^-$ through a point $a_i^-$ on the outgoing leaf $\gamma_i$, where $\gamma_0=\gamma_n$.
Let
$$
h_i^+:(\tau_i^+,a_i^+)\longrightarrow(\mathbb{R},0),
\qquad
h_i^-:(\tau_i^-,a_i^-)\longrightarrow(\mathbb{R},0)
$$
be local coordinates.

To ensure that a return map is defined on a chosen side of $\Gamma$, assume that the region between $\gamma_{i-1}$ and $\gamma_i$ in $V_i$ is a hyperbolic sector.
Thus, every leaf entering $V_i$ through a sufficiently small semi-transversal section $\sigma_i^+\subset\tau_i^+$ leaves it through a semi-transversal $\sigma_i^-\subset\tau_i^-$.
This gives a local correspondence map, or Dulac map, $T_i:\sigma_i^+\rightarrow\sigma_i^-.$
Along the connecting leaf $\gamma_i$, the holonomy construction described previously gives a locally analytic map $\Delta_{\sigma_i^-,\sigma_{i+1}^+}:
\sigma_i^-\rightarrow\sigma_{i+1}^+, $ where $\sigma_{n+1}^+=\sigma_1^+.$
The return map associated with $\Gamma$ is therefore
$$
\Delta_\Gamma
=\Delta_{\sigma_n^-,\sigma_1^+}\circ T_n\circ\cdots\circ
\Delta_{\sigma_1^-,\sigma_2^+}\circ T_1.
$$
It is defined on a sufficiently small semi-transversal contained in $\sigma_1^+$.
Although the holonomy maps along the connecting leaves are analytic, the Dulac maps $T_i$ are generally not analytic at the points on $\Gamma$.
For example, they may contain logarithmic terms. Thus, $\Delta_\Gamma$ is generally defined only on one side of $\Gamma$ and need not extend analytically to $a_1^+$.

\begin{definition}[Polycycle]
\label{def:polycycle}
Let $\mathcal F$ be an algebraic foliation of degree $d$ in $\mathbb{RP}^2$,
and let $\Gamma\subset\mathbb{RP}^2$ be a nonempty closed
connected set. $\Gamma$ is called a \emph{polycycle} of
$\mathcal F$ if the following two conditions hold:
\begin{itemize}
    \item[(1)] $\Gamma$ is of one of the following types:
    \begin{itemize}
        \item[(1.1)] $\Gamma$ is either a single singularity or a
        closed leaf of $\mathcal F$.

        \item[(1.2)] $\Gamma$ is the image of a $C^0$ immersion
        $\gamma:S^1\to\mathbb{RP}^2$ and has the form
        \[
            \Gamma=\{p_1,\ldots,p_m\}\cup
            \bigcup_{i=1}^{m}\ell_i,\qquad m\ge 1,
        \]
        where $p_1,\ldots,p_m$ are isolated singularities of
        $\mathcal F$, not necessarily distinct, and
        $\ell_1,\ldots,\ell_m$ are pairwise distinct non-closed
        regular leaves of $\mathcal F$. For each $i=1,\ldots,m$,
        the two ends of $\ell_i$ converge to $p_i$ and $p_{i+1}$,
        with $p_{m+1}=p_1$.
    \end{itemize}

    \item[(2)] $\Gamma$ admits a well-defined one-sided Poincar\'e
    return map on a semi-transversal section.
    If $\Gamma$ is a single singularity, then the section has this singularity as its endpoint.
\end{itemize}
\end{definition}

\begin{remark}
    Note that Definition \ref{def:polycycle} allows a polycycle to be a single singularity, such as a focus, or a periodic orbit.
    For a singularity, the semi-transversal is a half-open arc ending at that singularity.
    For a projective cycle, the Poincar\'e return map in condition (ii) is understood to be its second iteration.
\end{remark}

As in Definition \ref{defconvex}, a polycycle $\Gamma$ is called \emph{strictly convex} if every projective line intersects $\Gamma$ in at most two distinct points. The affine realizability argument used for cycles extends to strictly convex polycycles.

{\begin{proposition}
\label{embed-polycycle}
Every strictly convex polycycle that is not a singularity in $\mathbb{RP}^2$ is a simple closed curve contained in an affine chart.
In that chart, it bounds a compact convex region.
\end{proposition}

\begin{proof}
Let $\Gamma$ be a strictly convex polycycle that is not a singularity.
If $\Gamma$ is a cycle, the conclusion follows from Proposition \ref{relationship} and Theorem \ref{embed_cycle}, using central projection in a hemisphere.
We may therefore assume that $\Gamma$ contains singularities and connecting leaves.
Among the closed sequences of connecting orbits in $\Gamma$, choose one with the
smallest number of orbits.
The singularities joining successive orbits are pairwise distinct, for otherwise, a repeated singularity determines a shorter closed sequence.
The union of the chosen orbits and their endpoint singularities therefore forms
a simple closed curve $\eta\subset\Gamma$.
Since $\Gamma$ is strictly convex, $\eta$ is also strictly convex, and so, every projective line intersects $\eta$ in at most two distinct points.

We first show that $\mathbb{RP}^2\setminus\eta$ has two connected components.
Note that no orbit of $\eta$ is contained in a projective line, since such a line intersects $\eta$ in infinitely many points, contradicting strict convexity.
We can therefore choose two non-singularities $p,q\in\eta$ such that the projective line through $p,q$ intersects $\eta$ transversally at $p,q$.
Moreover, $p,q$ are its only intersections with $\eta$.
If $\mathbb{RP}^2\setminus\eta$ were connected, the argument in the proof of
Proposition \ref{non_projective} requires the number of transversal intersections to be odd, a contradiction.

Since $\eta$ is a simple closed curve whose complement in $\mathbb{RP}^2$ has two connected components, the inverse image of $\eta$ under $\pi$ consequently consists of two disjoint simple closed curves exchanged by the antipodal map, as in the affine cycle case of Proposition~\ref{cycle-lifts}.
The proofs of Proposition \ref{relationship} and Theorem \ref{embed_cycle} also can apply to $\eta$, since their geometric arguments use its strict convexity and the fact that it is a simple closed curve, but do not require smoothness at the finitely many singularities.
Thus $\eta$ lies in an affine chart and bounds a compact convex region $D$.

Finally, we claim that $\Gamma=\eta$. Indeed, if $q\in\Gamma\setminus\eta$, we can choose a projective line $L$ through $q$ and a point in the interior of $D$.
Since $D$ is compact and convex and $L$ passes through its interior, $L$ intersects its boundary $\eta$ in two distinct points.
Together with $q\notin\eta$, these give three distinct points of $L\cap\Gamma$, contradicting strict convexity.
Hence $\Gamma=\eta$, and the conclusion follows.
\end{proof}}

For later use, we also recall some terminology concerning singularities.
A singularity is called \emph{elementary} if its Jacobian matrix has at least one nonzero eigenvalue.
It is \emph{hyperbolic} if both eigenvalues have nonzero real parts and \emph{semi-hyperbolic} if one eigenvalue is zero and the other is nonzero.
A polycycle is called \emph{hyperbolic}, respectively \emph{elementary}, if all of its singularities are hyperbolic, respectively elementary.
An elementary polycycle is called \emph{balanced} if the number of its semi-hyperbolic singularities whose nonzero eigenvalue is positive is equal to the number of those whose nonzero eigenvalue is negative.
Otherwise, it is called \emph{unbalanced}.

The following theorem is a reformulation, in the setting of algebraic foliations in $\mathbb{RP}^2$, of the non-accumulation results originating with Dulac \cite{Dulac1923} and developed by Ilyashenko \cite{Ilyashenko1985,Ilyashenko1991,Ilyashenko2002}.
Its proof is based on the analysis of local correspondence maps near the singularities.
Since these maps can be constructed in affine charts of $\mathbb{RP}^2$, with the local vector fields chosen consistently with the direction used to define the Poincar\'e return map, the argument directly extends from $\mathbb S^2$ to $\mathbb{RP}^2$.

\begin{theorem}[Non-accumulation theorem in $\mathbb{RP}^2$]
\label{nonaccumulation}
Hyperbolic polycycles and unbalanced elementary polycycles of an algebraic foliation in $\mathbb{RP}^2$ cannot be accumulated by limit cycles.
\end{theorem}

\subsection{Hausdorff limits of strictly convex cycles}

We first transfer to $\mathbb{RP}^2$ the accumulation result for
analytic vector fields on $\mathbb S^2$.
Let $\pi:\mathbb{S}^2\to\mathbb{RP}^2$ be the canonical covering map.
For $p,q\in\mathbb{RP}^2$, define
$d_{\mathbb{RP}^2}(p,q)=\min\left\{d_{\mathbb{S}^2}(\tilde p,\tilde q),d_{\mathbb{S}^2}(-\tilde p,\tilde q)\right\},$
where
$\pi^{-1}(p)=\{\tilde p,-\tilde p\},$ and $\pi^{-1}(q)=\{\tilde q,-\tilde q\}.$
The Hausdorff distance $d_{\mathbb{RP}^2}^{H}$ between nonempty compact subsets of $\mathbb{RP}^2$ is defined from this metric in the same way as the Hausdorff distance on $\mathbb{S}^2$ defined in Appendix \ref{hausdorff-convergence}.

We shall use the finiteness result of limit cycles for analytic vector fields on $\mathbb{S}^2$ with isolated singularities, that is, an infinite sequence of distinct limit cycles has a subsequence whose Hausdorff limit is a polycycle (see \cite[Section~3.2]{Roussarie1998}).
We obtain the following projective version.

\begin{theorem}[Limit sets of limit cycles in $\mathbb{RP}^2$]
\label{limit set of limit cycles}
Let $\mathcal{F}$ be an algebraic foliation with finitely many singularities in $\mathbb{RP}^2$.
If $\{\gamma_n\}_{n\in\mathbb{N}}$ is an infinite sequence of distinct limit cycles of $\mathcal{F}$, then it has a subsequence converging in the Hausdorff distance to  $\Gamma$ that is a polycycle of $\mathcal{F}$.
\end{theorem}

\begin{proof}
Let $\mathcal{F}$ be defined by a $1$-form $\alpha$ of the form \eqref{1form}, and let $\tilde{\mathcal{F}}=\pi^*\mathcal{F}$ be its pull-back to $\mathbb{S}^2$.
By Proposition \ref{spherical-representation}, the pull-back foliation is generated by the analytic vector field $V$ in \eqref{vfs2}, and its singular set is $\pi^{-1}(\operatorname{Sing}(\mathcal{F}))$.
Therefore, the vector field $V$ has finitely many singularities.

For each $n$, choose a connected component $\tilde{\gamma}_n$ of $\pi^{-1}(\gamma_n)$. Since $\pi$ is a local diffeomorphism, Proposition \ref{cycle-lifts}, together with Definition \ref{definition-limit-cycle}, shows that $\tilde{\gamma}_n$ is a limit cycle of $V$.
Note that if $\gamma_n$ is a projective cycle, the Poincar\'e return map of $\tilde \gamma_n$ corresponds to the second iteration of the Poincar\'e return map of $\gamma_n$.
Thus, $\{\tilde{\gamma}_n\}$ is an infinite sequence of distinct limit cycles of
$V$.
By Blaschke's theorem \cite[Theorem~7.3.8]{BBI2001}, after passing to a subsequence, the compact sets $\tilde{\gamma}_n$ converge in the Hausdorff distance to a nonempty compact set $\tilde{\Gamma}\subset\mathbb{S}^2$.
Since $V$ has only finitely many singularities, the finiteness result for analytic vector fields on $\mathbb{S}^2$ \cite[Chapter~2]{Roussarie1998} implies that $\tilde{\Gamma}$ is a polycycle of $\tilde{\mathcal{F}}$.

Since $\pi$ is uniformly continuous, Hausdorff convergence is preserved under
$\pi$. Consequently,
\begin{equation}
    \label{equconverge}
    \Gamma :=\lim_{n\to\infty}\gamma_n =\lim_{n\to\infty}\pi(\tilde{\gamma}_n) =\pi(\lim_{n\to+\infty}\tilde{\gamma}_n)=\pi(\tilde{\Gamma}).
\end{equation}
Because $\pi$ maps leaves to leaves and singularities to singularities, $\Gamma$ satisfies condition (i) of Definition \ref{def:polycycle}.
Furthermore, since $\tilde{\G}$ is compact, it can be covered by finitely many open sets $\{\tilde{U}_i\}_{i=1}^n$ that serve as local charts.
Moreover, we can choose them sufficiently small that $\pi$ restricts to a diffeomorphism on each $\tilde U_i$.
Similar to the analysis preceding Definition \ref{def:polycycle}, within these local charts, we can choose sufficiently many intermediate transversal sections $\tilde{\tau}_i\subset \tilde{U}_i\cap \tilde{U}_{i+1}$ for $i=1,\dots,n$ (with $\tilde{U}_{n+1}=\tilde{U}_1$), such that there exists a local correspondence map $\Delta_{\tilde{\tau}_i,\tilde{\tau}_{i+1}}:\tilde{\tau}_i\to\tilde{\tau}_{i+1}$ in each chart, and $\Delta_{\tilde{\G}}=\Delta_{\tilde{\tau}_n,\tilde{\tau}_1}\circ\Delta_{\tilde{\tau}_{n-1},\tilde{\tau}_n}\circ\dots\circ\Delta_{\tilde{\tau}_1,\tilde{\tau}_2}$.
Since $\pi$ is a local diffeomorphism, it maps each chart $\tilde{U}_i$ to a chart $U_i = \pi(\tilde{U}_i)$ covering $\Gamma$ for $i=1,2,\dots,n$.
The transversal sections project to $\tau_i = \pi(\tilde{\tau}_i) \subset U_i \cap U_{i+1}$, and the correspondence map $\Delta_{\tilde{\tau}_{i},\tilde{\tau}_{i+1}}$ naturally descends to a well-defined transition map from $\tau_i$ to $\tau_{i+1}$, given by $\Delta_{{\tau}_{i},{\tau}_{i+1}} = \pi \circ \Delta_{\tilde{\tau}_{i},\tilde{\tau}_{i+1}} \circ (\pi|_{\tilde{U}_i})^{-1}$.
Consequently, the Poincar\'e return map on $\tau_1$ can be defined by
    \begin{align*}
        \Delta_\G&=\Delta_{{\tau}_{n},{\tau}_{1}}\circ\Delta_{{\tau}_{n-1},{\tau}_{n}}\circ\dots\circ\Delta_{{\tau}_{1},{\tau}_{2}}\\
        &= \left(\pi \circ \Delta_{\tilde{\tau}_{n},\tilde{\tau}_{1}} \circ (\pi|_{\tilde{U}_n})^{-1}\right) \circ \left(\pi \circ \Delta_{\tilde{\tau}_{n-1},\tilde{\tau}_{n}} \circ (\pi|_{\tilde{U}_{n-1}})^{-1}\right) \circ \dots \circ \left(\pi \circ \Delta_{\tilde{\tau}_{1},\tilde{\tau}_{2}} \circ (\pi|_{\tilde{U}_1})^{-1}\right)\\
        &= \pi \circ \Delta_{\tilde{\tau}_n,\tilde{\tau}_1} \circ \Delta_{\tilde{\tau}_{n-1},\tilde{\tau}_n} \circ \dots \circ \Delta_{\tilde{\tau}_1,\tilde{\tau}_2} \circ (\pi|_{\tilde{U}_1})^{-1}\\
        &= \pi \circ \Delta_{\tilde{\G}} \circ (\pi|_{\tilde{U}_1})^{-1}.
    \end{align*}
    Therefore, $\Gamma$ also satisfies condition (ii). Thus, we conclude that the Hausdorff limit set $\Gamma$ is a polycycle of $\mathcal{F}$ in $\mathbb{RP}^2$.
\end{proof}

Although each strictly convex cycle is contained in an affine chart, a Hausdorff limit of strictly convex cycles need not be strictly convex.
For example, a sequence of affine cycles may converge to a projective line.
Nevertheless, this polycycle still satisfies a weaker property, which is stated in the following theorem.

\begin{theorem}
\label{convex_polycycle}
Let $\mathcal{F}$ be an algebraic foliation in $\mathbb{RP}^2$, and let $\Gamma$ be a polycycle that is the Hausdorff limit of a sequence of strictly convex affine cycles of $\mathcal{F}$.
If $\Gamma$ is not strictly convex, then it has a leaf contained in an invariant line of $\mathcal{F}$.
\end{theorem}

\begin{proof}
Let $\{\gamma_n\}$ be a sequence of strictly convex affine cycles of $\mathcal{F}$ converging to $\Gamma$.
By Proposition \ref{cycle-lifts}, the inverse image of $\gamma_n$ on $\mathbb{S}^2$ via the canonical covering map $\pi$ consists of two connected components $\tilde{\gamma}_n$ and $-\tilde{\gamma}_n$.
Let $K_n$ be the closed region bounded by $\tilde{\gamma}_n$ that contains no pair of antipodal points.
By Proposition \ref{relationship}, $K_n$ is convex.

Since $\mathbb{S}^2$ is compact, Blaschke's theorem gives a subsequence of $\{K_n\}$, still denoted by $\{K_n\}$, converging in the Hausdorff distance to a compact set $K$.
Proposition \ref{converge} shows that $K$ is convex.
Since $K_n\cap(-K_n)=\varnothing$ for every $n$, convexity and Hausdorff convergence imply that $\operatorname{int}K\cap\operatorname{int}(-K)=\varnothing$, and hence $K\ne\mathbb S^2$.
Furthermore, passing to another subsequence if necessary, according to Proposition \ref{boundary}, one has
$$
\tilde{\Gamma} := \lim_{n\to\infty} \tilde{\gamma}_n = \lim_{n\to\infty} \partial K_n = \partial K.
$$
Moreover, taking into account \eqref{equconverge}, we have
$\Gamma=\pi(\tilde{\Gamma})=\pi(\partial K).$
Let $\tilde{L}$ be an arbitrary great circle on $\mathbb{S}^2$, and let $L=\pi(\tilde{L})$ be the corresponding projective line.
We distinguish two cases.

\smallskip
\noindent\emph{Case 1:}
$\tilde{L}\cap\operatorname{int}(K)\neq\varnothing$.
We claim that in this case, $\tilde{L} \cap \partial K$ contains at most two points.
Note that in this case, $\tilde{L} \cap K$ can not be a great circle, for otherwise $\operatorname{int}(K)$ contains a point on this great circle, which implies that $K=\mathbb{S}^2$ by the convexity of $K$ and Lemma \ref{interior}.
Thus, $\tilde{L}\cap K$ is a minor arc or a semicircle, denoted by $[x,y]$.
Choose $z\in\tilde{L}\cap\operatorname{int}(K)$.
By Lemma \ref{interior}, $(x,z]\subset\operatorname{int}(K)$ and $[z,y)\subset\operatorname{int}(K)$.
Therefore, $(x,y)\subset\tilde{L}\cap\operatorname{int} K,$ and so,
$\tilde{L}\cap\partial K\subset\{x,y\}.$
Hence, $L\cap\Gamma$ contains at most two distinct points $\pi(x),\ \pi(y)$.

\smallskip
\noindent\emph{Case 2:}
$\tilde{L}\cap\operatorname{int}(K)=\varnothing$.
Then $\tilde{L}\cap K\subset\partial K$.
If every projective line intersects $\Gamma$ in at most two distinct points, then $\Gamma$ is strictly convex by Definition \ref{defconvex}.
Therefore, we suppose that $L\cap\Gamma$ contains more than two distinct points, and choose distinct points $p,q\in L\cap\Gamma$ and lifts $\tilde p\in \pi^{-1}(p),\ \tilde q\in \pi^{-1}(q)$ such that $\tilde p,\ \tilde q\in \tilde L \cap\partial K$.
Since $p \neq q $, $\tilde p$ and $ \tilde q$ are not antipodal.
By convexity, the minor arc $[\tilde p,\tilde q]\subset K$.
Taking into account that $\tilde{L}\cap\operatorname{int}(K)=\varnothing$ and $[\tilde p, \tilde q]\subset \tilde L$, it gives $[\tilde p, \tilde q]\subset\tilde{L}\cap\partial K.$
Consequently, $\pi([\tilde p,\tilde q])\subset L\cap\Gamma$, and so, $\mathcal{F}$ is tangent to $\pi([\tilde p,\tilde q])\setminus\operatorname{Sing}(\mathcal{F})$, which is a subset of $L$.
Since $\mathcal F$ has only finitely many singularities, one has
$\#\bigl(\pi([\tilde p,\tilde q])\setminus
\operatorname{Sing}(\mathcal F)\bigr)=\infty.$
By Proposition \ref{contact}, a non-invariant line has only finitely many points of tangency with $\mathcal F$.
Hence $L$ is invariant.
Choose $r\in\pi([\tilde p,\tilde q])\setminus\operatorname{Sing}(\mathcal F)$,
and let $\ell$ be the leaf of $\mathcal F$ through $r$.
Since $\Gamma$ is a polycycle and $L$ is invariant, we have $\ell\subset\Gamma\cap L$ is a leaf contained in an invariant line of $\mathcal F$.

Thus, for every projective line $L$, either $L\cap\Gamma$ contains at most two points, or $\Gamma$ has a leaf contained in an invariant line.
If $\Gamma$ is not strictly convex, the latter case must occur.
\end{proof}

\subsection{The number of limit cycles for algebraic foliations of degrees zero and one in $\mathbb{RP}^2$}

We now determine the maximum number of limit cycles for algebraic foliations of degrees zero and one in $\mathbb{RP}^2$. We begin with a property about invariant lines
for degree-one foliations.

\begin{proposition}
\label{1-invariant}
An algebraic foliation of degree one in $\mathbb{RP}^2$ has at least one invariant line.
\end{proposition}

\begin{proof}
Let $\mathcal{F}$ be an algebraic foliation of degree one in $\mathbb{RP}^2$.
By Proposition \ref{spherical-representation}, its pull-back under $\pi:\mathbb{S}^2\to\mathbb{RP}^2$ is determined by a vector field $V$ on $\mathbb{S}^2$.
The Poincar\'e--Hopf Index Theorem \cite[Section 3.5]{GP1974} implies that $V$ has at least one singularity, because the Euler characteristic of $\mathbb{S}^2$ is two. Consequently, $\mathcal{F}$ has at least one singularity.
After a change of coordinates, we may suppose that this singularity is $[0:0:1]$.
In the affine chart $\{[X:Y:Z]\in\mathbb{RP}^2:Z=1\}$, it is the origin, and $\mathcal F$ is defined by the vector field
$$\mathcal{X}_1
=\bigl(a_1x+b_1y+x(ax+by)\bigr)\frac{\partial}{\partial x}
+\bigl(a_2x+b_2y+y(ax+by)\bigr)\frac{\partial}{\partial y}.$$
If $a=b=0$, then $\mathcal{X}_1$ is linear and the line at infinity $\{[X:Y:Z]\in\mathbb{RP}^2:Z=0\}$ is invariant.
We therefore consider the case $a^2+b^2\neq0$. Let
$$
A=D\mathcal{X}_1(0)
=\begin{pmatrix}a_1&b_1\\a_2&b_2\end{pmatrix}.
$$
We distinguish two cases according to $\det A$.

Suppose first that $\det A\neq0$. Define the affine line $l:\;kx+my+n=0, $
where
$$ k=a_2b-b_2a, \qquad m=b_1a-a_1b, \qquad n=a_2b_1-a_1b_2. $$
Since $\det A\neq0$ and $a^2+b^2\neq0$, one has $k^2+m^2\neq0$, and so, $l$ is well defined.
A direct calculation gives
$$\begin{aligned}
\mathcal{X}_1(kx+my+n)
&=k\bigl(a_1x+b_1y+x(ax+by)\bigr)\\
&\quad+m\bigl(a_2x+b_2y+y(ax+by)\bigr)\\
&=(kx+my+n)(ax+by).
\end{aligned}$$
Thus, the projective closure of $l$ is invariant under $\mathcal{F}$.

Suppose now that $\det A=0$.
Make the rotation $\bar{x}=x\cos\theta+y\sin\theta, \bar{y}=-x\sin\theta+y\cos\theta$ for $\theta\in [0,2\pi).$
In the $(\bar{x},\bar{y})$ coordinates, the vector field $\mathcal X_1$ has the form
$$
\mathcal{X}_1:
\begin{cases}
\dot{\bar{x}}
=\bar a_1\bar{x}+\bar b_1\bar{y}
+\bar{x}(\bar a\bar{x}+\bar b\bar{y}),\\
\dot{\bar{y}}
=\bar a_2\bar{x}+\bar b_2\bar{y}
+\bar{y}(\bar a\bar{x}+\bar b\bar{y}),
\end{cases}
$$
where
$$
\begin{aligned}
\bar a_1
&=a_1\cos^2\theta+(a_2+b_1)\sin\theta\cos\theta
+b_2\sin^2\theta,\\
\bar b_1
&=b_1\cos^2\theta-(a_1-b_2)\sin\theta\cos\theta
-a_2\sin^2\theta,\\
\bar a_2
&=a_2\cos^2\theta-(a_1-b_2)\sin\theta\cos\theta
-b_1\sin^2\theta,\\
\bar b_2
&=b_2\cos^2\theta-(a_2+b_1)\sin\theta\cos\theta
+a_1\sin^2\theta,\\
\bar a&=a\cos\theta+b\sin\theta,
\qquad
\bar b=-a\sin\theta+b\cos\theta.
\end{aligned}
$$
Since $A$ is a real matrix with $\det A=0$, there exists a nonzero real vector $v$ such that $Av=0$.
Choose a rotation through an angle $\theta$ so that the direction of $v$ becomes
the $\bar x$-axis.
Then $\bar a_2=0$, and
$$
\dot{\bar{y}}
=\bar{y}\bigl(\bar b_2+\bar a\bar{x}+\bar b\bar{y}\bigr).
$$
Consequently, the line $\{\bar{y}=0\}$ is invariant under $\mathcal{X}_1$.
Its projective closure is an invariant line of
$\mathcal{F}$.

The two cases complete the proof.
\end{proof}

Given Proposition \ref{1-invariant}, we determine the upper bound on the number of limit cycles, both affine and projective, for algebraic foliations of degrees zero and one in $\mathbb{RP}^2$.

\begin{theorem}
\label{01}
Algebraic foliations of degree zero in $\mathbb{RP}^2$ have no limit cycles.
Algebraic foliations of degree one in $\mathbb{RP}^2$ have at most one limit cycle.
Moreover, if it exists, it is a projective cycle and is itself a projective line.
\end{theorem}

\begin{proof}
Consider first that the algebraic foliation $\mathcal{F}$  in $\mathbb{RP}^2$ is of degree zero.
By Theorem \ref{even}, an algebraic foliation of degree zero in $\mathbb{RP}^2$ has no projective cycle.
Moreover, Theorems \ref{embed_cycle} and \ref{prop_convex} show that every affine cycle of a zero-degree foliation is contained in an affine chart.
By Proposition \ref{cases}, in this affine chart, the foliation $\mathcal{F}$ is given by
$$(a_0+\lambda x)\frac{\partial}{\partial x}
+(b_0+\lambda y)\frac{\partial}{\partial y}.$$
If $\lambda\neq0$, this vector field has a global star node and no affine cycles.
If $\lambda=0$, then $a_0^2+b_0^2\neq0$ by the coprimality condition, and the vector field is constant and has no affine cycles.
Therefore, algebraic foliations of degree zero in $\mathbb{RP}^2$ have no limit cycles.

Now suppose that $\mathcal{F}$ has degree one.
By Proposition \ref{1-invariant}, it has an invariant line $L$. Taking $L$ as the line at infinity $L_\infty$,  Proposition \ref{cases} and Corollary \ref{non-generic} show that the restriction of $\mathcal{F}$ to the affine chart $\mathbb{RP}^2\setminus L_\infty$ is defined by a linear vector field
$$(a_1+a_2x+a_3y)\frac{\partial}{\partial x}
+(b_1+b_2x+b_3y)\frac{\partial}{\partial y},$$
where $a_{1}+a_{2} x+a_{3} y$ and $b_{1}+b_{2} x+b_{3}y$ are coprime.
Note that a cycle distinct from $L_\infty$ cannot intersect $L_\infty$, since two distinct leaves cannot intersect.
Hence every affine cycle lies in this affine chart.
Since a linear vector field in $\mathbb{R}^2$ has no limit cycles, $\mathcal{F}$ has no affine limit cycles.

If $\mathcal{F}$ has a projective cycle $\gamma$, then $\gamma=L_\infty$ since distinct leaves cannot intersect.
Thus, there is at most one projective limit cycle.
Moreover, this upper bound is attained.
For instance, consider the foliation in $\mathbb{RP}^2$ determined in an affine chart by
$$(x-y)\frac{\partial}{\partial x}
+(x+y)\frac{\partial}{\partial y}.$$
Here, $x(x+y)-y(x-y)=x^2+y^2$, which has no real linear factor, so the invariant line at infinity contains no singularities.
The unique singularity in the affine chart is the focus at the origin, and the line at infinity is a projective limit cycle.
Therefore, the maximum number of limit cycles for a algebraic foliation of degree one in $\mathbb{RP}^2$ is one.
\end{proof}

\subsection{The number of limit cycles for algebraic foliations of degree two in $\mathbb{RP}^2$}
\label{proof_main}

We now turn to algebraic foliations of degree two in $\mathbb{RP}^2$.
We first record the projective form of Bamon's theorem \cite{Bamon1986}.

\begin{theorem}[Bamon's theorem in $\mathbb{RP}^2$]
\label{bamon}
An algebraic foliation of degree two in $\mathbb{RP}^2$ with an invariant line has only finitely many limit cycles.
\end{theorem}

\begin{proof}
Choose the invariant line as the line at infinity $L_\infty$.
By Proposition \ref{cases} and Corollary \ref{non-generic}, the foliation is given by a quadratic polynomial vector field in the affine chart $\mathbb{RP}^2\setminus L_\infty$.
By Bamon's theorem \cite{Bamon1986}, this affine vector field has only finitely many limit cycles.
Theorem \ref{even} excludes projective cycles for degree-two foliations, so these limit cycles in this affine chart are all the limit cycles of the foliation.
\end{proof}

We  now prove the finiteness of limit cycles for any algebraic foliation of degree two in $\mathbb{RP}^2$.

\begin{theorem}
\label{main2}
Every algebraic foliation of degree two in $\mathbb{RP}^2$ has only finitely many limit cycles.
\end{theorem}

Moreover, we establish a lower bound for the maximum number of limit cycles for  algebraic foliations of degree two in $\mathbb{RP}^2$.

\begin{theorem}
\label{main3}
The maximal number of limit cycles in algebraic foliation of degree two in $\mathbb{RP}^2$ is at least six.
\end{theorem}

Note that  an algebraic foliation of degree two in $\mathbb{RP}^2$ has no projective cycles by Theorem \ref{even}.
Moreover, by Theorem \ref{bamon}, Theorem \ref{main2} already holds if the foliation has an invariant line.
Therefore, we only need to  consider the generic case.

If a degree-two foliation $\mathcal{F}$ has infinitely many limit cycles, then, by Theorem \ref{prop_convex}, all of these limit cycles are strictly convex.
Moreover, Theorem \ref{limit set of limit cycles} shows that there exists a subsequence of these limit cycles converging in the Hausdorff distance to a compact invariant set $\Gamma$, which is a polycycle.
Let $\Gamma$ be such a polycycle.
Theorem \ref{convex_polycycle} shows that either $\Gamma$ is strictly convex, or it has a leaf contained in an invariant line.
The latter case is proved by Theorem \ref{bamon}.
In the former case, if $\Gamma$ is a singularity, we can choose an affine chart containing it.
Otherwise,  by Proposition \ref{embed_cycle} and Proposition \ref{embed-polycycle}, $\Gamma$ is contained in an affine chart.
Our strategy is to analyze all possible configurations of polycycles in the affine chart, including singularities, cycles, and complicated polycycles with at least one singularity and one leaf.
We will show that none of these configurations can serve as the limit set for an infinite sequence of strictly convex limit cycles.
We now proceed with a case-by-case analysis of all possible polycycles to establish this contradiction.

By Proposition \ref{cases}, in an affine chart, the foliation for which the line at infinity is not invariant is given by a vector field of the form
\begin{equation}
\label{gvf}
\mathcal{X}:
\begin{cases}
\dot{x}=p(x,y)+xr(x,y)\\
\phantom{\dot{x}}
=a_{00}+a_{10}x+a_{01}y+a_{20}x^2+a_{11}xy+a_{02}y^2
+x(ax^2+bxy+cy^2),\\
\dot{y}=q(x,y)+yr(x,y)\\
\phantom{\dot{y}}
=b_{00}+b_{10}x+b_{01}y+b_{20}x^2+b_{11}xy+b_{02}y^2
+y(ax^2+bxy+cy^2),
\end{cases}
\end{equation}
where $a^2+b^2+c^2\neq0$ since the line at infinity is not invariant.
Note that this form is preserved by translations and rotations.
Since a limit cycle in the affine plane must enclose at least one singularity of the foliation by the Index Theory \cite[Section~3.12]{Perko2013}, we can translate one singularity to the origin and hence assume $a_{00}=b_{00}=0$.


\subsubsection{Accumulation at a singularity or a cycle}

We first consider elementary singularities and cycles.
The following proposition holds for algebraic foliations of arbitrary degree in $\mathbb{RP}^2$.

\begin{proposition}
\label{elementary-cycle-nonaccumulation}
Limit cycles of an algebraic foliation in $\mathbb{RP}^2$ cannot accumulate to an elementary singularity or a cycle.
\end{proposition}

\begin{proof}
We prove the proposition by contradiction.

Assume first that an elementary singularity of an algebraic foliation in $\mathbb{RP}^2$ is the Hausdorff limit of a sequence of distinct limit cycles.
Choose an affine chart containing this singularity.
By Hausdorff convergence, all sufficiently large terms of the sequence are contained in a sufficiently small neighborhood of the singularity in this chart.
By a translation, we move the singularity to the origin $O$.
Since the singularity is elementary, the Jacobian matrix of the vector field determining the foliation has at least one nonzero eigenvalue.
Moreover, it is surrounded by periodic orbits arbitrarily close to it.
By the classification of elementary singularities, its Jacobian matrix must therefore have a pair of nonreal complex conjugate eigenvalues.
Hence, after a suitable linear change of coordinates, the leaves of the foliation near the origin are given by
\begin{equation}\label{elementary-dxdy}
(fx+gy+F_1(x,y))\,dy-(-gx+fy+F_2(x,y))\,dx=0,
\end{equation}
where $f,g\in\mathbb R$, $g\neq0$, and $F_1,F_2$ are polynomials vanishing to order at least two at the origin.
The polar coordinate transformation $(x,y)=(r\cos\theta,r\sin\theta)$, followed by division by $r$, gives that for $0<r\ll1$,
\begin{equation}\label{elementary-drdtheta}
(g+O(r))\,dr+(fr+O(r^2))\,d\theta=0.
\end{equation}
For $r>0$, the polar coordinate transformation is an analytic diffeomorphism from $(0,+\infty)\times\mathbb S^1$ onto $\R^2\setminus\{O\}$.
Hence the orbits of \eqref{elementary-dxdy} and \eqref{elementary-drdtheta} are equivalent in a sufficiently small punctured neighborhood of $O$.
Consequently, nearby limit cycles correspond to isolated zeros of the displacement function $d_1(r_0):=r(2\pi;r_0)-r_0$ for $r_0>0$, where $r(\theta;r_0)$ is the solution of \eqref{elementary-drdtheta} satisfying $r(0;r_0)=r_0$.
Our initial assumption gives distinct positive numbers $r_n\to0$ such that $d_1(r_n)=0$.
More precisely, since $g\ne0$, $-g+O(r)\neq0$ for $|r|\ll1$ uniformly in $\theta$, equation \eqref{elementary-drdtheta} change to
$$\frac{dr}{d\theta}=\frac{fr+O(r^2)}{-g+O(r)},\ 0<r\ll1.$$
which admits a natural real-analytic extension across $r=0$ to $r<0$.
Analytic dependence over $0\le\theta\le2\pi$ therefore shows that $d_1$ is analytic on $(-R_1,R_1)$ for some $R_1>0$.
By the Identity Theorem for analytic functions \cite[Chapter~IV, Theorem~3.7]{conway11}, a sequence of zeros $\{r_n\}$ converging to the origin forces $d_1(r_0)\equiv0$ in $(0, R_1)$.
Thus the origin $O$ is a center, and a punctured neighborhood of $O$ is filled with non-isolated closed leaves.
This contradicts the fact that these cycles are limit cycles.
Therefore, limit cycles cannot accumulate to an elementary singularity.

Assume next that a sequence of distinct limit cycles converges in the Hausdorff distance to an affine cycle $\gamma_0$.
Let $\Sigma$ be a transversal section to $\gamma_0$, and let $\Delta_{\gamma_0}(h)$ be the Poincar\'e return map on $\Sigma$, where $h$ is a local coordinate such that $h=0$ corresponds to $\Sigma\cap\gamma_0$.
Since the foliation is analytic and $\gamma_0$ is a cycle, $\Delta_{\gamma_0}$ is analytic in a neighborhood of $h=0$, with $\Delta_{\gamma_0}(0)=0$.
If limit cycles accumulate to $\gamma_0$, then there exists a sequence of distinct zeros $\{h_n\}$ for the displacement function $d_2(h):=\Delta_{\gamma_0}(h)-h$, such that $h_n \to 0$ as $n\to\infty$.
By the Identity Theorem for analytic functions again, $d_2\equiv0$ in $(-R_2, R_2)$ for some $R_2>0$, so that $\Delta_{\gamma_0}(h)=h$.
This implies that a neighborhood of $\gamma_0$ is filled with closed leaves, contradicting the isolation of limit cycles.
If $\gamma_0$ is a projective cycle, its Poincar\'e return map reverses orientation on $\Sigma$.
Nearby closed leaves therefore correspond to fixed points of its second iteration.
Applying the same argument to $d_2(h):=\Delta_{\gamma_0}^{\,2}(h)-h,$ where $\Delta_{\gamma_0}^{\,2}=\Delta_{\gamma_0}\circ\Delta_{\gamma_0}$, gives the same contradiction.
Thus, limit cycles cannot accumulate to a cycle.
\end{proof}

It remains to consider nilpotent and degenerate singularities.
A singularity is called \emph{nilpotent} if its Jacobian matrix is a nonzero nilpotent matrix, and \emph{degenerate} if its Jacobian matrix is the zero matrix.
Recall that a singularity is \emph{monodromic} if the leaves in a punctured neighborhood induce a well-defined Poincar\'e return map on a semi-transversal section.
In fact, we prove the following proposition, thereby generalizing the classic result that nilpotent and degenerate singularities of a quadratic vector field in $\mathbb{R}^2$ cannot be monodromic singularities \cite{Ye1986}.

\begin{proposition}
\label{non-elementary-monodromic}
An algebraic foliation of degree two in $\mathbb{RP}^2$ has no non-elementary monodromic singularities.
\end{proposition}

\begin{proof}
Let $\mathcal{F}$ be an algebraic foliation in $\mathbb{RP}^2$ given by \eqref{gvf}.
The proof is divided into the nilpotent and degenerate cases.

Suppose first that the origin is a nilpotent singularity.
After a linear change of coordinates, we may assume that
$$a_{10}=b_{10}=b_{01}=0, \qquad a_{01}=1. $$
We distinguish two cases according to $b_{20}$.

\begin{enumerate}
    \item[(i)] If $b_{20}=0$, then $\dot y=0$ on the line $\{y=0\}$.
    Hence, $\{y=0\}$ is invariant under $\mathcal{F}$.
    As leaves of a foliation cannot intersect, the line $\{y=0\}$ acts as a local barrier, preventing nearby leaves from surrounding the origin.
    Consequently, the origin is not monodromic.

    \item[(ii)]If $b_{20} \neq 0$, we rewrite the vector field \eqref{gvf} in the following form
    \begin{equation}
    \label{vfcusp}
        \dot{x}=y+M(x,y),\dot{y}=N(x,y).
    \end{equation}
    By setting $\dot{x}=0$, we obtain the implicit function
    $$y=\phi(x)=-a_{20}x^2+O(x^3),$$
    and
    $$F(x)=N(x,\phi(x))=b_{20}x^2+O(x^3),$$
    $$G(x)=M_x(x,\phi(x))+N_y(x,\phi(x))=(2a_{20}+b_{11})x+O(x^2).$$
    Since the degree of the leading term of $F(x)$ is two, and either $G(x)\equiv0$ or the degree of the leading term of $G(x)$ is no less than one, according to the Nilpotent Singular Points Theorem \cite{DLA2006, Zhang1992}, the origin is a cusp, a type of singularity that is inherently non-monodromic.
\end{enumerate}

Thus, a nilpotent singularity cannot be monodromic.

Now suppose that the origin is a degenerate singularity.
Then $p$ and $q$ in \eqref{gvf} are homogeneous quadratic polynomials.
Define
$$E(x,y)=xq(x,y)-yp(x,y).$$
If $E(x,y)\equiv0$, then $x\dot y-y\dot x\equiv0$, so every line through the origin is invariant and the origin is not monodromic.
We may therefore assume that $E(x,y)\not\equiv0$.
Then we have
\begin{align*}
    \mathcal{X}(E) &= \frac{\partial E}{\partial x}(p+xr) + \frac{\partial E}{\partial y}(q+yr) \\
                   &= (q + x\frac{\partial q}{\partial x} - y\frac{\partial p}{\partial x})(p+xr) + (x\frac{\partial q}{\partial y} - p - y\frac{\partial p}{\partial y})(q+yr).
\end{align*}
Since in this case $p$ and $q$ are quadratic homogeneous polynomials, Euler's Theorem states that $x\frac{\partial p}{\partial x} + y\frac{\partial p}{\partial y} = 2p$ and $x\frac{\partial q}{\partial x} + y\frac{\partial q}{\partial y} = 2q$.
Applying Euler's theorem and regrouping the expression for $\mathcal{X}(E)$ shows that
\begin{align*}
    \mathcal{X}(E) = ( 3r + \frac{\partial p}{\partial x} + \frac{\partial q}{\partial y} ) E.
\end{align*}
This confirms that $E(x,y)$ is a Darboux polynomial with cofactor
$K(x,y) = 3r + \frac{\partial p}{\partial x} + \frac{\partial q}{\partial y}.$
Since $p(x,y)$ and $q(x,y)$ are quadratic homogeneous polynomials, $E(x,y)$ is a cubic homogeneous polynomial.
By the Fundamental Theorem of Algebra, this polynomial must have at least one real linear factor, which corresponds to an invariant line passing through the origin.
The existence of such an invariant line precludes the singularity from being monodromic.
Therefore, no degenerate singularity in the foliation $\mathcal{F}$ given by \eqref{gvf} can be a monodromic singularity.
Thus, we complete the proof.
\end{proof}

\begin{corollary}
\label{non-elementary-not-limit}
A nilpotent or degenerate singularity of an algebraic foliation of degree two in $\mathbb{RP}^2$ cannot be the Hausdorff limit of a sequence of limit cycles.
\end{corollary}

\begin{proof}
A singularity accumulated by cycles is necessarily monodromic. The conclusion
therefore follows from Proposition \ref{non-elementary-monodromic}.
\end{proof}

Proposition \ref{elementary-cycle-nonaccumulation} and Corollary \ref{non-elementary-not-limit} exclude the first two possible types of polycycles $\Gamma$, namely a singularity and a cycle.
It remains to exclude more complicated polycycles.

\subsubsection{Accumulation to a polycycle with one singularity}

We first analyze a polycycle with exactly one singularity and at least one leaf.

\begin{lemma}
\label{one-singularity-polycycle}
Limit cycles cannot accumulate to a polycycle with exactly one singularity and at least one leaf of an algebraic foliation of degree two in $\mathbb{RP}^2$.
\end{lemma}

\begin{proof}
Suppose that such a polycycle is accumulated by limit cycles.
By Theorem \ref{nonaccumulation}, its singularity cannot be elementary, since an elementary polycycle with one singularity is either hyperbolic or unbalanced.
Hence the singularity is nilpotent or degenerate.
We analyze these two cases separately.

Let $\mathcal{F}$ be given by \eqref{gvf}.
Suppose first that the singularity is nilpotent and located at the origin.
After a linear change of coordinates, we may assume that
$$
a_{10}=b_{10}=b_{01}=0,
\qquad a_{01}=1.
$$
If $b_{20}=0$, then the line $\{y=0\}$ is invariant.
Theorem \ref{bamon} then gives the finiteness of the limit cycles, contradicting the assumed accumulation.
If $b_{20}\neq0$, then as shown in the proof of Proposition \ref{non-elementary-monodromic}, the singularity is a cusp, and the vector field has the form \eqref{vfcusp}.
If necessary, the linear change $(x,y)\mapsto(-x,-y)$ preserves the linear part $\dot x=y$ and reverses the sign of $b_{20}$, and so, we can assume that $b_{20}>0$.
Therefore, the two characteristic orbits of the cusp have Puiseux expansions
$$f^\pm(x)=\pm c_0x^{\frac{3}{2}}+c^\pm_1x^2+c_2^\pm x^{\frac{5}{2}}+o(x^\frac{5}{2}),
\qquad c_0>0.$$
We claim that the polycycle cannot be strictly convex.
Choose a sufficiently small $x_0>0$, and consider the line $l$ passing through the two points, $(x_0,f^+(x_0))$ and $(2x_0,f^+(2x_0))$, see Figure \ref{fig_cusp}.
The equation of $l$ is
$$y(x)=\frac{f^+(2x_0)-f^+(x_0)}{x_0}x+b,$$
where
\begin{align*}
        b&=2f^+(x_0)-f^+(2x_0)\\&
        =2( c_0x_0^{\frac{3}{2}}+c^+_1x_0^2+o(x_0^2))-(c_0(2x_0)^{\frac{3}{2}}+c^+_1(2x_0)^2+o(x_0^2))\\
        &=2(1-\sqrt{2})c_0x_0^\frac{3}{2}+O(x_0^2).
    \end{align*}
Since $c_0 > 0$ and $1 - \sqrt{2} < 0$, for sufficiently small $x_0$, the dominant term is negative, implying $b < 0$.
Hence, $y(0) = b < 0 = f^-(0)$.
On the other hand, by construction, we have $y(x_0) = f^+(x_0)$.
Since $c_0 > 0$, for small $x_0$ we have $f^+(x_0) > 0$ and $f^-(x_0) < 0$.
Thus, $y(x_0)=f^+(x_0)>0>f^-(x_0)$.
The Intermediate Value Theorem gives an $x_*\in(0,x_0)$ such that $y(x_*)=f^-(x_*)$. Thus, the line $l$ intersects the polycycle in at least three distinct points, namely, $(x_*,f^-(x_*)),\ (x_0,f^+(x_0))$, and $(2x_0,f^+(2x_0)).$
Therefore, the polycycle is not strictly convex.
By Theorem \ref{prop_convex}, the limit cycles are strictly convex, and Theorem \ref{convex_polycycle} claims that the non-strictly-convex polycycle limit must have a leaf contained in an invariant line.
However, the connecting leaves of the polycycle cannot be contained in an invariant line, since locally they coincide with the two characteristic cusp orbits displayed above, whose leading terms are $\pm c_0x^\frac32$.
This contradiction excludes the cusp case.

\begin{figure}[h]
    \centering
    \includegraphics[width=0.4\linewidth]{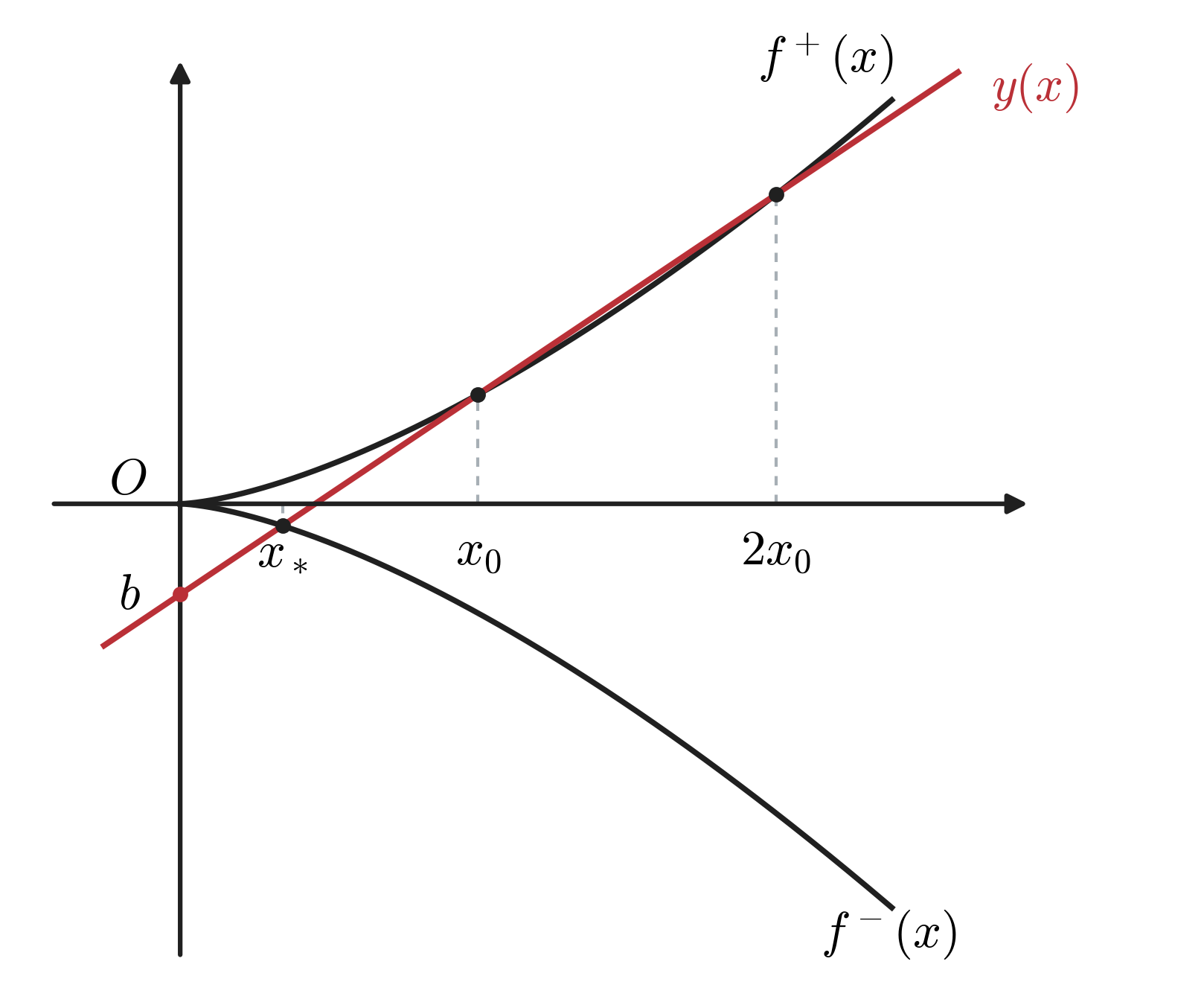}
    \caption{A cusp of $\mathcal{F}$.}
    \label{fig_cusp}
\end{figure}

Finally, if the singularity is degenerate, the proof of Proposition \ref{non-elementary-monodromic} gives an invariant line.
Applying Theorem \ref{bamon} again yields a contradiction of the assumed accumulation.
Thus, both cases have been excluded.
\end{proof}

\subsubsection{Accumulation to a polycycle with several singularities}

We finally consider a polycycle containing at least two singularities.
The following lemma is proved based on Proposition \ref{contact}, by showing that the line through two singularities is already invariant.

\begin{lemma}
\label{several-singularities-polycycle}
Let $\Gamma$ be a polycycle of an algebraic foliation of degree two in $\mathbb{RP}^2$ and suppose that $\Gamma$ contains at least two singularities.
Then either $\mathcal{F}$ has an invariant line or $\Gamma$ is not strictly convex.
\end{lemma}

\begin{proof}
Suppose that $\mathcal{F}$ has no invariant line, and for the sake of contradiction, that $\Gamma$ is strictly convex.
By Proposition \ref{embed-polycycle}, $\Gamma$ is contained in an affine chart $(U,\varphi)$.
By Proposition \ref{cases}, $\mathcal{F}$ can be given by a polynomial vector field $\mathcal X$ of the form \eqref{gvf}.
Let $\mathsf p$ and $\mathsf q$ be two singularities on $\Gamma$, and let $L$ be the projective line through them.
In the affine plane corresponding to $U$, the line segment $(\varphi(\mathsf p),\varphi(\mathsf q))$ is contained in the closed convex region bounded by $\varphi(\Gamma)$, by Proposition \ref{embed-polycycle} and a suitable central projection.

By assumption, the line $L$ is not invariant.
Thus, the two singularities $\mathsf p$ and $\mathsf q$ already account for the two tangency points allowed in Proposition \ref{contact}.
Hence, the vector field $\mathcal X$ crosses the segment $(\varphi(\mathsf p),\varphi(\mathsf q))$ with a constant transversal orientation.
Since $\mathcal X$ has the form \eqref{gvf} and vanishes at $\varphi(\mathsf p)$ and $\varphi(\mathsf q)$, it also crosses $\varphi(L\cap U)\setminus(\varphi(\mathsf p),\varphi(\mathsf q))$ with the same orientation.
However, the existence of a Poincar\'e return map requires $\mathcal X$ to cross $L$ with opposite orientations near $\mathsf p$ and $\mathsf q$ on the side where the return map is defined, as in the proof of \cite[Theorem~2.8]{CS1983}.
This gives a contradiction whether the Poincar\'e return map is defined on the interior or the exterior side of $\Gamma$.
Therefore, $\Gamma$ cannot be strictly convex.
\end{proof}

We now complete the proof of Theorem \ref{main2}.

\begin{proof}[Proof of Theorem \ref{main2}]
By Theorem \ref{even}, all limit cycles of a degree-two foliation are affine cycles.
If the foliation has an invariant line, Theorem \ref{bamon} gives the finiteness.
Therefore, suppose that the foliation has no invariant line and has infinitely many limit cycles.
By Theorem \ref{prop_convex}, these cycles are strictly convex, and Theorem \ref{limit set of limit cycles} gives a subsequence converging to a polycycle $\Gamma$.

Proposition \ref{elementary-cycle-nonaccumulation} and Corollary \ref{non-elementary-not-limit} exclude a singularity or a cycle as the limit.
Lemma \ref{one-singularity-polycycle} excludes a polycycle with one singularity and at least one leaf.
Finally, suppose that the limiting polycycle has at least two singularities.
Lemma \ref{several-singularities-polycycle} shows that either the foliation has an invariant line or the polycycle is not strictly convex.
The first case is covered by Theorem \ref{bamon}.
In the second case, Theorem \ref{convex_polycycle} produces a leaf of the polycycle $\Gamma$ contained in an invariant line, again reducing the problem to Theorem \ref{bamon}.
Every possible limit has therefore been excluded. This contradicts the existence of infinitely many limit cycles.
\end{proof}

Finally, we prove Theorem \ref{main3}, which consists of the following Lemma \ref{six:result:focus} and Proposition \ref{six:result:six}.

By Proposition \ref{cases}, consider a family of quadratic algebraic foliation $\mathcal{F}_\mu$ in $\mathbb{RP}^2$ which is determined by the following system in $U_Z$
\begin{equation}\label{quadratic_like}
\begin{aligned}
\dot x&=\lambda x-y-2x^2+xy+y^2+x(ax^2+bxy+cy^2),\\
\dot y&=x+\lambda y+sx^2-2xy+y(ax^2+bxy+cy^2).
\end{aligned}
\end{equation}
Let $\mathcal{X}_\mu$ denote the vector field given by system \eqref{quadratic_like}, where $\mu=(\lambda,a,b,c,s)\in\mathbb R^5$.
For $\lambda=0$, the Poincar\'{e}--Lyapunov method yields
a formal power series
\[
H(x,y)=\frac{1}{2}x^{2}+\frac{1}{2}y^{2}+\sum_{n=3}^{\infty}H_n(x,y),
\]
where $H_n$ is a homogeneous polynomial of degree $n$, such that
\[
\mathcal{X}_{\mu}H(x,y)
=\sum_{j=0}^{\infty}V_j(x^{2}+y^{2})^{j+1}.
\]
The coefficients $V_j$ are called \emph{Lyapunov constants} (or \emph{focus value}) of vector field $X_\mu$, see
\cite[Chapter~4]{DLA2006} for detail. It is easy to see that $V_0=\lambda$. Applying the classical Poincar\'e--Lyapunov procedure
\cite[Chapter~4]{DLA2006}, we give the first five Lyapunov constants $V_i$ of vector field $X_\mu$ under certain conditions by direct algebraic computation, where $i=1,2,3,4,5$.

\begin{lemma}[Lyapunov constants of vector field $X_\mu$]\label{six:result:focus}
For $V_0=\lambda=0$, the first five Lyapunov constants of \eqref{quadratic_like} satisfy the following formulas.
\begin{itemize}
\item[(1)] The first Lyapunov constant $V_1$ is
$$
V_1=a+c+\frac{6s-1}{4}.
$$
\item[(2)] When $V_1=0$, The second Lyapunov constant $V_2$ is
$$
V_2=\frac{3s-1}{6}b+\frac{5s^2+2s-1}{6}a+\frac{11-55s+38s^2+120s^3}{96}.
$$
Moreover, if $s\ne\frac{1}{3}$, the conditions $V_1=V_2=0$ are equivalent to
$$
c=\frac{1-6s}{4}-a,\quad b=-\frac{5s^2+2s-1}{3s-1}a-\frac{11-55s+38s^2+120s^3}{16(3s-1)}.
$$
\item[(3)] When $V_1=V_2=0$, if $s\ne\frac{1}{3}$ and $D(s)\ne0$, then the third Lyapunov constant $V_3$ is
$$
V_3=-\frac{K(s)D(s)}{24(3s-1)^2}
\bigl(a-a_0(s)\bigr)\bigl(a-a_*(s)\bigr),
$$
where
$$
D(s)=1-8s+5s^2,\quad K(s)=1-2s-3s^2-2s^3
$$
and
$$
a_0(s)=\frac{1-12s}{8},\quad a_*(s)=\frac{37-215s+482s^2-330s^3}{16D(s)}.
$$
\item[(4)] When $V_3=0$ along the curve $a=a_*(s)$, and $s\ne\frac{1}{3}$ and $D(s)\ne0$, the conditions $V_1=V_2=0$ are equivalent to
$$
\begin{aligned}
b&=b_*(s)=\frac{175s^4-105s^3+5s^2+34s-13}{8D(s)},\\
c&=c_*(s)=\frac{3(70s^3-90s^2+53s-11)}{16D(s)}.
\end{aligned}
$$
When $(a,b,c)=(a_*(s),b_*(s),c_*(s))$, the fourth Lyapunov constant $V_4$ is
$$
V_4=-\frac{245K(s)U(s)W(s)}{8192(D(s))^3}Q(s)
$$
where
$$
U(s)=1-2s+2s^2, W(s)=2-2s+3s^2,
$$
and
$$
Q(s)=150s^6-455s^5+459s^4-623s^3+435s^2-76s+2.
$$
\item[(5)] When $(a,b,c)=(a_*(s),b_*(s),c_*(s))$ and $V_4=0$ along the condition $Q(s)=0$, then the fifth Lyapunov constant $V_5$ is
$$
V_5=\frac{245K(s)U(s)W(s)}{1179648D(s)^4}N(s),
$$
where
$$
\begin{aligned}
N(s)=&813000s^{10}-2362525s^9+10201865s^8-39839268s^7\\
&+67997760s^6-78664973s^5+82816473s^4-49428298s^3\\
&+11747410s^2-1048836s+24808.
\end{aligned}
$$
\end{itemize}
\end{lemma}

It is easily checked that $Q(s)=0$ has uniquely a simple root (denoted by $\sigma$) in interval $(2.26,2.27)$, and $D(\sigma)\neq 0$. Let
$$
\mu_0=(0,a_*(\sigma),b_*(\sigma),c_*(\sigma),\sigma),
$$
then it is directly calculated that the first non-zero Lyapunov constant of $\mathcal{X}_{\mu_0}$ is $V_5$. Theorem \ref{main3} follows directly from the following result.
\begin{proposition}\label{six:result:six}
There are parameters $\mu$ arbitrarily close to $\mu_0$ for which the algebraic foliation
$\mathcal F_\mu$ has at least six distinct limit cycles in $\mathbb{RP}^2$.
\end{proposition}

\begin{proof}
We first use a degenerate Hopf bifurcation to produce five limit cycles, and then apply the Poincar\'e-Bendixson  annulus theorem to establish the existence of the remaining one.

By Lemma~\ref{six:result:focus}, the Lyapunov constants of $\mathcal{X}_{\mu_0}$ at the origin satisfies
\[
V_1=V_2=V_3=V_4=0,\qquad -30525<V_5<-30524,
\]
which means the origin is a weak focus of order five. By Lemma~\ref{six:result:focus}, It is straightforward to verify that
\begin{equation}
\det\frac{\partial(V_0,V_1,V_2,V_3,V_4)}{\partial(\lambda,a,b,c,s)}|_{\mu=\mu_0}
=\frac{8575K^2U^2WQ'}{18874368D^3}|_{s=\sigma}\ne0.
\end{equation}
The weak-focus unfolding criterion
\cite{YangYu2017}
therefore gives five hyperbolic limit cycles in any prescribed neighborhood of $O$, for suitable parameters arbitrarily close to $\mu_0$.

On the other hand, by transformation $(x,y)=(\frac{1}{u},\frac{v}{u})$, in the chart $U_X$, with coordinates given by
$\varphi_X(u,v)=[1:v:u]$, the foliation $\mathcal F_\mu$
is represented by
\begin{equation}\label{six:eq:otherchart}
\begin{aligned}
\dot u&=R(u,v)=(\lambda-v)u^2+(v^2+v-2)u+a+bv+cv^2,\\
\dot v&=S(u,v)=-u(1+v^2)-s+v^2+v^3,
\end{aligned}
\end{equation}
which correspond the vector field $\mathcal Z_\mu$. 
Define the compact elliptical disk
$$
\mathcal D=\{(u,v)\in\mathbb{R}^2|\ E(u,v)=(12u-14v-3)^2+225(2v-3)^2-8100\leq0\}.
$$
The boundary $\partial\mathcal D$ admits the parametrization
$$
u(t)=2+\frac{15(1-t^2)+14t}{2(1+t^2)},\quad
v(t)=\frac32+\frac{6t}{1+t^2}.
$$
except one point which corresponds the limit $(u(t),v(t))$ when $t\to\infty$.
When $(\lambda,a,b,c)=(0,a_*(s),b_*(s),c_*(s))$, we have
\begin{equation}\label{six:eq:boundary}
\langle (R(u(t),v(t)),S(u(t),v(t))),(\frac{\partial E}{\partial u}(u(t),v(t)),\frac{\partial E}{\partial v}(u(t),v(t)))\rangle=-\frac{45\sum_{j=0}^6 n_j(s)t^j}{4D(s)(1+t^2)^3},
\end{equation}
where $R$ and $S$ are the components of \eqref{six:eq:otherchart} on this branch, and
\begin{align*}
n_0=&-6300s^4+950s^3+88838s^2-140313s+18059,\\
n_1=&-25200s^4-25440s^3+500360s^2-758984s+99088,\\
n_2=&-6300s^4-89770s^3+568718s^2-790185s+104963,\\
n_3=&9600s^3-6960s^2-11520s+1680,\\
n_4=&6300s^4+89770s^3-486518s^2+658665s-88523,\\
n_5=&25200s^4+35040s^3-507320s^2+747464s-97408,\\
n_6=&6300s^4-950s^3-82238s^2+129753s-16739.
\end{align*}
Let
$$
\frac{\sum_{j=0}^6 n_j(s)t^j}{50000}=z_t^TM(s)z_t,
$$
where
$f_j(s)=\frac{n_j(s)}{50000},\quad z_t=(1,t,t^2,t^3)^T$, and
\begin{equation*}
M=\begin{pmatrix}
f_0&f_1/2&-1/9&-1/45\\
f_1/2&f_2+2/9&f_3/2+1/45&1/10\\
-1/9&f_3/2+1/45&f_4-1/5&f_5/2\\
-1/45&1/10&f_5/2&f_6
\end{pmatrix}.
\end{equation*}
By exact rational arithmetic in computer, when $s$ is near $\sigma$, the leading principal minors $\Delta_k$ ($k=1,2,3,4$) of $M(s)$ is positive, which means $M(\sigma)$ is positive definite by Sylvester's criterion. Hence $\frac{\sum_{j=0}^6 n_j(\sigma)t^j}{50000}>0$ for every real $t$.
The omitted point on the boundary $\mathcal D$ is obtained by $t\to\infty$, where the sign is determined by $f_6(\sigma)>0$.
Together with $D(\sigma)>0$ and \eqref{six:eq:boundary}, this proves the flow of the system \eqref{six:eq:otherchart} crosses transversally into $\mathcal D$ in forward time.

By exact real root isolation algorithm, system \eqref{six:eq:otherchart} has the unique real equilibrium $p=(u_p,v_p)$ in $\operatorname{int}\mathcal D$, where $u_p\in(-2.13168,-2.13154)$ and $v_p\in(0.1968,0.1969)$. A direct computation gives
$$
\operatorname{tr}D\mathcal{Z}_{\mu_0}(p)>0,\quad \det D\mathcal Z_{\mu_0}(p)>0,\quad (\operatorname{tr}D\mathcal Z_{\mu_0}(p))^2-4\det D\mathcal{Z}_{\mu_0}(p)<0
$$
Thus $p\in\operatorname{int}\mathcal D$ is a hyperbolic repelling focus.

Choose a sufficiently small ellipse around $p$ on which the vector field points strictly outward.
Let $B\subset\operatorname{int}\mathcal D$ be the closed disk bounded by this ellipse.
Then $\mathcal A=\mathcal D\setminus\operatorname{int}B$ is a compact positively invariant annulus containing no equilibrium.
Fix these two boundary curves. Their strict transversality and the absence of equilibria in the annulus persist for all $\mu$ sufficiently close to $\mu_0$.
The Poincar\'e annulus theorem then gives a limit cycle $\gamma$ in $\operatorname{int}\mathcal A$. Note that $\widehat{\mathcal D}=\varphi_X(\mathcal D)$ is compact in
$\mathbb{RP}^2$ and does not contain $O$ which lies in $U_Z$. Hence $\gamma$ is distinct from the five limit cycles previously bifurcating from the singularity $O$.

\end{proof}

\appendix
\section{Proof of Proposition \ref{cases} and Corollary \ref{non-generic}}\label{proof2}
\begin{proof}[Proof of Proposition \ref{cases}]
    After a change of coordinates, we may assume that $L=\{[X:Y:Z]\in\mathbb{RP}^2:Z=0\}$.
    We analyze the foliation $\mathcal{F}$ in the affine chart $U_Z$, and then the line at infinity $L_\infty$ of $U_Z$ is $L$.
    By \eqref{affine-restriction}, the foliation on $U_Z$ is given by the vector field
    \begin{equation}
    \label{dehomogeneous-vector-field}
    B(x,y,1)\frac{\p}{\p x}
    -A(x,y,1)\frac{\p}{\p y}.
    \end{equation}
    Dehomogenizing the projective condition \eqref{projective-condition} gives
    \begin{equation}
    \label{dehomogeneous-projective-condition}
    xA(x,y,1) + yB(x,y,1) + C(x,y,1) = 0.
    \end{equation}

    Let $A_{d+1}(x,y)=A(x,y,0),\ B_{d+1}(x,y)=B(x,y,0)$, then $A_{d+1}$, $B_{d+1}$ are the homogeneous components of degree $d+1$ of the polynomials $A(x,y,1),\ B(x,y,1)$, respectively.
    The homogeneous component of degree $d+2$ in \eqref{dehomogeneous-projective-condition} is
    $$x A_{d+1}(x,y) + y B_{d+1}(x,y) = 0.$$
    Since $x$ and $y$ are coprime, there is a homogeneous polynomial $r(x,y)$ of degree $d$ such that
    \begin{equation}
    \label{top-degree-components}
    A_{d+1}(x,y)=-yr(x,y),
    \qquad
    B_{d+1}(x,y)=xr(x,y).
    \end{equation}
    Define
    \begin{equation}
    \label{definition-PQ}
    p(x,y)= B(x,y,1)-xr(x,y),
    \qquad
    q(x,y)=- A(x,y,1)-yr(x,y).
    \end{equation}
    Then $\max\{\deg p,\deg q\}\leq d$, and \eqref{dehomogeneous-vector-field} becomes \eqref{vf_cases}.
    Moreover, taking into account \eqref{dehomogeneous-projective-condition}, if $p+xr$ and $q+yr$ have a common nonconstant factor, then $C(x,y,1)$ has the same factor.
    Homogenizing this factor produces a common factor of $A,B,C$, which leads to a contradiction.
    Thus, $\gcd(p+xr,q+yr)=1$.

    Suppose first that $r\equiv0$.
    Then $\deg p,\deg q\leq d$ by \eqref{top-degree-components} and \eqref{definition-PQ}.
    In fact, we claim that $\max\{\deg p,\deg q\}=d.$
    Otherwise, both $A(x,y,1)$ and $B(x,y,1)$ have degree at most $d-1$, and so, the homogeneous polynomials $A$ and $B$ can be divisible by $Z^2$.
    It then follows from \eqref{projective-condition} that $C$ is divisible by $Z$.
    This makes $Z$ a common factor of $A, B, C$, a contradiction.

    It remains to determine when $L_\infty$ is invariant.
    From
    \eqref{top-degree-components},
    $$A(X,Y,0)=-Yr(X,Y),
    \qquad
    B(X,Y,0)=Xr(X,Y).$$
    Therefore,
    \begin{equation}
    \label{line-at-infinity-invariance}
    \left.(\alpha\wedge dZ)\right|_{Z=0}
    =r(X,Y)(-Y\,dX+X\,dY)\wedge dZ.
    \end{equation}
    Thus, $\alpha\wedge dZ=0$ on $L_\infty$ if and only if $r\equiv0$. Hence $r\equiv0$ gives case (i), whereas $r\not\equiv0$ gives case (ii).
\end{proof}

\begin{proof}[Proof of Corollary \ref{non-generic}]
    Choose coordinates in which $L=\{[X:Y:Z]\in\mathbb{RP}^2:Z=0\}$.
    Note that changes of coordinates preserve the degree of the foliation and the invariance of a line.
    Since $L$ is invariant, equation \eqref{line-at-infinity-invariance} implies that $r\equiv0$.
    Thus, by Proposition \ref{cases}, $\mathcal F$ on the affine chart $\mathbb{RP}^2\setminus L$ is given by a vector field of the non-generic case.
\end{proof}

\begin{proof}[Proof of Proposition \ref{contact}]
    Choose a line at infinity $L_\infty\neq L$ and coordinates in which $L_\infty=\{[X:Y:Z]\in\mathbb{RP}^2:Z=0\}$.
    By Proposition~\ref{cases}, the foliation $\mathcal{F}$ is given by the vector field \eqref{vf_cases} in the affine chart $U_Z$, given by $\{[X:Y:Z]\in\mathbb{RP}^2:Z=1\}$.
    After an affine change of coordinates within this chart, we may assume that the restriction of $L$ has equation $-lx+y-n=0$.
    Note that an appropriate rotation of the affine coordinates can avoid $x=0$.
    A tangency point of $\mathcal F$ with $L$ in this chart satisfies the following two conditions:
        \begin{itemize}
        \item The point lies on the line $L$:
        \begin{equation}\label{eq-1}
        -lx+y-n=0;
        \end{equation}
        \item the normal vector of $L$, namely $(-l,1)$, is orthogonal to the vector field \eqref{vf_cases}, that is,
        \begin{equation}\label{eq-2}
        -l(p(x,y)+xr(x,y))+(q(x,y)+yr(x,y))=0.
        \end{equation}
        \end{itemize}
        Substituting equation \eqref{eq-1} into equation \eqref{eq-2}, we obtain that
        \begin{equation}
            \label{num}
            F(x):=-lp(x,lx+n)+q(x,lx+n)+nr(x,lx+n)=0.
        \end{equation}
        Note that $p(x,lx+n)$, $q(x,lx+n)$, and $r(x,lx+n)$ are polynomials of $x$ with degree at most $d$.
        Moreover, $F\not\equiv0$.
        Otherwise, the polynomial $-l(p(x,y)+xr(x,y))+(q(x,y)+yr(x,y))$ has the factor $-lx+y-n$, and $L$ is invariant under $\mathcal{F}$.
        Consequently, the equation \eqref{num} has at most $d$ roots.
        Therefore, there are at most $d$ tangency points in the affine part of $L$.

        It remains to analyze the point $\mathsf p=L\cap L_\infty$.
        We choose the local coordinate chart  $U_X$, given by $\{[X:Y:Z]\in\mathbb{RP}^2:X=1\}$. In $U_Z\cap U_X$, the coordinate transformation from the coordinates $(x,y)$ in $U_Z$ to the coordinates $(v,u)$ in $U_X$ is
        \begin{equation}
        \label{cotrans_foliation}
        \varphi_{ZX}(x,y)=(v,u)
        =\left(\frac{1}{x},\frac{y}{x}\right),
        \qquad x\neq0.
        \end{equation}
        Under $\varphi_{ZX}$, the vector field \eqref{vf_cases} is represented on $v\neq0$, after multiplication by the factor $v^{d-1}$, by
        \begin{align}\label{vf-uv}
        &(-v^{d+1}p(v^{-1},v^{-1}u)-v^dr(v^{-1},v^{-1}u))\frac{\p}{\p v}\\\notag
        &\qquad +(v^dq(v^{-1},v^{-1}u)-v^dup(v^{-1},v^{-1}u))\frac{\p}{\p u}.
        \end{align}
        Since $\deg p,\deg q\leq d$ and $r$ is homogeneous of degree $d$, the coefficients of the vector field \eqref{vf-uv} extend polynomially to $v=0$.
        Moreover, in the chart $\{[X:Y:Z]\in\mathbb{RP}^2:X=1\}$, the line $L$ is given by $-l-nv+u=0$, and the point $\mathsf p$ corresponds to $(0,l)$.
        Therefore, $\mathsf p$ is a tangency point if and only if the normal vector $(-n,1)$ is orthogonal to the vector field \eqref{vf-uv} at $\mathsf p$, that is,
        \begin{align*}
            &-n(-v^{d+1}p(v^{-1},v^{-1}u)-v^dr(v^{-1},v^{-1}u))|_{(v,u)=(0,l)}\\&+(v^dq(v^{-1},v^{-1}u)-v^dup(v^{-1},v^{-1}u))|_{(v,u)=(0,l)}\\
            =&nr(1,l)+q_d(1,l)-lp_d(1,l)=0,
        \end{align*}
        where $p_d$ and $q_d$ are the homogeneous parts of degree $d$ of $p$ and $q$, respectively.
        Note that this expression is precisely the coefficient of $x^d$ in $F(x)$.
        If there are exactly $d$ tangency points in the affine part of $L$, then $F$ has $d$ distinct roots.
        Since $F\not\equiv0$ and $\deg F\leq d$, its coefficient of $x^d$ is nonzero, and hence $\mathsf p$ is not a tangency point.
        Otherwise, there are at most $d-1$ such points in the affine part of $L$, and including the possible point $\mathsf p$ at infinity gives at most $d$ tangency points in total. This proves the proposition.
\end{proof}

\section{Convexity on $\mathbb{S}^2$}
\label{convexity-s2}

For $p\in\mathbb{S}^2$ and $r>0$, let $\mathcal{B}(p,r)=\{q\in\mathbb{S}^2:d_{\mathbb{S}^2}(p,q)<r\}$ be the open spherical ball with center $p$ and radius $r$.
For a subset $S\subset\mathbb{S}^2$, its interior, closure, and boundary are denoted by $\operatorname{int}(S)$, $\operatorname{cl}(S)$, and $\partial S$, respectively.
We denote the minor arc connecting two non-antipodal points $p,q$ by $[p,q]$.
If $p=-q$, the same notation denotes one of the semicircles specified by the context.
We also set
$$
(p,q]=[p,q]\setminus\{p\},\qquad
[p,q)=[p,q]\setminus\{q\},\qquad
(p,q)=[p,q]\setminus\{p,q\}.
$$

\begin{lemma}
\label{interior}
Suppose that $K$ is a closed convex subset of $\mathbb{S}^2$.
If $z\in\operatorname{int}(K)$ and $x\in\partial K$, then $(x,z]\subset\operatorname{int}(K).$
\end{lemma}

\begin{proof}
Since $z$ is an interior point, there exists $\delta>0$ such that the open spherical ball $\mathcal{B}(z,\delta)\subset K.$
Note that $x\neq-z$.
Otherwise, if $x$ were the antipode of $z$, then any point $v \in \mathbb{S}^2\setminus\{z\}$ lies on a minor arc connecting $x$ to some $p \in \mathcal{B}(z, \delta) \subset K$.
By the convexity of $K$, this implies $v\in K$, that is, $K = \mathbb{S}^2$, which contradicts $x \in \partial K$.

After decreasing $\delta$ if necessary, we may assume that $-x\notin\mathcal{B}(z,\delta)$.
This ensures that $\|k_1 x + k_2 p\| \neq 0$ for any $p \in \mathcal{B}(z, \delta)$ and $k_1, k_2 \geq 0$ (not both zero), allowing us to uniquely define the following convex region $T$ generated by $x$ and $\mathcal{B}(z, \delta)$:
    $$T = \left\{ q \in \mathbb{S}^2 : q = \frac{k_1 x + k_2 p}{\|k_1 x + k_2 p\|}, k_1 \geq 0, k_2 \geq 0, k_1+k_2 > 0, p \in \mathcal{B}(z, \delta) \right\}.    $$
Since $K$ is convex and $\mathcal{B}(z, \delta)\cup \{x\}\subset K$, this region $T$ must be contained in $K$.
For any point $u \in (x,z]$, we define the following radius $r$ based on the spherical law of sines:
    $$ r = \arcsin\left(\frac{\sin d_{\mathbb{S}^2}(u,x)}{\sin d_{\mathbb{S}^2}(z,x)}\sin\delta\right).$$
It follows that $\mathcal{B}(u, r) \subset T \subset K$, see Figure \ref{fig_r}. This implies that every point $u \in (x,z]$ is an interior point of $K$.
    \begin{figure}[htpb]
        \centering
        \includegraphics[width=0.5\linewidth]{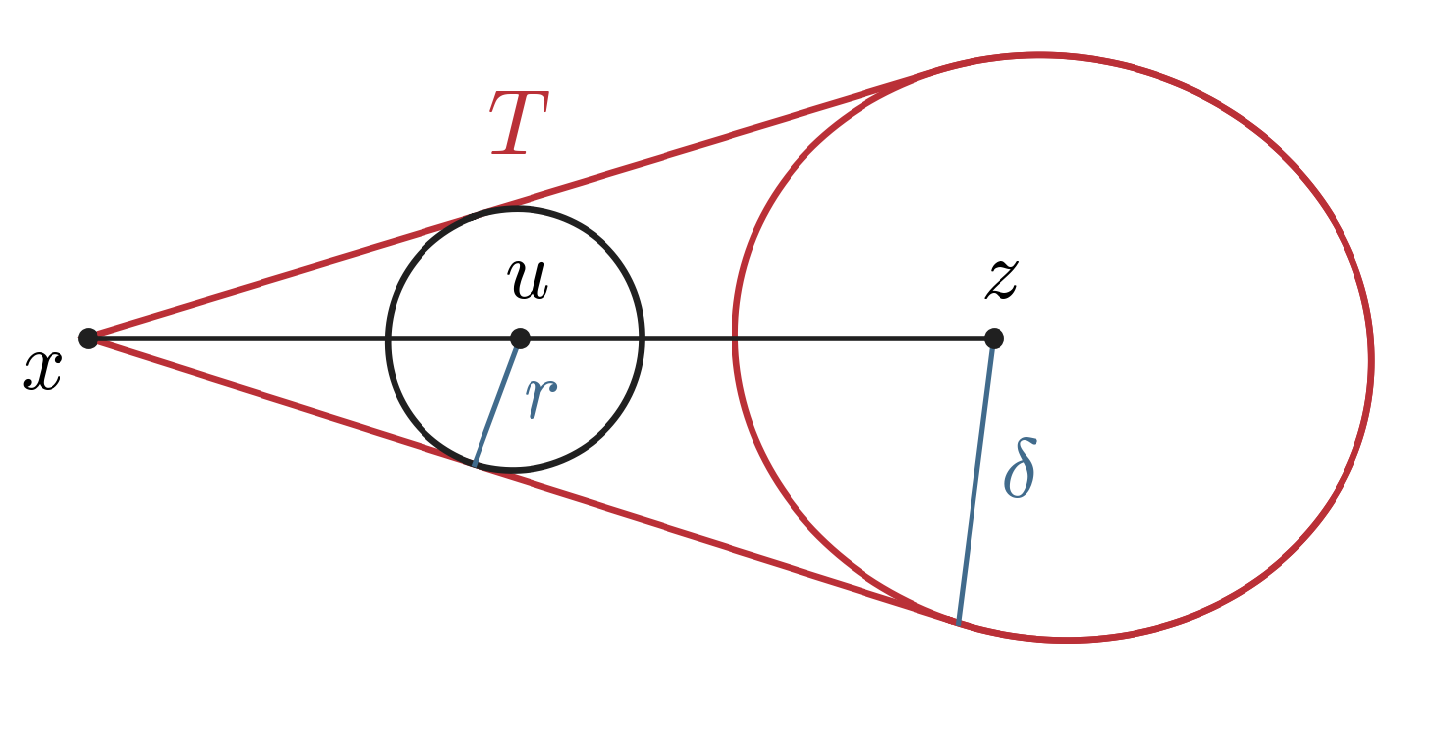}
        \caption{The configuration of the spherical ball $\mathcal{B}(u, r)$.}
        \label{fig_r}
    \end{figure}
\end{proof}

\section{Convergence in the Hausdorff distance}
\label{hausdorff-convergence}

For a point $x\in\mathbb{S}^2$ and a nonempty subset $D\subset\mathbb{S}^2$, let
$d_{\mathbb{S}^2}(x,D) =\inf_{p\in D}d_{\mathbb{S}^2}(x,p). $
For two nonempty compact subsets $A,B\subset\mathbb{S}^2$, their Hausdorff distance is
$d_{\mathbb{S}^2}^{H}(A,B)=\max\left\{\sup_{a\in A}d_{\mathbb{S}^2}(a,B),\sup_{b\in B}d_{\mathbb{S}^2}(b,A)\right\}.$
The following two propositions show that the Hausdorff limit of closed convex sets is convex, and that for such convex subsets of the sphere, the boundaries converge to the boundary of the limit set.

\begin{proposition}
\label{converge}
Suppose that $\{K_n\}$ is a sequence of nonempty closed convex subsets of $\mathbb{S}^2$ that converges to a nonempty closed set $K$ in the Hausdorff distance.
Then $K$ is convex.
\end{proposition}

\begin{proof}
Let $p,q\in K$. It suffices to show that $K$ contains the minor arc $[p,q]$, or one of the semicircles connecting them if they are antipodal.
Since $d_{\mathbb{S}^2}^{H}(K_n,K)\rightarrow0,$ there exist sequences $p_n,q_n\in K_n$ such that
$p_n\rightarrow p \ \text{and}\ q_n\rightarrow q. $

First suppose that $p$ and $q$ are not antipodal, that is, $d_{\mathbb{S}^2}(p, q) < \pi$. If $p=q$, there is nothing to prove. Otherwise, let $r=d_{\mathbb{S}^2}(p,q)$ and $r_n=d_{\mathbb{S}^2}(p_n,q_n).$
Then $0<r<\pi$, and for all sufficiently large $n$, $p_n$ and $q_n$ are not  antipodal.
Since $K_n$ is convex, $[p_n,q_n]\subset K_n$. The minor arc $[p,q]$ can be parameterized by the spherical linear interpolation $\gamma:[0,1]\to\mathbb{S}^2$, where
\begin{equation}
\label{gamma}
\gamma(s)
=\frac{\sin((1-s)r)}{\sin r}\,p
+\frac{\sin(sr)}{\sin r}\,q,
\end{equation}
and $[p_n,q_n]$ can be parameterized similarly by $\gamma_n$, using $p_n,q_n$, and $r_n$.
For any $x\in[p,q]$, choose $s_0\in[0,1]$ such that $x=\gamma(s_0)$, and let $x_n=\gamma_n(s_0)\in[p_n,q_n]$.
Since $p_n\to p$, $q_n\to q$, and $r_n\to r\in(0,\pi)$, formula \eqref{gamma} gives $x_n\to x$ as $n\to+\infty$.
Moreover, by the triangle inequality,
$$ d_{\mathbb{S}^2}(x, K) \le d_{\mathbb{S}^2}(x, x_n) + d_{\mathbb{S}^2}(x_n, K) \le d_{\mathbb{S}^2}(x, x_n) + d_{\mathbb{S}^2}^H(K_n, K). $$
Letting $n\to\infty$ gives $d_{\mathbb{S}^2}(x,K)=0$. Since $K$ is closed, $x\in K$. Hence $[p,q]\subset K$.

Now suppose that $p$ and $q$ are antipodal, that is, $d_{\mathbb{S}^2}(p, q) = \pi$. For each $n$, choose a midpoint $m_n$ of $[p_n,q_n]$ that is contained in $K_n$.
Since $\mathbb{S}^2$ is compact, a subsequence, still denoted by $\{m_n\}$, converges to some $m\in\mathbb{S}^2$.
By continuity of the distance function, we have that
    $$ d_{\mathbb{S}^2}(p, m) = \lim_{n\to\infty} d_{\mathbb{S}^2}(p_n, m_n) = \lim_{n\to\infty} \frac{1}{2}d_{\mathbb{S}^2}(p_n, q_n) = \frac{1}{2}d_{\mathbb{S}^2}(p, q)=\frac{\pi}{2}. $$
Similarly, $d_{\mathbb{S}^2}(q, m) = \frac{\pi}{2}$.
For $m_n \in K_n$, the triangle inequality gives $d_{\mathbb{S}^2}(m, K) \le d_{\mathbb{S}^2}(m, m_n) + d_{\mathbb{S}^2}^H(K_n, K)$.
Taking the limit gives $d_{\mathbb{S}^2}(m, K) = 0$, implying $m \in K$.
Applying the preceding non-antipodal case to the pairs $(p,m)$ and $(m,q)$ gives
$[p,m]\cup[m,q]\subset K.$
This union is a semicircle connecting $p$ and $q$. Thus, $K$ contains a semicircle connecting $p$ and $q$. Consequently, $K$ is convex.
\end{proof}

\begin{proposition}
\label{boundary}
Suppose that $\{K_n\}$ is a sequence of nonempty closed convex proper subsets of $\mathbb{S}^2$ that converges in the Hausdorff distance to a nonempty closed convex proper subset $K$. Then
$d_{\mathbb{S}^2}^{H}(\partial K_n,\partial K)\rightarrow0 \ \text{as }n\longrightarrow\infty.$
\end{proposition}

\begin{proof}
For $\epsilon>0$ and $D\subset\mathbb{S}^2$, let $U_\epsilon(D) =\{x\in\mathbb{S}^2: d_{\mathbb{S}^2}(x,D)<\epsilon\}$ denote the open $\epsilon$-neighborhood of $D$.
Recall that convergence of $\{K_n\}$ to $ K$ in the Hausdorff distance means that for any $\epsilon > 0$, there exists an integer $N$ such that for all $n \ge N$, $K_n \subset U_\e(K)$ and $K \subset U_\e(K_n)$. To prove that the boundaries converge in the Hausdorff distance, we need to show that for any given $\epsilon > 0$ and sufficiently large $n$, $\partial K_n \subset U_\epsilon (\partial K)$ and $\partial K \subset U_\epsilon (\partial K_n)$.

\smallskip
\noindent{(I) We first show that $\partial K_n\subset U_\epsilon(\partial K)$ for all sufficiently large $n$.}

Suppose, for the sake of contradiction, that there exist $\epsilon_0>0$, a subsequence $\{K_{n_j}\}$, and points $x_{n_j}\in\partial K_{n_j}$ such that $d_{\mathbb{S}^2}(x_{n_j},\partial K)\geq\epsilon_0.$
By compactness of $\mathbb{S}^2$, there exists a convergent subsequence (which we still denote by $x_{n_j}$) such that $x_{n_j}\to x$. Since $x_{n_j}\in K_{n_j}$ and $K_{n_j}\to K$, we have $x\in K$.
Moreover, by the triangle inequality,
$$d_{\mathbb{S}^2}(x,\partial K)\geq d_{\mathbb{S}^2}(x_{n_j},\partial K)
-d_{\mathbb{S}^2}(x,x_{n_j})
>\frac{\epsilon_0}{2}$$
for all sufficiently large $j$. Thus, $x\in\operatorname{int}(K)$.

Choose $\epsilon_1>0$ such that $\mathcal{B}(x,\epsilon_1)\subset\operatorname{int}(K)$.
We can choose three points $y_1,y_2,y_3\in\mathcal{B}(x,\epsilon_1)$ such that $x$ lies in the
interior of their spherical convex hull $T$ formed by them, that is,
$$T=\left\{\frac{\sum_{i=1}^{3}k_iy_i}{\left\|\sum_{i=1}^{3}k_iy_i\right\|}:k_i\geq0,\ \sum_{i=1}^{3}k_i=1\right\}.$$
Equivalently, the three points may be chosen so that the matrix $Y=(y_1,y_2,y_3)$ is invertible and $x=\sum_{i=1}^{3}c_iy_i$ for some $c_i>0$.
Since $K_{n_j}\to K$, for each $i$ there are points $y_i^{n_j}\in K_{n_j}$ such that $y_i^{n_j}\to y_i$. Let $Y_{n_j}=(y_1^{n_j},y_2^{n_j},y_3^{n_j}).$
For all sufficiently large $j$, $Y_{n_j}$ is invertible, and
$$c^{n_j}=(c_1^{n_j},c_2^{n_j},c_3^{n_j}):=Y_{n_j}^{-1}x\longrightarrow Y^{-1}x=c.$$
Since $c_i>0$ for $i=1,2,3$, $c_i^{n_j}$ are
also positive for all sufficiently large $j$.
Hence $x$ lies in the interior of the spherical convex hull $T_{n_j}$ generated by
$y_1^{n_j},y_2^{n_j},y_3^{n_j}$.
By convexity, $T_{n_j}\subset K_{n_j}.$ And so, $x \in \operatorname{int}(T_{n_j}) \subset \operatorname{int}(K_{n_j})$.

Since $x \in \operatorname{int}(T)$, the distance $2r = d_{\mathbb{S}^2}(x, \partial T) > 0$ is a constant independent of $n_j$.
Because the convergence $y_i^{n_j} \to y_i$ occurs as $n_j \to +\infty$, we can choose a sufficiently large $N_1$ such that for $n_j>N_1$, $\partial T_{n_j}$ lies entirely within the $r$-neighborhood of $\partial T$, which guarantees $\mathcal{B}(x, r) \subset \operatorname{int}(T_{n_j})$.
Similarly, since $x_{n_j} \to x$ as $n_j \to +\infty$, we can choose a sufficiently large $N_2$ such that for $n_j>N_2$, $d_{\mathbb{S}^2}(x_{n_j}, x) < r$, which ensures $x_{n_j} \in \mathcal{B}(x, r)$.
It follows that $x_{n_j} \in \mathcal{B}(x, r) \subset \operatorname{int}(T_{n_j}) \subset \operatorname{int}(K_{n_j})$ for $n_j>\max\{N_1,N_2\}$.
This contradicts the assumption that $x_{n_j} \in \partial K_{n_j}$.
Therefore, $\partial K_n\subset U_\epsilon(\partial K)$ for all sufficiently large $n$.

\smallskip
\noindent{(II) We next show that
$\partial K\subset U_\epsilon(\partial K_n)$ for all sufficiently large
$n$.}
We distinguish two cases.

Suppose first that $\operatorname{int}(K)=\varnothing$. Then $\partial K=K$.
We claim that $d_{\mathbb{S}^2}^{H}(K_n,\partial K_n)\rightarrow0.$
Otherwise, there exists $\delta>0$, a subsequence $\{K_{n_j}\}$, and points $z_{n_j}\in K_{n_j}$ such that $d_{\mathbb{S}^2}(z_{n_j},\partial K_{n_j})\geq\delta.$
Thus $\mathcal{B}(z_{n_j},\delta)\subset K_{n_j}$.
By compactness, there exists a subsequence, still denoted $\{z_{n_j}\}$, and a point $z$ such that $z_{n_j}\to z$.
Then the Hausdorff convergence implies that $K$ contains $\mathcal{B}(z,\frac\delta2)$, contradicting $\operatorname{int}(K)=\varnothing$.
Thus, for all sufficiently large $n$, $K_n\subset U_{\frac\epsilon2}(\partial K_n)$.
Additionally, since $d^H_{\mathbb{S}^2}(K_n,K) \to 0$, we have $K \subset U_{\frac\e2}(K_n)$ for sufficiently large $n$.
It follows that $\partial K=K\subset U_\epsilon(\partial K_n).$

Now suppose that $\operatorname{int}(K)\neq\varnothing$, and choose
$z\in\operatorname{int}(K)$.
By compactness of $\partial K$, there are points $x^{(1)},\ldots,x^{(m)}\in\partial K$ such that
$\partial K\subset\bigcup_{i=1}^{m}\mathcal{B}\left(x^{(i)},\frac{\epsilon}{2}\right).$
For each $i$, consider the geodesic ray originating from $z$ and passing
through $x^{(i)}$. Choose a point $p^{(i)}$ on $(z,x^{(i)})$, and a point $q^{(i)}$ on the extension of the ray past $x^{(i)}$, such that
$d_{\mathbb{S}^2}(p^{(i)},x^{(i)})<\frac{\epsilon}{2}$ and
$d_{\mathbb{S}^2}(q^{(i)},x^{(i)})<\frac{\epsilon}{2}.$
By Lemma \ref{interior}, $p^{(i)}\in\operatorname{int}(K)$.
Moreover, $q^{(i)} \notin K$, for otherwise $x^{(i)}$ is an interior point of $K$, which leads to a contradiction.
Since $\{p^{(i)}\}_{i=1}^m$ are finitely many points strictly inside $K$, the argument in the second paragraph of part (I) ensures that $p^{(i)}\in\operatorname{int}(K_n)$ for sufficiently large $n$.
Meanwhile, since $\{q^{(i)}\}_{i=1}^m$ are finitely many points strictly outside the closed set $K$, the Hausdorff convergence $d^H_{\mathbb{S}^2}(K_n,K) \to 0$ and the triangle inequality ensure $q^{(i)}\not\in K_n$ for sufficiently large $n$.
Therefore, the segment $[p^{(i)},q^{(i)}]$ intersects $\partial K_n$ at some point $w_n^{(i)}$. Since $w_n^{(i)} \in[p^{(i)}, q^{(i)}]$, we have $d_{\mathbb{S}^2}(x^{(i)},w_n^{(i)})<\frac{\epsilon}{2}$, which means $x^{(i)} \in U_{\frac\e2}(\partial K_n)$.
For any $x\in\partial K$, choose $i$ such that $d_{\mathbb{S}^2}(x,x^{(i)})<\frac\epsilon2$. Then
$$
d_{\mathbb{S}^2}(x,\partial K_n)
\leq d_{\mathbb{S}^2}(x,x^{(i)})
+d_{\mathbb{S}^2}(x^{(i)},w_n^{(i)})
<\epsilon.
$$
Hence $\partial K\subset U_\epsilon(\partial K_n)$ for all sufficiently large $n$.

Combining (I) and (II), we obtain $d_{\mathbb{S}^2}^{H}(\partial K_n,\partial K)\rightarrow0.$
\end{proof}

\section*{Acknowledgements}

Authors are partially supported by National Key R $\&$ D Program of China (No. 2022YFA1005900).
H. He is also partially supported by the National Natural
Science Foundations of China  (No. 12401212) and the Fundamental Research Funds for the Central Universities (Project No. 2025CDJ-IAISYB-013). D. Xiao is also partially supported by the National Natural
Science Foundations of China (Nos. 12271353,  12671200).


\begin{thebibliography}{9}

\bibitem{Bamon1986}
R.~Bamon.
\newblock \emph{Quadratic vector fields in the plane have a finite number of limit cycles.}
\newblock Publications Math\'ematiques de l'IH\'ES, 64:111--142, 1986.

\bibitem{BBI2001}
D.~Burago, Y.~Burago, and S.~Ivanov.
\newblock \emph{A course in metric geometry}.
\newblock Graduate Studies in Mathematics, Vol. 33.
\newblock American Mathematical Society, Providence, RI, 2001.

\bibitem{bochnak2013real}
J.~Bochnak, M.~Coste, and M.-F.~Roy.
\newblock \emph{Real algebraic geometry}.
\newblock Ergebnisse der Mathematik und ihrer Grenzgebiete (3), Vol.~36.
\newblock Springer, 2013.


\bibitem{CamachoLinsNeto1985}
C.~Camacho and A.~Lins Neto,
\newblock \emph{Geometric Theory of Foliations},
\newblock Birkh{\"a}user Boston, 1985.


\bibitem{Cerveau1991}
D.~Cerveau, A.~Lins Neto.
\newblock \emph{Holomorphic foliations in $\mathbb{C}P(2)$ having an invariant algebraic curve.}
\newblock Annales de l'institut Fourier, 41(4): 883--903, 1991.


\bibitem{Chen1979}
L.~Chen and M.~Wang.
\newblock \emph{The relative position and number of limit cycles of a quadratic differential system.}
\newblock Acta Mathematica Sinica, 22(6):751--758, 1979 (in Chinese).

\bibitem{CS1983}
C.~Chicone,  D.~S.~Shafer.
\newblock \emph{Separatrix and limit cycles of quadratic systems and Dulac's theorem.}
\newblock Transactions of the American Mathematical Society, 278(2): 585-612, 1983.


\bibitem{ChristopherLi2007}
C.~Christopher and C.~Li.
\newblock \emph{Limit Cycles of Differential Equations}.
\newblock Advanced Courses in Mathematics -- CRM Barcelona.
\newblock Birkh{\"a}user Basel, 2007.

\bibitem{conway11}
\newblock J.~B.~Conway,
\newblock \textit{Functions of one complex variable}.
\newblock Graduate Texts in Math. 11.
\newblock Springer-Verlag, New York, 1978.

\bibitem{Danzer1963}
L.~Danzer, B.~Gr\"unbaum, V.~Klee.
\newblock \emph{Helly's theorem and its relatives.}
\newblock Proceedings of Symposia in Pure Mathematics, 7: 101--180, 1963.


\bibitem{DesrochesJeffrey2011}
M.~Desroches and M.~R.~Jeffrey.
\newblock \emph{Canards and curvature: the ``smallness of $\varepsilon$''
in slow--fast dynamics.}
\newblock Proc. R. Soc. A, 467(2132):2404--2421, 2011.

\bibitem{Dulac1923}
H.~Dulac.
\newblock \emph{Sur les cycles limites.}
\newblock {Bulletin de la Soci\'{e}t\'{e} Math\'{e}matique de France}, 51:45--188, 1923.

\bibitem{DLA2006}
F.~Dumortier, J.~Llibre, and J.~C. Art{\'e}s.
\newblock \emph{Qualitative theory of planar differential systems}.
\newblock Universitext. Springer-Verlag, Berlin, 2006.



\bibitem{GP1974}
V.~Guillemin, A.~Pollack.
\newblock \emph{Differential topology.}
\newblock Prentice-Hall, 1974.



\bibitem{harnack}
A.~Harnack.
\newblock \emph{Ueber die vieltheiligkeit der ebenen algebraischen curven}.
\newblock Mathematische Annalen, 10:189--198, 1876.

\bibitem{Hilbert1902}
D.~Hilbert.
\newblock \emph{Mathematical problems.}
\newblock {Bulletin of the American Mathematical Society}, 8(10):437--479, 1902.

\bibitem{Ilyashenko1985}
Y.~S.~Ilyashenko.
\newblock \emph{Dulac's memoir ``on limit cycles'' and related problems of the local theory of differential equations.}
\newblock {Russian Mathematical Surveys}, 40(6):1--49, 1985.

\bibitem{Ilyashenko1991}
Y.~S.~Ilyashenko.
\newblock \emph{Finiteness theorems for limit cycles}.
\newblock Translations of Mathematical Monographs, Vol. 94, American Mathematical Society, 1991.

\bibitem{Ilyashenko2002}
Y.~S.~Ilyashenko.
\newblock \emph{Centennial history of Hilbert's 16th problem}.
\newblock {Bulletin of the American Mathematical Society}, 39(3): 301 - 354, 2002.

\bibitem{Ilyashenko2008}
Y.~S.~Ilyashenko, S.~Yakovenko.
\newblock \emph{Lectures on analytic differential equations.}
\newblock Graduate Studies in Mathematics, Vol.~86.
\newblock American Mathematical Society, 2008.

\bibitem{Kaloshin2003}
V.~Kaloshin.
\newblock \emph{The existential Hilbert 16-th problem and an
estimate for cyclicity of elementary polycycles}.
\newblock Inventiones Mathematicae, 151: 451--512, 2003.

\bibitem{LL2016}
I.~E.~Leonard, J.~E.~Lewis.
\newblock \emph{Geometry of convex sets.}
\newblock John Wiley \& Sons, 2015.

\bibitem{Li1986}
C.~Li.
\newblock \emph{Non-existence of limit cycle around a weak focus of order three for any quadratic system.}
\newblock Chinese Annals of Mathematics, Series B, 7(2):174--190, 1986.

\bibitem{Li2003}
J.~Li.
\newblock \emph{Hilbert's 16th problem and bifurcations of planar polynomial vector fields.}
\newblock International Journal of Bifurcation and Chaos, 13(01):47--106, 2003.


\bibitem{LinsNeto1988}
A.~Lins Neto.
\newblock \emph{Algebraic solutions of polynomial differential equations and foliations in dimension two.}
\newblock In X.~Gomez-Mont, J.~Seade, and A.~Verjovsky, editors, \emph{Holomorphic Dynamics}, Vol.~1345, pp.~192--232.
\newblock Springer-Verlag, Berlin, 1988.

\bibitem{Movasati2008}
H.~Movasati, E.~Vieira.
\newblock \emph{Projective limit cycles}.
\newblock Moscow Mathematical Journal, 9(4): 855--866, 2009.

\bibitem{Movasati2022}
H.~Movasati.
\newblock \emph{Introduction to algebraic curves and foliations}.
\newblock Lecture Notes, version of November 10, 2022.


\bibitem{Perko2013}
L. Perko.
\newblock\emph{Differential equations and dynamical systems.}
\newblock  Vol.~7.
\newblock Springer Science \& Business Media, 2013.

\bibitem{Poincare}
 H.~Poincar\'e,
\newblock \emph{M\'emoire sur les courbes d\'efinies par les \'equations diff\'erentielles},
\newblock Journal de Math\'ematiques Pures et Appliqu\'ees, s\'erie 3, 7:375--422, 1881; 8:251--296, 1882.
\newblock Reprinted in \emph{Oeuvres de Henri Poincar\'e},
Vol.~1, pp.~3--84, Gauthier-Villars, Paris, 1951.

\bibitem{Roussarie1998}
R.~Roussarie.
\newblock \emph{Bifurcation of planar vector fields and Hilbert's sixteenth problem}.
\newblock Progress in Mathematics, Vol. 164,
Birkh{\"a}user Basel, 1998.

\bibitem{Shi1980}
S.~Shi.
\newblock \emph{A concrete example of the existence of four limit cycles for plane quadratic systems.}
\newblock Scientia Sinica, 23(2):153--158, 1980.


\bibitem{YangYu2017}
J.~Yang and P.~Yu.
\newblock \emph{Nine limit cycles around a singular point by perturbing a cubic Hamiltonian system with a nilpotent center.}
\newblock Applied Mathematics and Computation, 298: 141--152, 2017.

\bibitem{Ye1986}
Y. Ye.
\newblock \emph{Theory of limit cycles}.
\newblock Translated from the Chinese by Chi Y. Lo,
Translations of Mathematical Monographs, vol. 66,
American Mathematical Society, Providence, RI, 1986.

\bibitem{Yeung2025}
M. Yeung.
\newblock \emph{Dulac's theorem revisited.}
\newblock Qualitative Theory of Dynamical Systems, 24(2): Paper No. 57, 27~pp., 2025.



\bibitem{Zhang2017}
X. Zhang.
\newblock \emph{Integrability of dynamical systems: algebra and analysis}.
\newblock Developments in Mathematics, Vol. 47, Springer, 2017.

\bibitem{Zhang1992}
Z. Zhang, T. Ding, W. Huang, and Z. Dong.
\newblock \emph{Qualitative theory of differential equations}.
\newblock Translations of Mathematical Monographs, Vol. 101, American Mathematical Society, 1992.
\end{thebibliography}
\end{document}